\input{set-pdf-version.tex}

\def\santtu{Santtu Söderholm}

\def\joonas{Joonas Lahtinen}

\def\carsten{Carsten H. Wolters}

\def\sampsa{Sampsa Pursiainen}

\def\mytitle{The Effects of Peeling on Finite Element Method -based EEG Source Reconstruction}

\def\myjournal{Biomedical Signal Processing and Control}

\def\mykeywords{%
	electroencephalography\sep%
	EEG\sep%
	finite element method\sep%
	lead field\sep%
	source localisation\sep%
	sLORETA\sep%
	Dipole Scan\sep%
	H(div)%
}

\def\mypacscodes{ 41.20.Cv \sep 42.30.Wb }

\def\mymsccodes{ 65N30 \sep 92C55 }

\def\mysubject{%
	The problem of reconstructing brain activity from electric potential
	measurements performed on the surface of a human head is not an easy task:
	not just because the solution of the related inverse problem is
	fundamentally ill-posed (not unique), but because the methods utilized in
	constructing a synthetic forward solution themselves contain many
	inaccuracies. One of these is the fact that the usual method of modelling
	primary currents in the human head via dipoles brings about at least \(2\)
	modelling errors: one from the singularity introduced by the dipole, and
	one from placing such dipoles near conductivity discontinuities in the
	active brain layer boundaries.\\\indent
	In this article we observe how the removal of possible source locations
	from the surfaces of active brain layers affects the localisation accuracy
	of two inverse methods, \sLORETA\ and \DipoleScan, at different
	signal-to-noise ratios (\SNR), when the \(\Hdiv\) source model is used. We
	also describe the \fe\ forward solver used to construct the synthetic
	\EEG\ data, that was fed to the inverse methods as input, in addition to
	the meshes that were used as the domains of the forward and inverse
	solvers. Our results suggest that there is a slight general improvement in
	the localisation results, especially at lower noise levels. The applied
	inverse algorithm and brain compartment under observation also affect the
	accuracy.%
}

\documentclass[5p,twocolumn]{elsarticle}

\pdfoutput=1

\usepackage[utf8]{inputenc}
\usepackage[T1]{fontenc}
\usepackage{mathtools} % patched amsmath
\usepackage{fourier}
\usepackage{siunitx}
\usepackage[labelfont=bf]{caption}
\usepackage[table]{xcolor}
\usepackage{booktabs}
\usepackage{pgfmath}
\usepackage{xfp}

\def\zi{Zeffiro Interface}

\def\fe{finite element}
\def\FE{FE}

\def\eeg{electroencephalography}
\def\EEG{EEG}
\def\meg{magnetoencephalography}
\def\MEG{MEG}

\def\PBO{PBO}

\def\MPO{MPO}

\def\sLORETA{sLORETA}
\def\DipoleScan{Dipole Scan}
\def\csf{cerebrospinal fluid}
\def\CSF{CSF}

\def\MRI{MRI}

\def\CT{CT}

\def\SNR{\text{SNR}}

\DeclarePairedDelimiter{\arcs}{(}{)}

\DeclarePairedDelimiter\norm{\lVert}{\rVert}
\DeclarePairedDelimiter\angles{\langle}{\rangle}
\DeclarePairedDelimiter\qunit{[}{]}
\DeclarePairedDelimiter\abs{|}{|}

\RenewDocumentCommand\vec{m}{\mathbf{#1}}

\NewDocumentCommand\subi{m}{_{#1}}

\NewDocumentCommand\subt{m}{_{\mathrm{#1}}}

\NewDocumentCommand\Jp{}{\vec J\subi{\mathrm p}}

\NewDocumentCommand\tetrahedralization{}{\mathcal{T}}
\NewDocumentCommand\modelerror{}{E_{\tetrahedralization + \Jp}}
\NewDocumentCommand\measurementerror{}{E_{M}}
\NewDocumentCommand\inverseerror{}{E_{K}}

\NewDocumentCommand\est{sm}{%
	\IfBooleanTF{#1}{%
		\angles*{#2}%
	}{%
		\angles{#2}%
	}%
}

\NewDocumentCommand\rec{sm}{%
	\IfBooleanTF{#1}{%
		\arcs*{#2}^\wedge%
	}{%
		#2^\wedge%
	}%
}

\NewDocumentCommand\Hdiv{}{H(\text{div})}

\NewDocumentCommand\MAG{}{MAG}

\NewDocumentCommand\RDM{}{RDM}

\NewDocumentCommand\restr{smm}{%
	\IfBooleanTF{#1}{%
		\arcs*{#2\mid#3}%
	}{%
		\arcs{#2\mid#3}%
	}%
}

\NewDocumentCommand\pdist{sm}{%
	\IfBooleanTF{#1}{%
		p\arcs*{#2}%
	}{%
		p\arcs{#2}%
	}%
}

\NewDocumentCommand\inverse{sm}{%
	\IfBooleanTF{#1}{%
		\arcs*{#2}^{-1}%
	}{%
		#2^{-1}%
	}%
}

\NewDocumentCommand\pinverse{sm}{%
	\IfBooleanTF{#1}{%
		\arcs*{#2}^\dagger%
	}{%
		#2^\dagger%
	}%
}

\NewDocumentCommand\transpose{sm}{%
	\IfBooleanTF{#1}{%
		\arcs*{#2}^{\mathrm T}%
	}{%
		#2^{\mathrm T}%
	}%
}

\NewDocumentCommand\regularization{sm}{%
	\IfBooleanTF{#1}{%
		\arcs*{#2}^\ddagger%
	}{%
		#2^\ddagger%
	}%
}

\NewDocumentCommand\onesvec{}{\mathbf{1}}

\NewDocumentCommand\nnorm{smm}{%
	\IfBooleanTF{#1}{%
		\norm*{#2}_{#3}%
	}{%
		\norm{#2}_{#3}%
	}%
}

\DeclareMathOperator*\indmax{ind\,max}

\DeclareMathOperator*\inv{inv}
\DeclareMathOperator*\SD{SD}

\NewDocumentCommand\trace{sm}{%
	\IfBooleanTF{#1}{%
		\operatorname{trace}\arcs*{#2}%
	}{%
		\operatorname{trace}#2%
	}%
}

\NewDocumentCommand\bff{}{W}

\NewDocumentCommand\La{}{L\subt{a}}

\NewDocumentCommand\Ln{}{L\subt{n}}

\NewDocumentCommand\rrv{sm}{%
	\IfBooleanTF{#1}{%
		\operatorname{RRV}\arcs*{#2}%
	}{%
		\operatorname{RRV}#2%
	}%
}

\NewDocumentCommand\gof{sm}{%
	\IfBooleanTF{#1}{%
		\operatorname{GoF}\arcs*{#2}%
	}{%
		\operatorname{GoF}#2%
	}%
}

\NewDocumentCommand\tsvd{}{truncated singular value decomposition}
\NewDocumentCommand\tSVD{}{tSVD}

\NewDocumentCommand\source{smmm}{%
	\IfBooleanTF{#1}{%
		\arcs*{#2, #3, #4}%
	}{%
		\arcs{#2, #3, #4}%
	}%
}

\DeclareMathOperator\mMAG{MAG}

\DeclareMathOperator\mRDM{RDM}

\NewDocumentCommand\domain{}{\Omega}

\NewDocumentCommand\peeld{}{d\subt{p}}

\NewDocumentCommand\figref{m}{Figure~\ref{#1}}

\NewDocumentCommand\tabref{m}{Table~\ref{#1}}

\NewDocumentCommand\colorhmportion{}{50}

\NewDocumentCommand\hmcell{mm}{%
	\def\division{\fpeval{#1 / #2}}%
	\def\totalpercentage{\fpeval{100 * \division}}%
	\ifdim\totalpercentage pt < 50 pt%
		\def\mycolor{blue}%
		\def\colorpercentage{\fpeval{\colorhmportion - \colorhmportion * \division}}%
	\else%
		\def\mycolor{red}%
		\def\colorpercentage{\fpeval{\colorhmportion * \division}}%
	\fi%
	\edef\x{\noexpand\cellcolor{\mycolor!\colorpercentage!white}}\x%
	\num{#1}%
}

\NewDocumentCommand\code{m}{\texttt{#1}}

\makeatletter
\ifx\c@figure\undefined
	
\fi
\ifx\c@section\undefined
	
\fi
\ifx\c@table\undefined
	
\fi
\makeatother

\numberwithin{figure}{section}

\DeclareCaptionType{subfigure}

\DeclareCaptionLabelFormat{subfigure}{\textbf{(#2)}}

\RenewDocumentEnvironment{subfigure}{O{\textwidth}mm}{%
	\begin{minipage}[t]{#1}%
	\captionsetup{type=figure}%
	\centering%
}{%
	\captionof{subfigure}{#2}%
	\label{#3}%
	\end{minipage}%
}

\NewDocumentCommand\subfigref{m}{\textbf{(\ref{#1})}}

\numberwithin{equation}{section}
\counterwithin{table}{section}

\RequirePackage[a-3u,mathxmp]{pdfx}[2021/12/31]
\RequirePackage{url}
\hypersetup{hidelinks}

\counterwithin*{subfigure}{figure}

\journal{\myjournal}

\begin{document}

\begin{frontmatter}


\def\tampere{tampere}

\affiliation[\tampere]{%
	organization={Computing Sciences, Tampere University},%
	addressline={Korkeakoulunkatu 3},%
	city={Tampere},%
	postcode={33014},%
	state={Pirkanmaa},%
	country={Finland}%
}

\def\munsterbiomag{munster-biomag}

\affiliation[\munsterbiomag]{%
	organization={Institute for Biomagnetism and Biosignalanalysis, University of Münster},%
	addressline={Schlossplatz 2},%
	city={Münster},%
	postcode={48149},%
	state={Nordrhein-Westfalen},%
	country={Germany}%
}

\def\munstercreutzfeldt{munster-creutzfeldt}

\affiliation[\munstercreutzfeldt]{%
	organization={Otto Creutzfeldt Center for Cognitive and Behavioral Neuroscience, University of Münster},%
	addressline={Schlossplatz 2},%
	city={Münster},%
	postcode={48149},%
	state={Nordrhein-Westfalen},%
	country={Germany}%
}

\author[\tampere]{\santtu}

\author[\tampere]{\joonas}

\author[\munsterbiomag,\munstercreutzfeldt]{\carsten}

\author[\tampere]{\sampsa}

\title{\mytitle}

\begin{abstract}
\mysubject
\end{abstract}

\begin{keyword}
	\mykeywords\sep
	\PACS\mypacscodes\sep
	\MSC\mymsccodes
\end{keyword}

\end{frontmatter}

\section{Introduction}\label{sec:introduction}

The \eeg\ (\EEG) forward problem of attempting to mathematically construct the
electric potentials \(u\) produced by given electrical activity \(\Jp\) in the
human brain is almost a century-old
endeavour~\cite{berger-1929}\cite{brette--destexhe-2012}, the requirements of
which are fairly well known. For simplified cases, such as homogenous and
unbounded conductors, layered conductors and spherical head models, there
exist analytical or formulaic
solutions~\cite{sarvas-1987}\cite{ary-etal-1981}\cite{de-munck--peters-1993}.
However, realistic head geometries require the application of numerical
approaches, such as the finite difference method~\cite{hallez-2008},
boundary or surface element
method~\cite{hamalainen-etal-1993}\cite{kybic--etal-2005} and finite or volume
element
method~\cite{wolters--etal-2007}\cite{drechsler--etal-2009}\cite{he-rezaei-pursiainen-2018}
to solve the same problem.

The primary currents \(\Jp\) themselves are often modelled as electrical
dipoles. One of the problems that arises in formulating the forward problem in
this way is the appearance of singularities in the potential field \(u\)
generated by the dipoles~\cite{wolters--etal-2007}\cite{drechsler--etal-2009}:
\(u\) is inversely proportional to the distance \(r\) from the dipole position
\(\vec x\), where \(u\) is singular. This has implications on the convergence
of numerical methods, which attempt to build the forward solution based on a
finite element model of the human head, with conductivity jumps between the
different brain compartments. More specifically, in the finite element
formulation~\cite{miinalainen-2019}, the load vector \(\vec f\) is not
well-defined in the case of a singular source. The accuracy of the forward
solution also ends up being reduced, as a dipole is placed near a boundary of
an active
compartment~\cite{wolters--etal-2007}\cite{pursiainen--vorwerk--wolters-2016}\cite{miinalainen-2019}.

To tackle the issue of potential singularities themselves, a so
called~\emph{subtraction
method}~\cite{wolters--etal-2007}\cite{drechsler--etal-2009}\cite{beltrachini-2019}\cite{source-reconstruction-book}
has been recently utilized. It involves splitting the potential field \(u\)
into a sum of two separate potentials, a problematic singularity potential
\(u^\infty\) and a correction potential \(u\subt{corr}\), and solving the
\EEG\ forward problem in the case of \(u\subt{corr}\). This amounts to
removing the singularity from the dipole model.

Another approach to handling the dipolar singularity is the so-called
\(\Hdiv\) model~\cite{pursiainen--vorwerk--wolters-2016}\cite{miinalainen-2019},
which assumes a higher smoothness or regularity at the primary source level,
being neurophysiologically well motivated~\cite{hamalainen-etal-1993}%
\cite{murakami--okada-2015}. Assuming such smoothness, the \(\Hdiv\) approach
resolves the ill-definedness of the load vector \(\vec f\) by replacing the
theoretical dipole with a non-singular function. This is achieved by requiring
that the model of a dipole is square-integrable in the finite element domain
\(\domain\), or belongs to the space \(\Hdiv = H\arcs{\mathrm{div},\domain}\),
whose elements can be constructed as linear combinations of
divergence-conforming basis functions. Here dipolar moments \(\vec d\) are
approximated by a vector, whose orientation is defined by a set of finite
element nodes surrounding a central node, taking into account the \emph{a
priori} information about the primary currents which are normally oriented to
the surface of the gray matter layer. Since the vector is only supported by a
few nodes, it is very focal and is therefore able to fit inside the thin gray
matter layer, if the resolution of the mesh is fine enough.

Both approaches suffer from inaccuracies near conductivity jumps. In the case
of \(\Hdiv\), which is used in this work, a dipole might still be placed at
the very edge of an active brain compartment, with some of the supporting
finite element nodes in a neighboring compartment. Due to the different
conductivities in the different compartments, the volumetric current between
the compartments ends up being altered, resulting in a change in the modelled
potential field \(u\)~\cite{pursiainen--vorwerk--wolters-2016}. For the
subtraction method, its numerical accuracy decreases as a dipole is placed
near a conductivity jump, as the upper bound of its error is a function of the
distance from the conductivity jump~\cite{wolters--etal-2007}.

The main purpose of this article is to find out how limiting the possible
source positions to a given distance from conductivity jumps affects the
localisation accuracy of two inverse methods,
\sLORETA~\cite{pascual-marqui-2002} and
\DipoleScan~\cite{neugebauer-etal-2022}, when a set of synthetic \EEG\
measurements has been given to them as input.  This is all performed in an
adapted high-resolution finite element mesh, with a base resolution of
\SI{2}{\milli\meter}, with additional refinements performed on the active
surfaces, to accommodate the \(\Hdiv\) approach in the
\num{1}--\SI{4}{\milli\meter}-thick gray matter layer. The motivation behind
choosing the two inverse methods lies in the fact that \sLORETA\ has been
suggested to localise \emph{distributed} patch-like sources in the entire head
volume, if the number of sources and the amount of measurement noise remains
low, with \(\SNR > \SI{10}{\decibel}\)~\cite{saha-2015}\cite{carboni-2022}.
\DipoleScan\ has been shown to work well in the case of cortically constrained
single-dipole \EEG\ reconstructions, with only a few millimeters' spatial
deviations between reconstructions and original source dipoles, in the
simulation studies of~\cite{fuchs-1998}.

In this article, Section~\ref{sec:methods} will focus on discussing the
anatomy of the forward problem, how the utilized inverse methods localise
sources based on the forward solution and how the \emph{peeling} or removal of
possible source positions from surfaces of active brain layers is performed at
the start of the forward algorithm. In Section~\ref{sec:results}, we then
present our results, evaluate the performance of the peeling algorithm,
compare the numerical forward solver to an analytical one and finally see how
the peeling of possible source locations affects the localisation accuracy in
a realistic head model. In Section~\ref{sec:discussion} results are discussed,
Section~\ref{sec:summary} summarizes our results and
Section~\ref{sec:future-prospects} presents possible future directions of
study.

\section{Methods and models}\label{sec:methods}

\subsection{Mathematical methods}\label{ssec:methods}

\begin{figure}[h!]
	\centering
	\includegraphics[width=0.8\linewidth]{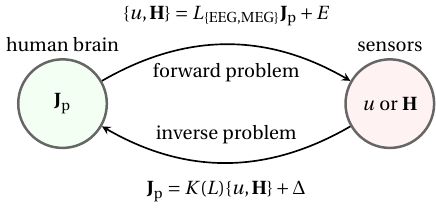}
	\caption{A simplified~\cite{scales--snieder-2001} illustration of the duality
	between forward and inverse problems in \EEG\ and \MEG\ imaging. Here
	\(L\) is the lead field matrix and \(K = K(L)\) an inversion method
	-specific kernel, that provides a reconstruction of the original
	activities \(\Jp\) based on the measurements at the
	sensors.~\cite{atena-2021}}
	\label{fig:forward-vs-inverse}
\end{figure}

For the purpose of testing inverse reconstructions of brain activity, we apply
the Matlab-based software suite \zi~\cite{he-rezaei-pursiainen-2018}, which
builds \(\Hdiv\)-based lead field matrices \(L\) for different brain imaging
modalities, such as \eeg\ (\EEG) or \meg\ (\MEG) in a volumetric domain
\(\domain\), discretized with a finite element mesh\cite{zeffiro-fem-algo}.
Here a forward solver refers to finding a scalar or vector field, generated by
a (synthetic) set of dipole-like sources or the \emph{primary current
distribution} \(\Jp = \Jp\arcs{\vec x}\) at positions \(\vec x\) in the domain
\(\domain\). This is done by mapping \(\Jp\) to an electric or magnetic set of
sensors placed on the surface of \(\domain\), i.e., by multiplying \(\Jp\)
with \(L\)~\cite{hamalainen-etal-1993}, as
in~\cite{he-rezaei-pursiainen-2018}\cite{miinalainen-2019}
\begin{equation}\label{eq:lead-field}
	M = L\Jp + E\,.
\end{equation}
Here \(M\) could consist of electric potentials \(u\) or magnetic field
strenghts \(\vec H\) at the sensors \(S\), and \(E\) is the error arising from
numerical approximation and measurement noise. Once \(L\) is constructed, it
plays a pivotal role in producing the inverted dipoles or primary current
distribution \(\rec\Jp\) in a given
volume~\cite{hauk-2022-neuroimage}\cite{dpm--ga-inverse-1998}. This is
illustrated in~\figref{fig:forward-vs-inverse}.

Regardless of the modality of \(L\) and in the absense of coarse artefacts, the
origins of the \emph{localisation error} \(\Delta\) of
\figref{fig:forward-vs-inverse} in the entire inversion process can be split
roughly into \num 3 components:
\begin{equation}\label{eq:localisation-error-components} \Delta = \modelerror
+ \measurementerror + \inverseerror\,. \end{equation} Here \(\modelerror\) is
the modelling error resulting from the discretization or tetrahedralization
\(\tetrahedralization\) of the head model and the configuration such as
positioning and orientation of the source space
\(\Jp\)~\cite{wang-2001}\cite{lanfer-2012}\cite{akalin--etal-2013}. The symbol
\(\measurementerror\) represents the error brought to the measurements \(u\)
and/or \(\vec H\) by (simulated) measurement
noise~\cite{phillips-2002}\cite{teplan-2002}\cite{meg-eeg-primer-2017} and
\(\inverseerror\) is an inverse method -specific error brought about by, among
other things, the biasing nature of each method, caused by prior assumptions
regarding the measurements \(M\) or the primary current distribution \(\Jp\)
in the derivation of said
methods~\cite{baillet--garnero-1997}\cite{scales--tenorio-2001}\cite{phillips-2002}\cite{calvetti--etal-2009}.

One of the most significant contributors to the term \(\modelerror\) in
\eqref{eq:localisation-error-components} is the actual shape and structure of
the head model used. It has been reported that using spherical models instead
of realistic ones generated by segmenting \MRI\ or \CT\ data increase the
localisation error \(\Delta\) by as much as
\SI{40}{\milli\meter}~\cite{roth-1993}. On the other hand, inaccuracies in the
meshing procedure, such as labeling finite elements into compartments they do
not belong to, or not modelling enough different compartments to accurately
simulate brain structure also contribute to the error. For example, it has
been concluded~\cite{vorwerk-2014}, that the absence of a \csf\ (\CSF) layer
would have detrimental effects on the accuracy of the forward solution \(L\)
and the ensuing inverse reconstruction
\(\rec\Jp\)~\cite{ramon--etal-2006}\cite{wendel--etal-2008}\cite{piastra--etal-2021}.

\begin{figure}[h!]
	\includegraphics[width=\columnwidth]{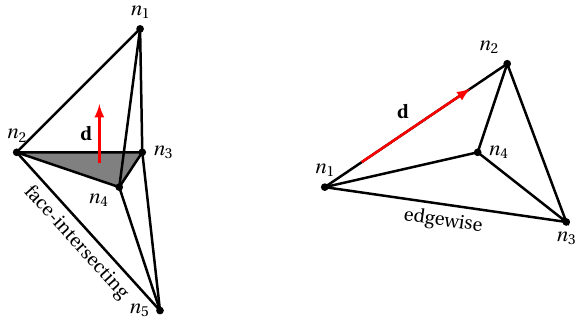}
	\caption{
		Face-intersecting and edgewise dipoles. In the first case, the
		opposing node pair \(\arcs{n_5,n_1}\) function as the ends of a
		dipole, whereas in the second case it is the pair \(\arcs{n_1, n_2}\).
		The red vectors \(\vec d\) denote the directions of the dipoles from
		negative to positive end. The \(\Hdiv\) source model places both kinds
		of sources into the \FE\
		mesh.~\cite{pursiainen--vorwerk--wolters-2016}
	}
	\label{fig:source-models}
\end{figure}

The tetrahedralization \(\tetrahedralization\) might also affect the
generation of the source space \(\Jp\), as it has the possibility of
restricting how synthetic sources can be placed into the head model. In the
case where \(\Jp\) forms a divergence-conforming field, the source space is
directly anchored to the tetrahedra, and interpolated across their faces and
edges, due to the presence of a face-intersecting and edgewise
divergence-conforming\ source model~\cite{pursiainen--vorwerk--wolters-2016}.
This is visualized in \figref{fig:source-models}, where the directions of the
dipoles are given by the vectors \(\vec d\). These align with the dipole
moments \(\vec p\) of the corresponding dipoles.

This holds significance, because as as a direct consequence of Maxwell's
equations the electric and magnetic components of an electromagnetic field are
perpendicular to each other~\cite{dassios--etal-2007}. Hence \EEG\ electrodes
are sensitive to the radial dipoles pointing in their direction, whereas \MEG\
sensors best detect tangential sources, that are perpendicular to the ideal
\EEG\ orientations~\cite{piastra--etal-2021}. Therefore, both \EEG\ and \MEG\
might fail to detect some source orientations, which are themselves subject to
modelling errors.

The possible causes for the measurement noise component \(\measurementerror\)
of the localisation error \(\Delta\) are numerous. These include insufficient
shielding of the imaging room, poor properties of the measuring equipment, the
ill-positioning of the sensors on the scalp, the lack of establishment of
proper electrical contact between electrodes and skin via the application of
electolyte gel, and so forth.~\cite{teplan-2002}\cite{meg-eeg-primer-2017} To
model the uncertainties related to \(\measurementerror\), we add Gau\ss ian
noise to the simulated signal. To understand how the added noise affects the
reconstruction of \(\Jp\) in the inversion phase, we investigate a sample of
inverse estimates obtained with random Gau\ss ian noise realizations.

To finally consider the \(\inverseerror\) component of \(\Delta\) in
\eqref{eq:localisation-error-components}, the ill-posed nature of the inverse
problem needs to be taken into
account~\cite{scales--tenorio-2001}\cite{benning-burger-2018}\cite{clason-2021}.
The problem is underdetermined, characterized by a larger number of variables
or degrees of freedom than there are restrictions or equations in the forward
model. This leads to the need to incorporate \emph{prior-assumptions}
regarding the primary current distribution \(\Jp\) of
equation~\eqref{eq:lead-field} into the model. This can be achieved by
applying \emph{regularization} or \emph{penalty functions} to the related cost
function~\cite{benning-burger-2018}\cite{kaipio--somersalo-2004}. In what
follows, we briefly review two inverse methods, \sLORETA\ and \DipoleScan, that
represent different approaches to favouring certain kinds of unique estimates
\(\rec\Jp\) of \(\Jp\).

\subsection{\sLORETA}\label{ssec:sloreta}

In \emph{Standardized Low Resolution Tomography} or
\sLORETA~\cite{pascual-marqui-2002}, the \emph{Moore--Penrose pseudoinverse}
\(\pinverse L = \inverse*{\transpose L L} \transpose L\) of the typical
least-squares solution is modified to produce an initial regularized solution
\(\rec\Jp\) as follows:
\begin{equation}\label{eq:sloreta-regularized-solution}
	\rec\Jp
	=
	\transpose L \pinverse*{ L \transpose L + \lambda I} L \Jp
	=
	L^\ddagger L \Jp
	=
	L^\ddagger  M\,.
\end{equation}
This initial solution is then standardized by dividing its entries by the
square roots of the corresponding diagonal entries in the \emph{resolution
matrix} \(L^\ddagger L\). As shown in~\cite{pascual-marqui-2002},
standardization balances out the initial solution that otherwise is \emph{a
priori} known to be biased towards the sensors.
Equation~\eqref{eq:sloreta-regularized-solution} is a solution to the
objective functional~\cite{pascual-marqui-2002}
\begin{equation}\label{eq:sloreta-objective-functional}
	F
	=
	\min_{\Jp,c}\nnorm{M - L\Jp - c\onesvec}{2}^{2}
	+
	\lambda\nnorm{\Jp}{2}^2
	\,,
\end{equation}
which represents a regularized fit between the measurements and the lead field
projection of the primary current \(\Jp\). The role of \(c\) is to set the
zero potential level of the electric field. The term \(\lambda\norm\Jp^2\) is
the regularization penalty function in which the coefficient \(\lambda\geq0\)
is the regularization parameter, chosen here according
to~\cite{rezaei--etal-2020}.

\subsection{\DipoleScan}\label{ssec:dipole-scan}

A second inverse method used in estimating source positions in this paper is
the so-called \emph{\DipoleScan} method~\cite{neugebauer-etal-2022}, where
following the inverse kernel idea of \figref{fig:forward-vs-inverse}, a filter
matrix \(\bff=\bff\arcs{\vec x}\) gives a local estimate \(\rec\Jp\) of the
primary current distribution \(\Jp\) at \(\vec x\) as
follows~\cite{van-veen--etal-1997}:
\begin{equation}\label{eq:beamformer-filter}
	\rec\Jp\arcs{\vec x} = \bff\arcs{\vec x} M\,.
\end{equation}
Here \(M\) is the measured data. The filter can be chosen to be the
Moore--Penrose pseudoinverse, which does not perform any kind of filtering on
the data~\cite{neugebauer-etal-2022}, or a \tsvd\ (\tSVD) of \(L\) at \(\vec
x\)~\cite{fuchs-1998}\cite{wolters--etal-1999}\cite{kaipio--somersalo-2004}.
The latter of these methods has a regularizing effect on the
solution~\cite{neugebauer-etal-2022}\cite{kaipio--somersalo-2004}, meaning the
high spatial noise components are not amplified.

This alone does not suffice for actually locating sources, since for that
purpose one needs an actual objective function to be optimized. \DipoleScan\
then minimizes the \emph{relative residual variance}
\begin{equation}\label{eq:rrv}
	\rrv\vec x = \frac{
		\norm{M\subt{avg} - L\bff\arcs{\vec x} M\subt{avg}}^2
	}{
		\norm{M\subt{avg}}^2
	} \,,
\end{equation}
with \(M\subt{avg}\) being an average measurement, or maximizes its
complement, the \emph{goodness of fit}
\begin{equation}\label{eq:gof}
	\gof\vec x = 1 - \rrv\vec x
\end{equation}
to determine the best source position.

\subsection{Head models and the computation of the lead field}\label{ssec:head-models}

\begin{figure}[h!]
	\begin{subfigure}[0.3\linewidth]{}{fig:ary-model-nomesh}
	\includegraphics[trim={4cm 1.5cm 6cm 3cm}, clip, width=\linewidth]{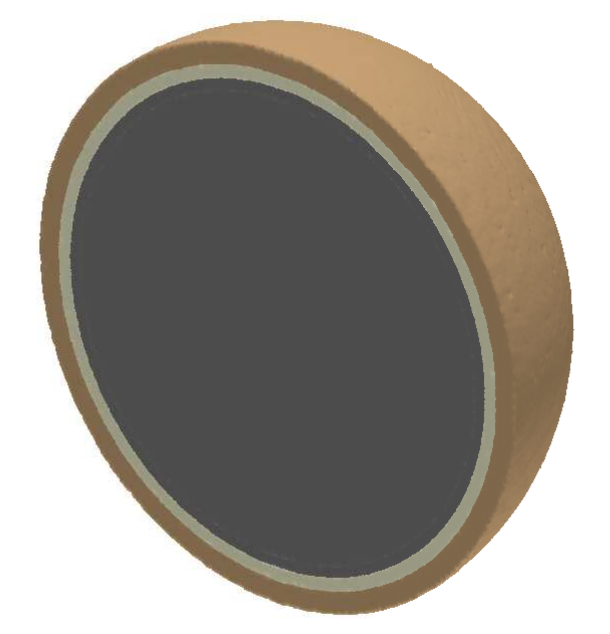}
	\end{subfigure}
	\hfill
	\begin{subfigure}[0.3\linewidth]{}{fig:ary-model-mesh}
		\includegraphics[trim={4cm 1.5cm 6cm 3cm}, clip, width=\linewidth]{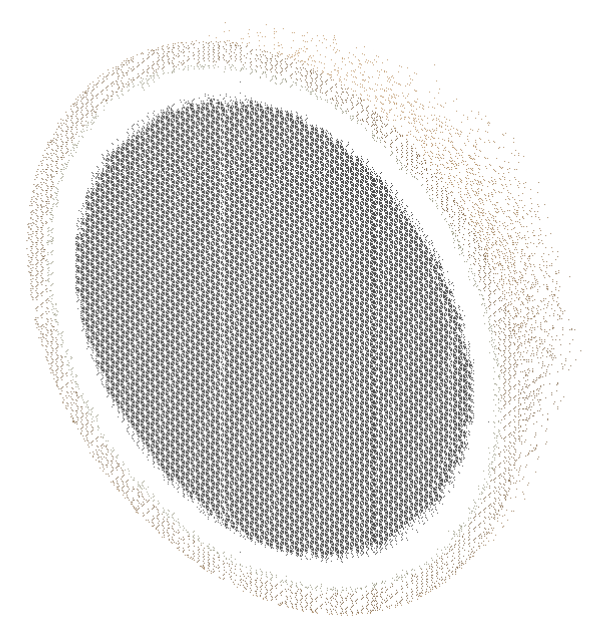}
	\end{subfigure}
	\hfill
	\begin{subfigure}[0.31\linewidth]{}{fig:ary-model-gray-matter-location}
		\includegraphics[width=\linewidth]{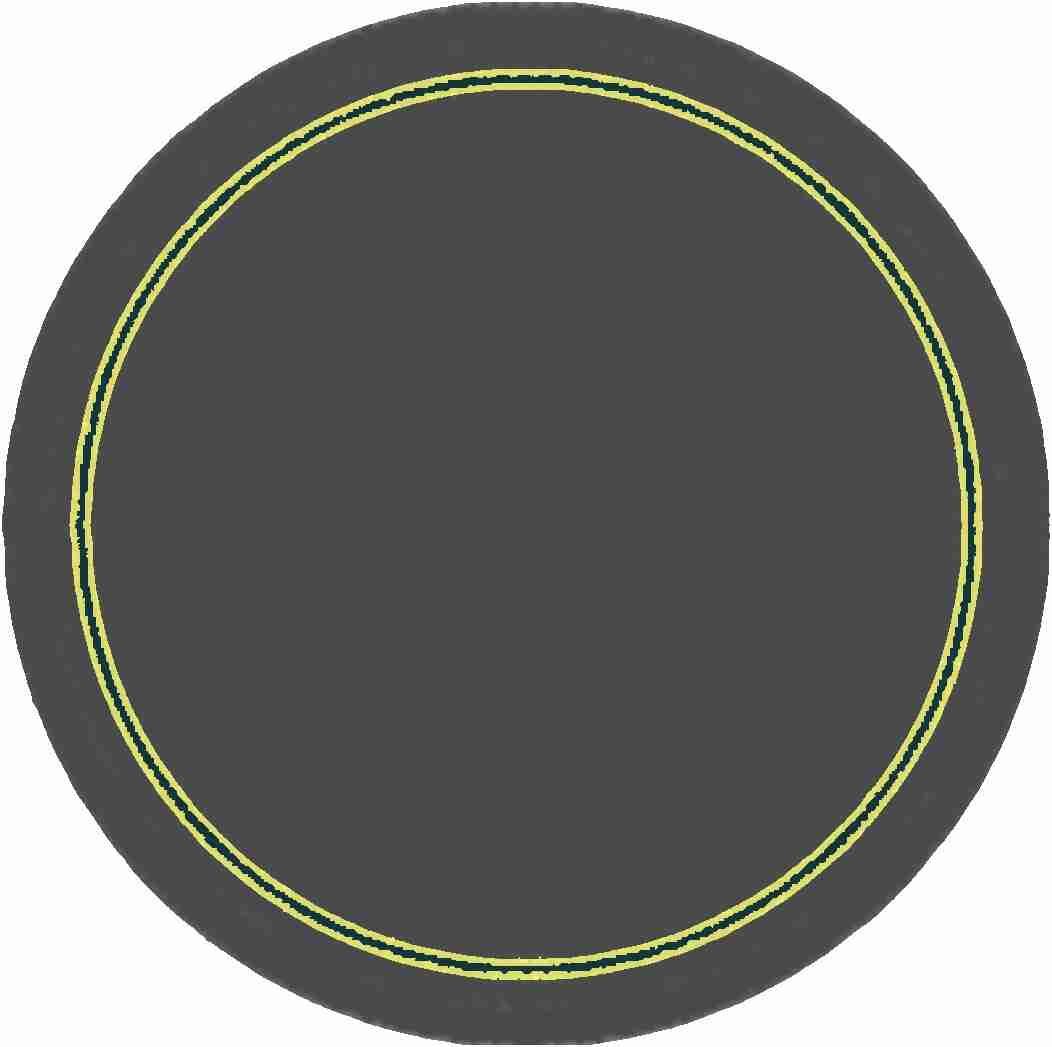}
	\end{subfigure}

	\caption{Cross sections of a \fe\ construction of an Ary sphere model. In
	\subfigref{fig:ary-model-nomesh}, the \FE\ edges are not displayed and the
	layers in the model are clearly visible. In
	\subfigref{fig:ary-model-mesh}, the edges and therefore the high
	resolution of the generated mesh are displayed.
	Subfigure~\subfigref{fig:ary-model-gray-matter-location} displays the
	location of the thin gray matter layer, with the yellow streaks indicating
	its inner and outer boundaries.}
	\label{fig:ary-model}
\end{figure}

To compute the forward and inverse solutions to the biomagnetic source
modelling problem presented in the previous subsection, one has to discretize
the volume conductor such that a computer can process it. To this end, \zi\
builds finite element domains from given MRI
segmentations~\cite{zeffiro-fem-algo}. \figref{fig:ary-model} presents a
cross-section of a spherical \num 3-layer finite element Ary
model~\cite{ary-etal-1981} generated in \zi.

\begin{figure}[h!]
	\begin{subfigure}[0.3\linewidth]{}{fig:real-head-a}
		\includegraphics[width=\linewidth]{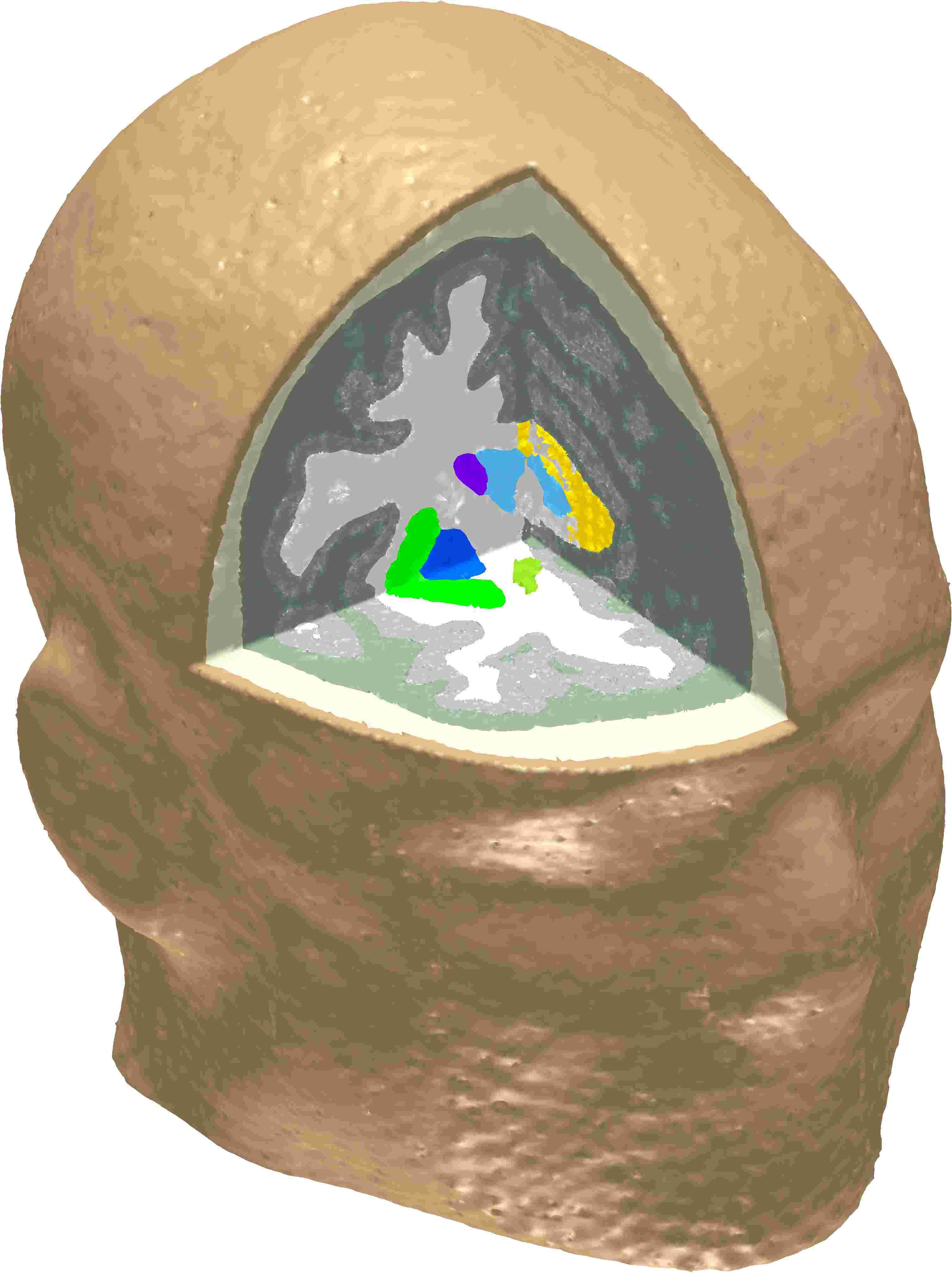}
	\end{subfigure}
	\hspace{1em}
	\begin{subfigure}[0.6\linewidth]{}{fig:real-head-b}
		\includegraphics[width=\linewidth]{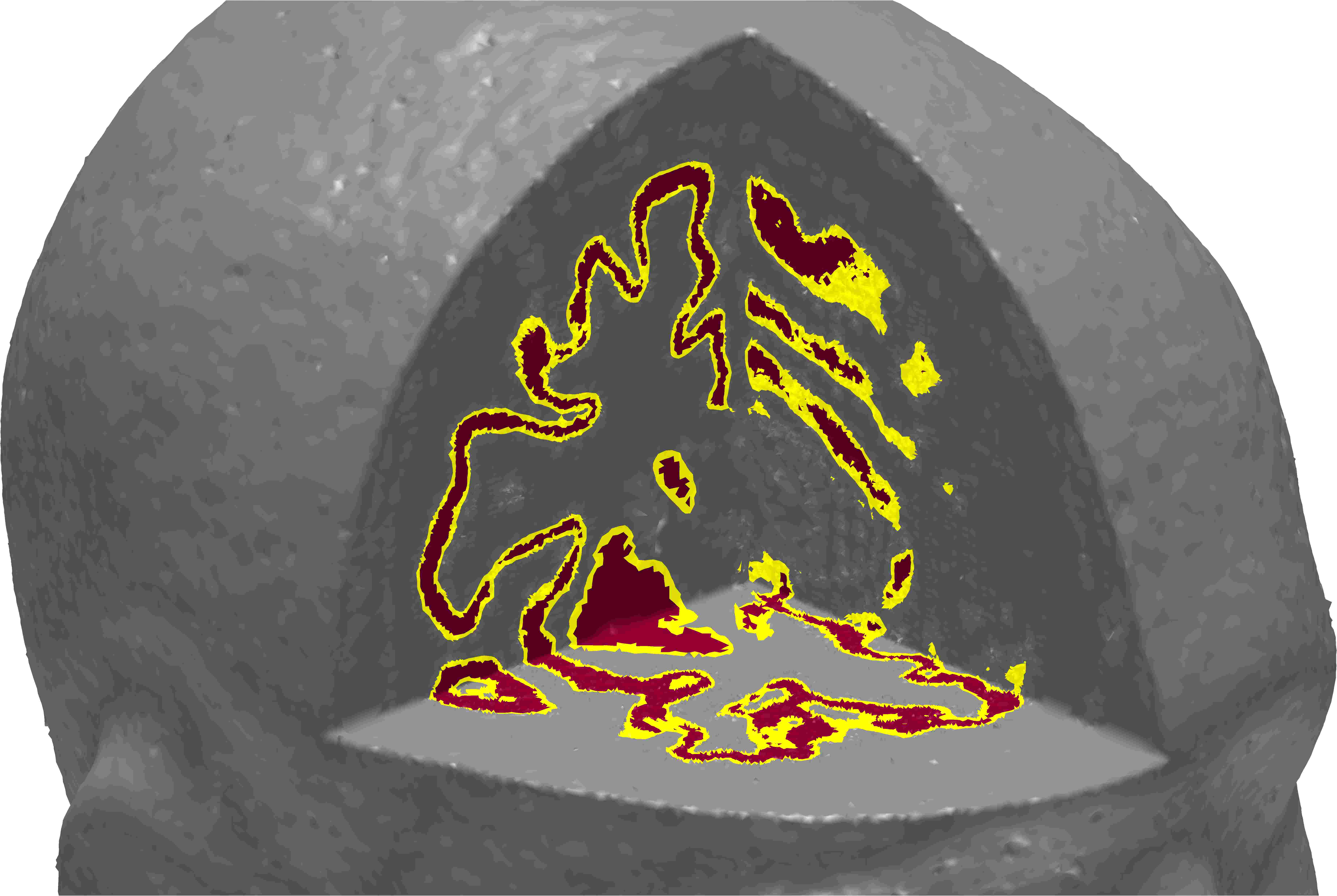}
	\end{subfigure}
	\caption{A realistic head model used in the computations with \zi.
	\subfigref{fig:real-head-a} shows a cross section of the head with
	different brain compartments visible, whereas \subfigref{fig:real-head-b}
	shows the active layers in dark red and yellow, with the yellow streaks
	indicating peeled off tetrahedra. The peeling depth used was
	\SI{0.1}{\milli\meter}. The active brain compartments include
	\emph{Cerebellum cortex}, \emph{Amygdala}, \emph{Thalamus},
	\emph{Caudate}, \emph{Nucleus accumbens}, \emph{Putamen},
	\emph{Hippocampus}, \emph{Pallidum}, \emph{Brain stem} and \emph{Ventral
	Diencephalon}.}
	\label{fig:real-head}
\end{figure}

The spherical volume conductor has an analytic or formulaic solution to the
computation of the lead field \(L\)~\cite{ary-etal-1981}, and hence it
provides a useful point of comparison when analysing the goodness of the
result. Later on, we analyse the differences between the numerical and
analytical solutions to the forward \EEG\ problem. In addition to the
spherical head model, we also take a look at the realistic human head
model~\cite{duneuro-head-model} seen in \figref{fig:real-head}.

It is also known that jump discontinuities in the conductivities of the brain
compartments affect the stiffness matrix needed in construction of the finite
element forward solution~%
\cite{wolters--etal-2007}\cite{drechsler--etal-2009}\cite{pursiainen--vorwerk--wolters-2016}\cite{piastra--etal-2021}.
Hence the amount of \emph{peeling} of the active brain layers seen in
\figref{fig:real-head}~\subfigref{fig:real-head-b} is also varied slightly, to
see how outliers in the reconstruction are affected.

To this end, the lead field routine of \zi\ was augmented to first peel off
the unwanted layers and only then start the lead field construction and source
positioning. The algorithm now has the following structure:
\begin{enumerate}
	\item
		use the peeling algorithm to select a true subset of the active brain
		elements and evenly distribute allowed source positions into those
		tetrahedra,
	\item
		build a system matrix \(A\)~\cite{braess-2007} for the \fe\ mesh,
	\item
		use \(A\) to build a transfer matrix \(T\)~\cite{weinstein-2000}%
		\cite{pursiainen--vorwerk--wolters-2016}\cite[9]{holtershinken-2021}
		by solving one linear system per sensor position,
	\item
		build a source space interpolation matrix \(D\) based on the \(\Hdiv\)
		source model~\cite{miinalainen-2019}\cite{atena-2021}%
		\cite{pursiainen--vorwerk--wolters-2016} with either
		\emph{position-based optimization} (PBO)~\cite{bauer--etal-2015} or
		\emph{mean position and orientation} (MPO)~\cite{atena-2021} used as
		optimization methods,
	\item
		compute the lead field as the matrix product \(L = T D\) and
	\item
		subtract the mean measurement from each column of \(L\) to set the
		zero potential level of the solution.
\end{enumerate}

To further disseminate on how the peeling algorithm
\code{zef\_deep\_nodes\_and\_tetra} functions, it
\begin{enumerate}
	\item
		determines the surface- and non-surface node indices of the given
		active layers by relying on the \zi\ function
		\code{zef\_surface\_mesh}.
	\item
		It then finds the indices of the nodes that are far enough from the
		surface nodes with Matlab's \code{rangesearch} functionality and
	\item
		observes which tetra have all \num 4 of their nodes far enough from
		the surface mesh by relying on the Matlab functions \code{ismember},
		\code{sum} and \code{find}.
\end{enumerate}

Note that step 3 results in the outermost layer being peeled off, regardless
of what the peeling depth is. This is desirable, as it removes the chance of
singularities due to discontinuities appearing in the solution.

\subsection{Measures for comparing the numerical and analytical lead fields}\label{ssec:analytical-error-measures}

To compare the \EEG\ lead fields produced by the forward solver in the case of
\figref{fig:ary-model}, the relative difference measure
\begin{equation}
	\label{eq:rdm}
	\mRDM
	=
	\nnorm*{
		\frac{\Ln}{\nnorm{\Ln}{2,1}}
		-
		\frac{\La}{\nnorm{\La}{2,1}}
	}{2,1}
\end{equation}
and magnitude measure
\begin{equation}
	\label{eq:mag}
	\mMAG
	=
	\abs*{1 - \frac{\nnorm{\La}{2,1}}{\nnorm{\Ln}{2,1}}}
\end{equation}
were employed to compute differences between the analytical~\cite{ary-etal-1981}
and numerical lead fields \(\La\) and \(\Ln\). Here \(\nnorm{L}{2,1}\) denotes
the \(2\)-norm of \(L\) along the rows of \(L\).

\subsection{Evaluation of the localisation error}\label{ssec:localisation-error}

To evaluate the localisation error \(\Delta\) of
equation~\eqref{eq:localisation-error-components}, \num{10 000} synthetic
sources were placed evenly into the active regions of the volume conductor:
\emph{Cerebellum cortex}, \emph{Amygdala}, \emph{Thalamus}, \emph{Caudate},
\emph{Nucleus accumbens}, \emph{Putamen}, \emph{Hippocampus}, \emph{Pallidum},
\emph{Brain stem} and \emph{Ventral Diencephalon}. Then a lead field \(L\)
corresponding to these sources was computed, as specified at the end of
subsection~\ref{ssec:head-models}. For each source given by the
position--direction--amplitude triplet \(\source {\vec x} {\vec d} a\),
synthetic measurements \(M\subi{\source {\vec x} {\vec d} a}\) were
constructed according to
\begin{equation}\label{eq:synthetic-measurements}
	M\subi{\source {\vec x} {\vec d} a}
	=
	a\,\restr{L}{ S, \vec x}\,\vec d + 10^{-\SNR / 20} N\,,
\end{equation}
with \(\norm{\vec d} = 1\), \(a > 0\), and the notation \(\restr{L}{S,\vec
x}\) indicating a restriction of \(L\) to the rows corresponding to all
sensors \(S\) and the columns corresponding to the dipole position \(\vec x\).
Here \(\SNR\) is the noise level in decibels, as in \(\qunit\SNR =
\si\decibel\), and \(N\) is a normally distributed random variable, with mean
\(\mu = 0\) and variance \(\sigma^2=1\).

With the simulated measurements \(M\subi{\source {\vec x} {\vec d} a}\) in
place, they were then inverted with \sLORETA\ and \DipoleScan, to produce a
distribution of reconstructions: one \(\inv\vec d\subi{i}\) for each \(\vec
x_i\) in the original source space%
\footnote{Constituting an inverse crime.~\cite{kaipio--somersalo-2004}}.
The position of the most focal reconstruction \(\vec x\subi{I}\) was sought by
finding the index \(I\) of the direction with the largest norm or dipole
moment, as in
\begin{equation}\label{eq:max-norm-index}
	I = \indmax_i\nnorm{\inv\vec d\subi{i}} 2 \,,
\end{equation}
following the \emph{maximum principle of source
localization}~\cite{hamalainen--ilmoniemi-1994}\cite{pascual-marqui-2002}\cite{mohamed-2013}.
The localisation error was then computed as
\begin{equation}\label{eq:localisation-error}
	\Delta = \frac 1 {\sqrt 3}\nnorm{\vec x_i - \vec x_I} 2 \,,
\end{equation}
as in the differences among the positions in the original source space.
The \(\sqrt 3\) in \eqref{eq:localisation-error} is the norm of a Cartesian dipole.

\subsection{Spatial dispersion}\label{ssec:spatial-dispersion}

As \figref{fig:forward-vs-inverse} implies, to construct an estimate
\(\rec\Jp\) of the original synthetic source distribution \(\Jp\) in a linear
fashion, one can multiply a lead field \(L\) by an inverse kernel \(K\) to
obtain a so-called \emph{resolution matrix} \(R\), such that \(\rec\Jp = KL\Jp
= R\Jp\). In such a case, the columns of \(R\) or \(\vec R_i\) describe how
the activity at the corresponding source positions \(\vec x_i\) are blurred in
this process~\cite{hauk-2022-neuroimage}. The width of this blurring is given
by the so-called \emph{spatial dispersion} measure~\cite{molins-2008}
\begin{equation}\label{eq:dispersion}
	\SD_i = \sqrt{
		\frac{
			\sum_k \arcs{d_{k,i} \norm{\vec p \subi{k}}} ^ 2
		}{
			\sum_k \norm{\vec p \subi{k}} ^ 2
		}
	} \,,
\end{equation}
where \(k\) ranges over the indices of the source positions \(\vec x_k\) in a
given ROI around \(\vec x_i\), \(\vec p_k\) are the dipole moments of the
reconstructions in the ROI and \(d\subi{k,i}\) are the distances of each
dipole in the ROI from \(\vec x_i\). In this study, it is observed how varying
the peeling depth \(\peeld\) affects the measure in
equation~\eqref{eq:dispersion}.

\section{Results}\label{sec:results}

\subsection{Evaluating the peeling algorithm}\label{ssec:evaling-the-peeling-algorithm}

\begin{figure}[h!]
	\begin{subfigure}[0.49\linewidth]{Unrefined, \(\peeld=\SI{0.5}{mm}\).}{fig:ary_zoomed_unrefined_r2mm_d0.5mm}
		\includegraphics[width=\linewidth]{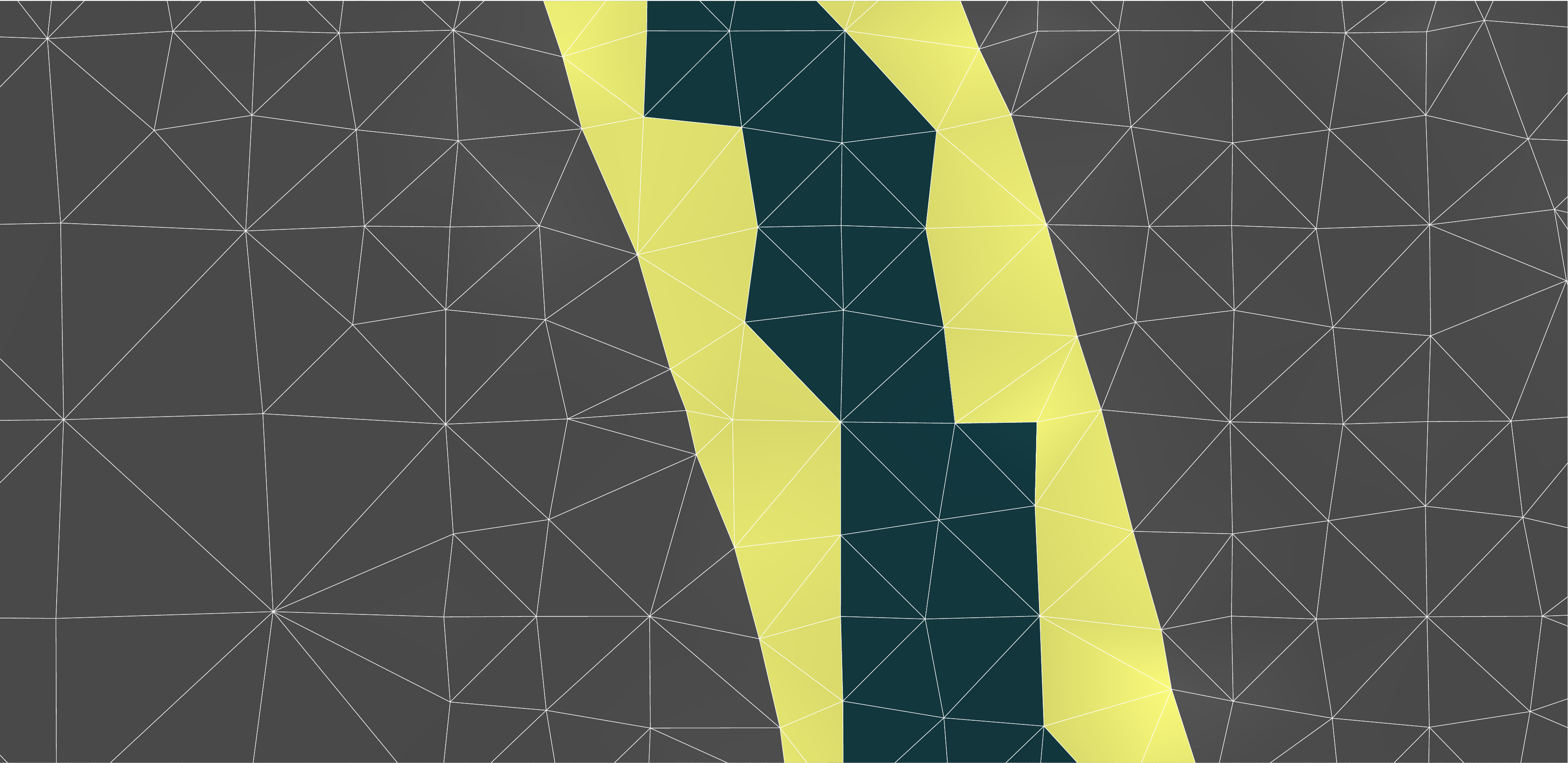}
	\end{subfigure}
	\hfill
	\begin{subfigure}[0.49\linewidth]{Refined, \(\peeld=\SI{0.5}{mm}\).}{fig:ary_zoomed_refined_r2mm_d0.5mm}
		\includegraphics[width=\linewidth]{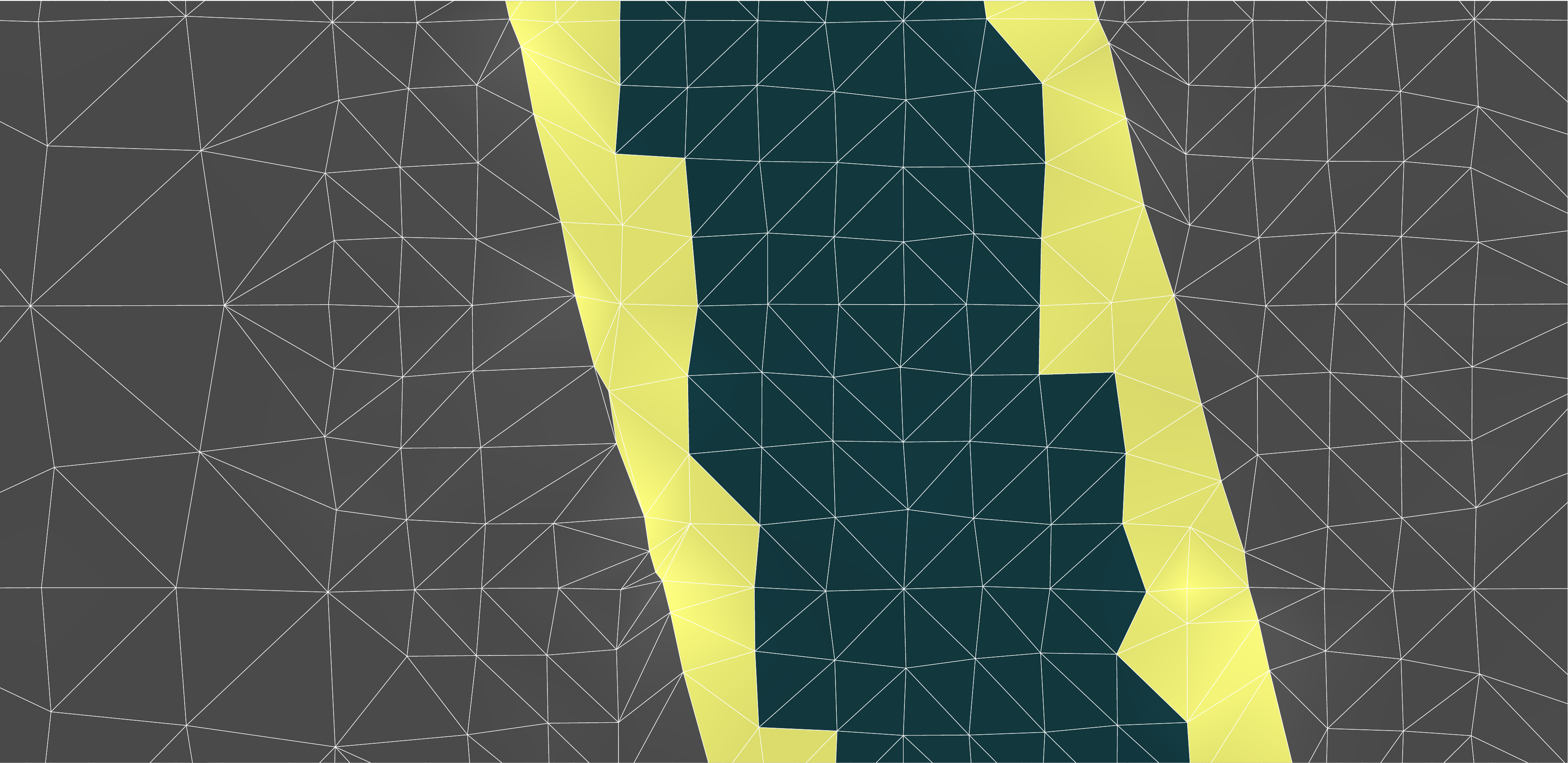}
	\end{subfigure}

	\vspace*{1ex}

	\begin{subfigure}[0.49\linewidth]{Unrefined, \(\peeld = \SI{1.0}{mm}\).}{fig:ary_zoomed_unrefined_r2mm_d1.0mm}
		\includegraphics[width=\linewidth]{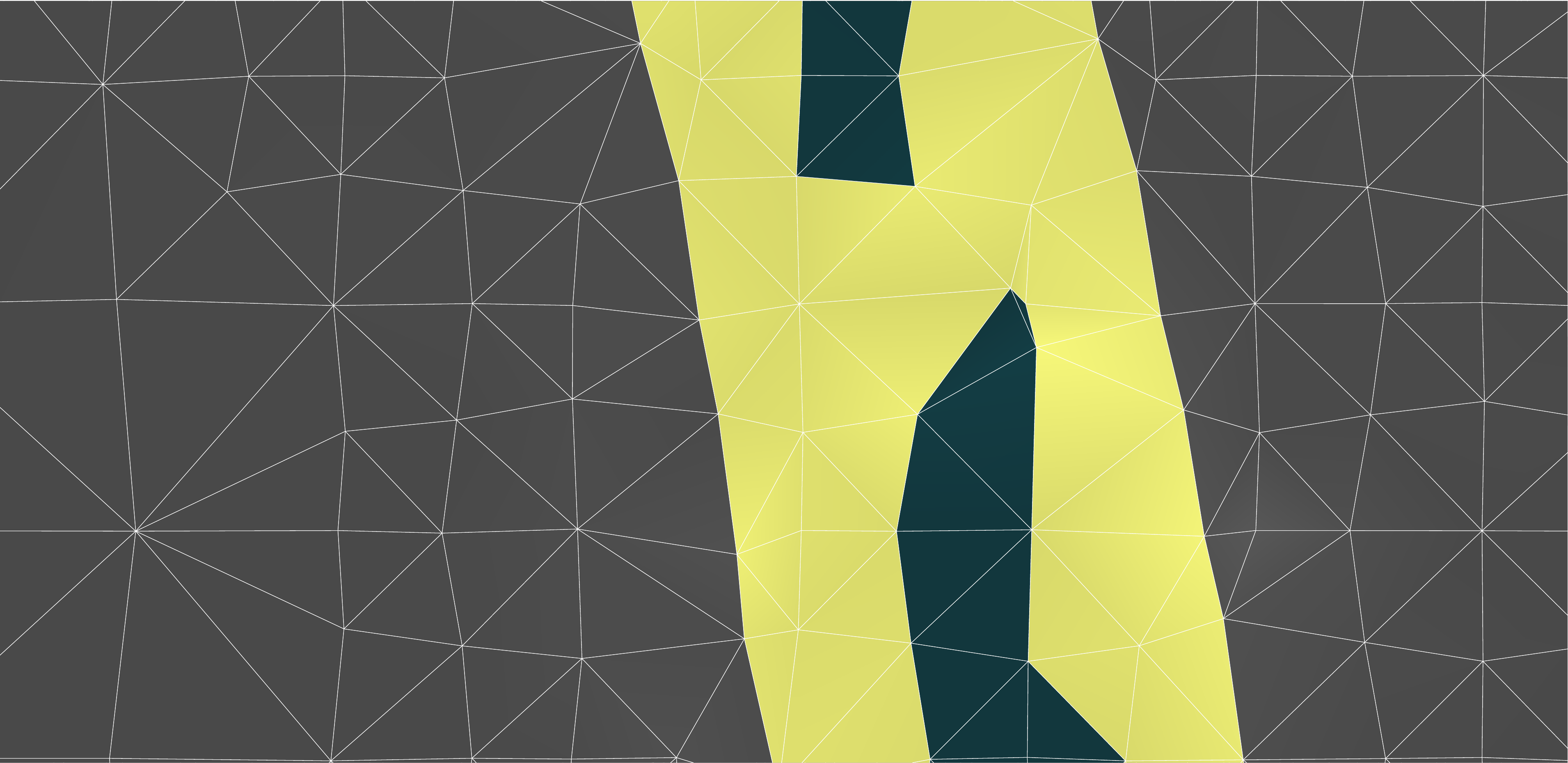}
	\end{subfigure}
	\hfill
	\begin{subfigure}[0.49\linewidth]{Refined, \(\peeld=\SI{1.0}{mm}\).}{fig:ary_zoomed_refined_r2mm_d1.0mm}
		\includegraphics[width=\linewidth]{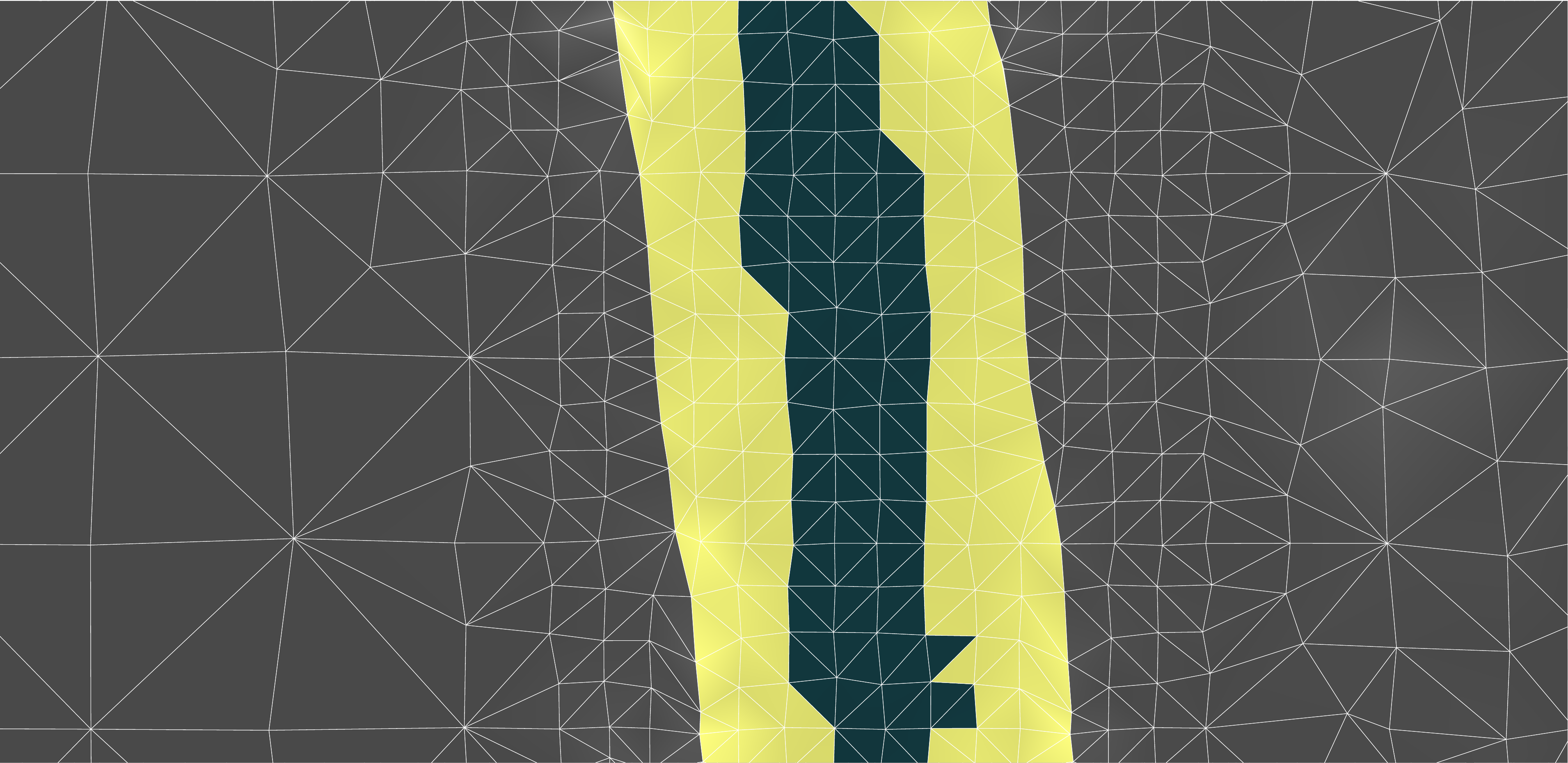}
	\end{subfigure}

	\caption{%
		The effects of mesh refinement on \emph{peeling}, or the disqualification of
		tetrahedra from the set of valid source positions. The yellow triangles
		indicate which tetrahedra were peeled off the inner and outer surfaces of the
		thin gray matter layer, after the peeling algorithm had been applied to the
		mesh. In \subfigref{fig:ary_zoomed_unrefined_r2mm_d0.5mm} and
		\subfigref{fig:ary_zoomed_refined_r2mm_d0.5mm} the peeling depth \(\peeld\)
		has been set at \SI{0.5}{\milli\meter} whereas in
		\subfigref{fig:ary_zoomed_unrefined_r2mm_d1.0mm} and
		\subfigref{fig:ary_zoomed_refined_r2mm_d1.0mm} we have \(\peeld =
		\SI{1.0}{\milli\meter}\). A more refined mesh produces a more consistent
		peeling outcome.%
	}

	\label{fig:ary_zoomed_r2mm_d0.5mm}

\end{figure}

We start off with a discussion on how well the peeling algorithm itself works.
\figref{fig:ary_zoomed_r2mm_d0.5mm} presents the effects of refining the gray matter
layer of the Ary model of \figref{fig:ary-model} on its \emph{peeling}, the
disqualification of nodes and associated tetrahedra from the set of valid source
positions. The peeling was done within the distances of \num{0.5} and
\SI{1.0}{\milli\meter} from the inner and outer surfaces of the thin gray matter
layer. We consider peeling to be necessary, because placing \(\Hdiv\) dipoles in
positions of conductivity discontinuities would cause significant forward
errors~\cite{pursiainen--vorwerk--wolters-2016}.

\begin{figure}[h!]
	\includegraphics[width=\linewidth]{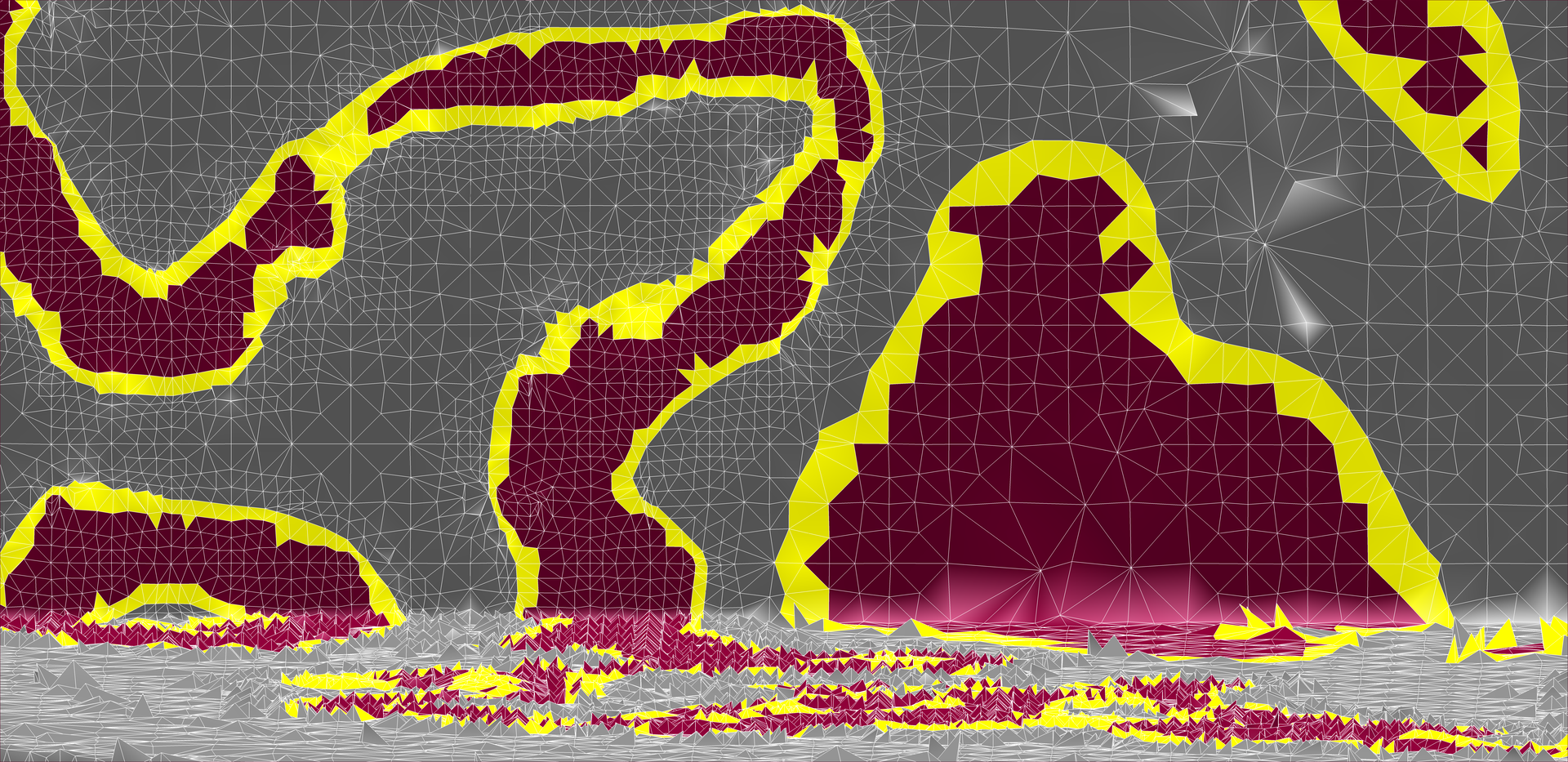}
	\caption{A portion of the \SI{4}{\milli\meter} thick active layers (dark red +
	yellow) of \figref{fig:real-head} \subfigref{fig:real-head-b}, displayed
	with the edges of the finite element mesh visualized, but with only
	\SI{0.1}{\milli\meter} peeling (yellow).}
	\label{fig:real-head-active-layer-edges}
\end{figure}

As can be observed, refinement plays an important role in the peeling process:
it prevents an excessive reduction of possible source positions from the
active layer. This can be observed in the unrefined \SI{1.0}{\milli\meter}
case shown in Subfigure~\subfigref{fig:ary_zoomed_unrefined_r2mm_d1.0mm},
where a hole is punched through the gray matter layer. This means that a
perfectly valid dipolar source location is excluded from the possible set of
source positions during forward modelling, or the computation of the lead
field matrix \(L\). Taking a closer look at the realistic head model of
\figref{fig:real-head}~\subfigref{fig:real-head-b} also displays a similar
effect with a peeling depth of \SI{0.1}{\milli\meter}, which is displayed in
\figref{fig:real-head-active-layer-edges}.

Peeling will remove all those tetrahedra with one or more nodes closer to
surface than a given peeling depth. The reason for this is to make sure that
at least one layer of tetrahedra is removed, so that sources absolutely cannot
be placed in tetrahedra, that are right next to another compartment with
possibly differing conductivity.

It turns out that even with the refinement performed on the surface of the
active gray matter layer, there are still parts of the compartment which are
not fine enough with \SI{0.1}{\milli\meter} peeling, and hence holes in
possible source positions, such as the ones seen in the upper left and right
corners of \figref{fig:real-head-active-layer-edges}, are formed. This is
again due to the requirement that at least one tetrahedral surface layer is
removed from the set of possible source positions, which results in the
effective peeling depth being greater than the low numerical value provided by
the user. The peeling algorithm then seems to function as intended.

\subsection{Measuring the goodness of the solver against an analytical model}\label{ssec:inverting-ary-model}

\begin{figure*}[h!]
	\begin{subfigure}{}{fig:pbo-eccentricity-rdm}
		\includegraphics[width=0.9\linewidth,trim={0pt 1.2cm 0pt 0pt}, clip]{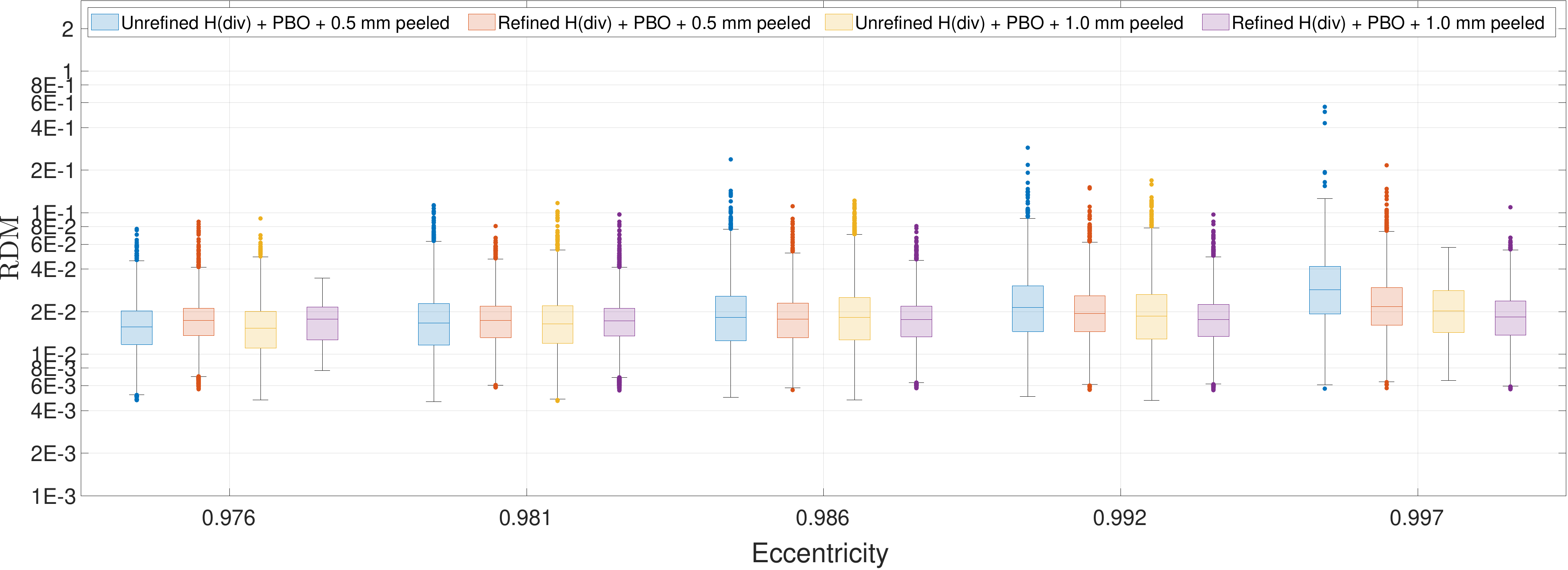}
	\end{subfigure}

	\medskip

	\begin{subfigure}{}{fig:pbo-eccentricity-mag}
		\includegraphics[width=0.9\linewidth]{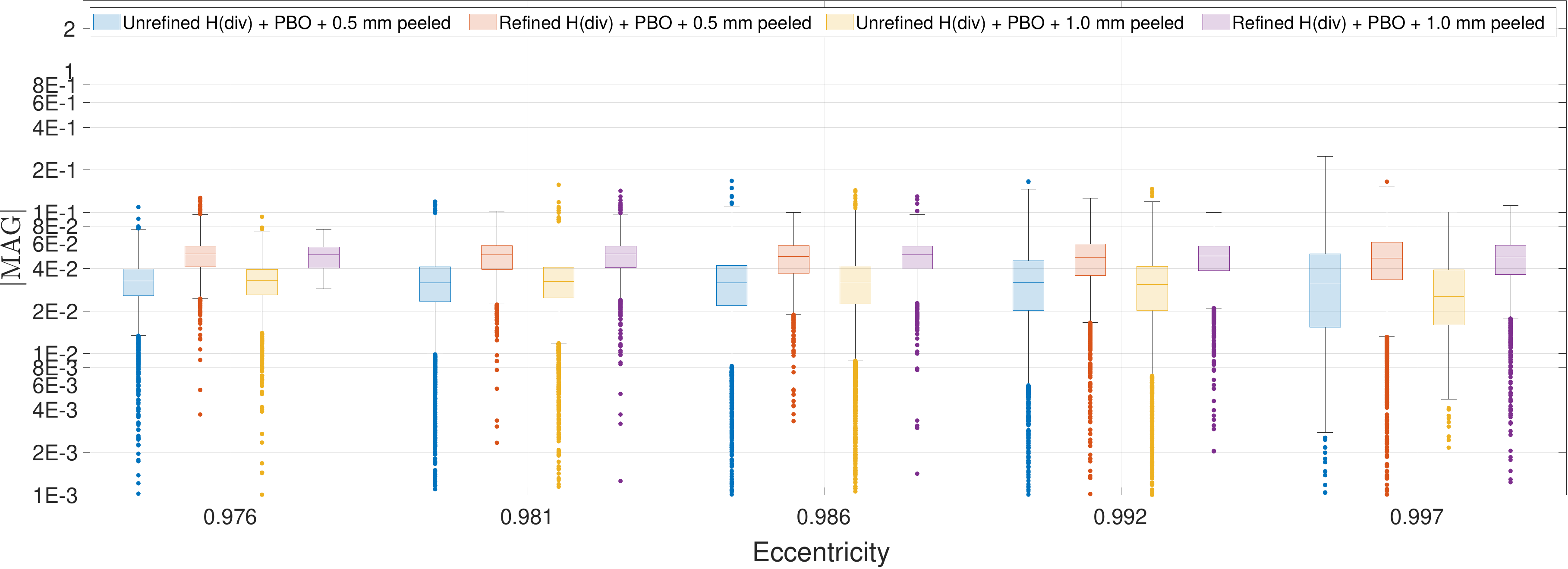}
	\end{subfigure}

	\caption{High-eccentricity \eqref{eq:rdm} \subfigref{fig:pbo-eccentricity-rdm} and \eqref{eq:mag} \subfigref{fig:pbo-eccentricity-mag} for PBO in the case of a spherical Ary model.}
	\label{fig:pbo-eccentricity}
\end{figure*}

Figure~\ref{fig:pbo-eccentricity} contains visualizations of the measures
\eqref{eq:rdm} and \eqref{eq:mag} for \PBO, in the case of the analytical
model of \figref{fig:ary-model}. The base resolution of the mesh was
\SI{2}{\milli\meter}, with refinement referring to the surface of the active
layer being refined as seen in \figref{fig:ary_zoomed_r2mm_d0.5mm}.
\emph{Eccentricity} refers to the relative radius of the position at which the
comparison was performed. As suggested by~\cite{drechsler--etal-2009},
eccentricities of over \SI{98}{\percent} are of special interest, as they
correspond to where a primary dipole might be physiologically located,
somewhere between the external granular and pyramidal layers (layers \num
2--\num 3) of the cerebral cortex, and hence observed here.

With \eqref{eq:rdm}, in almost all cases the \PBO\ medians remain below the
\num{0.02} limit, except at the highest two eccentricities. The upper outlier
quantiles \(q\subt{u} = q\subi{75\%} + 1.5\arcs{q\subi{75\%} -
q\subi{25\%}}\), which the whiskers of \figref{fig:pbo-eccentricity}
correspond to, mostly remain close to \num{0.1}, with the case of unrefined
\(\peeld = \SI{0.5}{\milli\meter}\) approaching \num{0.02} at the two highest
eccentricities. With \eqref{eq:mag}, the median remains below \num{0.04}, with
the upper quantile \(q\subt{u}\) behaving similarly to what was observed with
\eqref{eq:rdm}.

Increasing the resolution of the finite element mesh near active layer
boundaries via mesh refinement seems to be mostly reducing \eqref{eq:rdm},
especially towards the higher eccentricities. For \eqref{eq:mag}, the
refinement actually seems to increase the median differences between the
analytical and numerical lead fields by roughly \num{0.02}, which is seen in
the horizontal middle lines of the box plots of \figref{fig:pbo-eccentricity}.
However, the quantiles \(q\subt u\) were lowered by a few percent with
\(\peeld = \SI{0.5}{\milli\meter}\), suggesting that refinement has a net
positive statistical impact on the result.

\subsection{Source localisation in a realistic head model}\label{sec:inverting-head-models}

In the spirit of Cuffin et al.~\cite{cuffin-2001},
Tables~\ref{tab:sloreta-unfiltered-inverse-statistics}--\ref{tab:dipolescan-unfiltered-inverse-statistics}
present mean values \(\mu\) of localisation errors \(\Delta\) of
equation~\eqref{eq:localisation-error} and their standard deviations
\(\sigma\) for \num{10 000}-source lead fields \(L\), corresponding to
different peeling depths \(\peeld\). The lead fields were inverted \num{20}
times with \sLORETA~\cite{pascual-marqui-2002} and
\DipoleScan~\cite{neugebauer-etal-2022} at different measurement noise levels,
each inversion having a different white noise realization. The cells of the
Tables~\ref{tab:sloreta-unfiltered-inverse-statistics}--\ref{tab:dipolescan-unfiltered-inverse-statistics}
are color mapped based on the largest \(\Delta\) and \(\sigma\) in each table.

\begin{table}[!h]
	\def\avgmax{41.76}%
	\def\stdmax{18.90}%
	\caption{\sLORETA\ average localisation error \(\Delta\) statistics. Here \(\peeld\) is the peeling depth, \(\mu\) is the sample mean and \(\sigma\) is the sample standard deviation.}
	\label{tab:sloreta-unfiltered-inverse-statistics}
	\resizebox{\linewidth}{!}{
		\begin{tabular}{|c|cc|cc|cc|}
			\toprule
			\(\peeld\) (\unit{mm}) & \multicolumn{2}{c|}{0.0} & \multicolumn{2}{c|}{\num{0.5}} & \multicolumn{2}{c|}{1.0} \\
			\midrule
			SNR (\unit{dB}) & \(\mu\) (\unit{mm}) & \(\sigma\) (\unit{mm}) & \(\mu\) (\unit{mm}) & \(\sigma\) (\unit{mm}) & \(\mu\) (\unit{mm}) & \(\sigma\) (\unit{mm}) \\
			\midrule
			5  & \hmcell{41.76}{\avgmax} & \hmcell{18.88}{\stdmax} & \hmcell{41.57}{\avgmax} & \hmcell{18.90}{\stdmax} & \hmcell{41.57}{\avgmax} & \hmcell{18.72}{\stdmax} \\
			10 & \hmcell{30.84}{\avgmax} & \hmcell{15.62}{\stdmax} & \hmcell{30.62}{\avgmax} & \hmcell{15.55}{\stdmax} & \hmcell{30.77}{\avgmax} & \hmcell{15.43}{\stdmax} \\
			15 & \hmcell{20.45}{\avgmax} & \hmcell{09.57}{\stdmax} & \hmcell{20.35}{\avgmax} & \hmcell{09.52}{\stdmax} & \hmcell{20.66}{\avgmax} & \hmcell{09.54}{\stdmax} \\
			20 & \hmcell{13.81}{\avgmax} & \hmcell{05.69}{\stdmax} & \hmcell{13.73}{\avgmax} & \hmcell{05.71}{\stdmax} & \hmcell{14.01}{\avgmax} & \hmcell{05.80}{\stdmax} \\
			25 & \hmcell{10.76}{\avgmax} & \hmcell{03.36}{\stdmax} & \hmcell{10.70}{\avgmax} & \hmcell{03.36}{\stdmax} & \hmcell{10.78}{\avgmax} & \hmcell{03.46}{\stdmax} \\
			30 & \hmcell{09.87}{\avgmax} & \hmcell{02.04}{\stdmax} & \hmcell{09.81}{\avgmax} & \hmcell{02.05}{\stdmax} & \hmcell{09.79}{\avgmax} & \hmcell{02.12}{\stdmax} \\
			\bottomrule
		\end{tabular}
	}
\end{table}

\begin{table}[!h]
	\def\avgmax{26.88}%
	\def\stdmax{13.60}%
	\caption{\DipoleScan\ average localisation error \(\Delta\) statistics. The meaning of notations is the same as in \tabref{tab:sloreta-unfiltered-inverse-statistics}.}
	\label{tab:dipolescan-unfiltered-inverse-statistics}
	\resizebox{\linewidth}{!}{
		\begin{tabular}{|c|cc|cc|cc|}
			\toprule
			\(\peeld\) (\unit{mm}) & \multicolumn{2}{c|}{0.0} & \multicolumn{2}{c|}{\num{0.5}} & \multicolumn{2}{c|}{1.0} \\
			\midrule
			SNR (\unit{dB}) & \(\mu\) (\unit{mm}) & \(\sigma\) (\unit{mm}) & \(\mu\) (\unit{mm}) & \(\sigma\) (\unit{mm}) & \(\mu\) (\unit{mm}) & \(\sigma\) (\unit{mm}) \\
			\midrule
			5  & \hmcell{26.87}{\avgmax} & \hmcell{13.57}{\stdmax} & \hmcell{26.83}{\avgmax} & \hmcell{13.60}{\stdmax} & \hmcell{26.88}{\avgmax} & \hmcell{13.53}{\stdmax} \\
			10 & \hmcell{16.01}{\avgmax} & \hmcell{7.92}{\stdmax}  & \hmcell{15.97}{\avgmax} & \hmcell{7.92}{\stdmax}  & \hmcell{16.08}{\avgmax} & \hmcell{7.99}{\stdmax}  \\
			15 & \hmcell{9.13}{\avgmax}  & \hmcell{4.66}{\stdmax}  & \hmcell{9.11}{\avgmax}  & \hmcell{4.66}{\stdmax}  & \hmcell{9.17}{\avgmax}  & \hmcell{4.75}{\stdmax}  \\
			20 & \hmcell{5.09}{\avgmax}  & \hmcell{3.01}{\stdmax}  & \hmcell{5.07}{\avgmax}  & \hmcell{3.01}{\stdmax}  & \hmcell{5.07}{\avgmax}  & \hmcell{3.05}{\stdmax}  \\
			25 & \hmcell{2.65}{\avgmax}  & \hmcell{2.01}{\stdmax}  & \hmcell{2.64}{\avgmax}  & \hmcell{2.01}{\stdmax}  & \hmcell{2.64}{\avgmax}  & \hmcell{2.00}{\stdmax}  \\
			30 & \hmcell{1.18}{\avgmax}  & \hmcell{1.22}{\stdmax}  & \hmcell{1.18}{\avgmax}  & \hmcell{1.22}{\stdmax}  & \hmcell{1.20}{\avgmax}  & \hmcell{1.21}{\stdmax}  \\
			\bottomrule
		\end{tabular}
	}
\end{table}

Here \DipoleScan\ produces superior localisation results when compared to \sLORETA,
with both low and high SNR levels, which is an expected result in search of a single
source, matching the prior model of \DipoleScan. Especially in the case of
\DipoleScan, the mean error \(\Delta\approx\SI{10.5}{\milli\meter}\) reported by
Cuffin is reached already at \(\SNR = \SI{15}{\decibel}\), whereas with \sLORETA\ a
value as high as \(\SI{30}{\decibel}\) has to be used, for comparable values to
manifest themselves.

\begin{figure}[h!]
	\includegraphics[width=\linewidth]{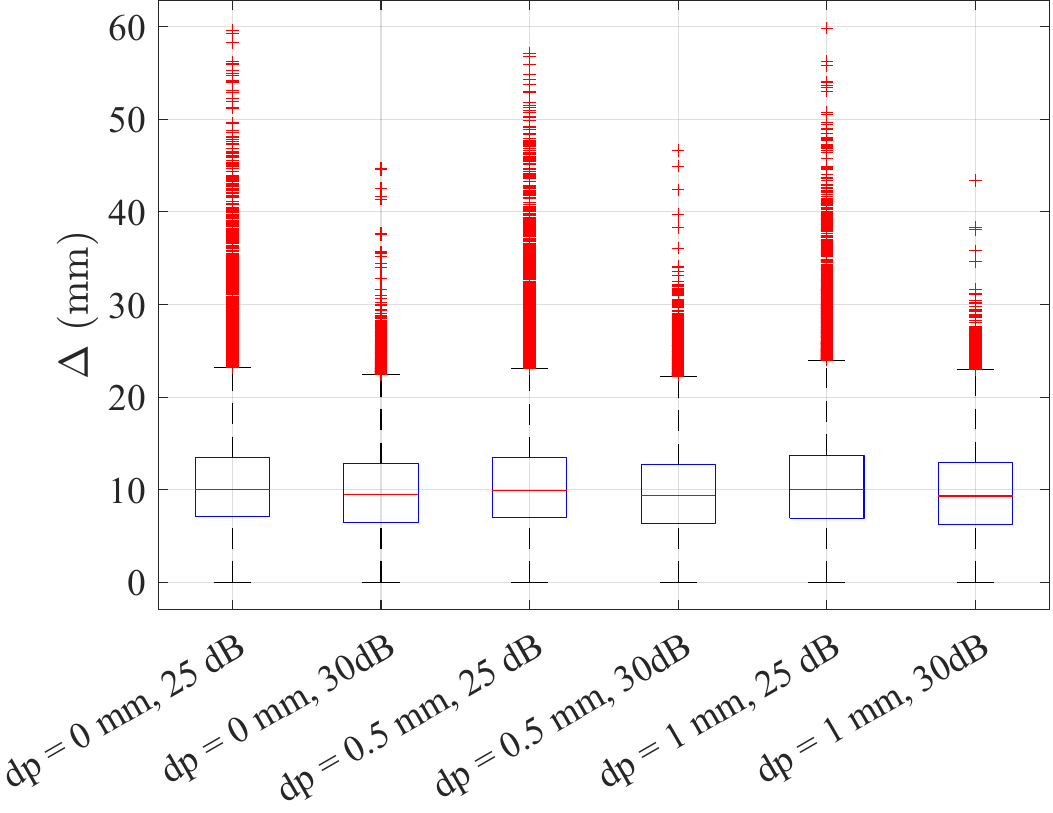}
	\caption{%
		\sLORETA\ localisation errors \(\Delta\) against two lowest noise
		levels and all peeling depths \(\peeld\).%
	}
	\label{fig:sloreta-outliers}
\end{figure}

\begin{figure}[h!]
	\includegraphics[width=\linewidth]{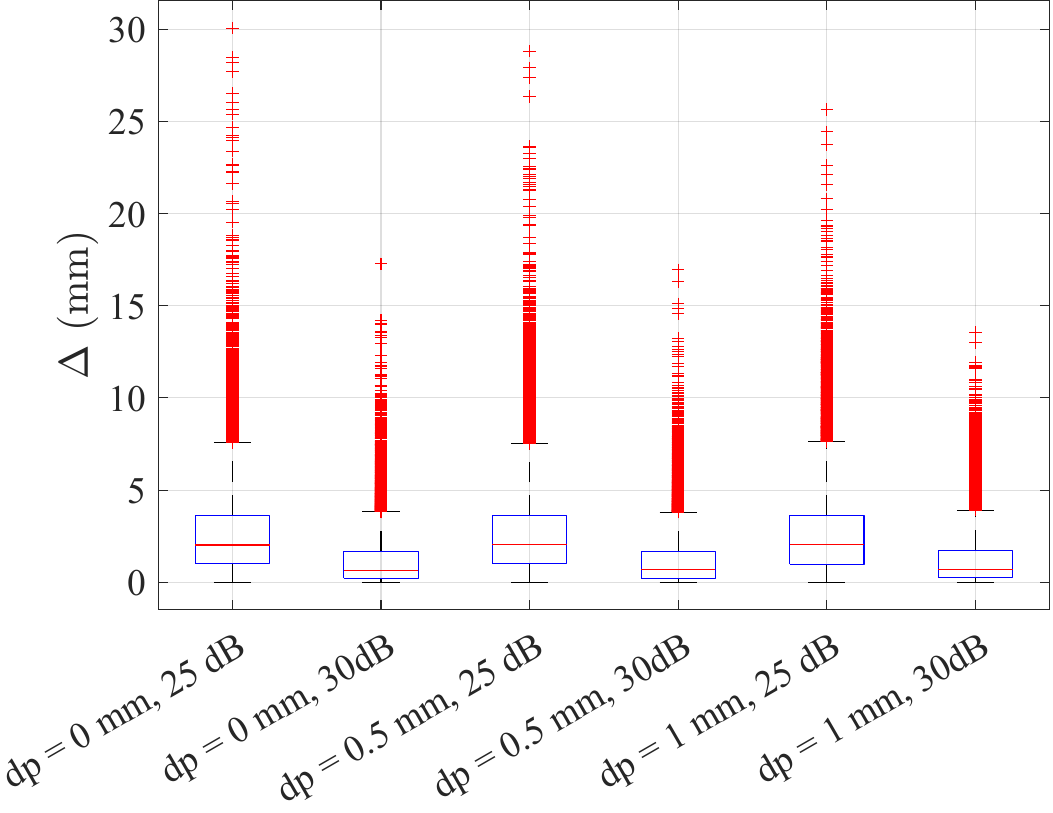}
	\caption{%
		\DipoleScan\ localisation errors \(\Delta\) against two lowest noise
		levels and all peeling depths \(\peeld\).%
	}
	\label{fig:dipolescan-outliers}
\end{figure}

To graphically observe how the peeling affects outliers of \(\Delta\),
Figures~\ref{fig:sloreta-outliers} and \ref{fig:dipolescan-outliers} were
formed. They display box plots of \(\Delta\) against the two lowest noise
levels in the case of \sLORETA\ and \DipoleScan, respectively.

In the case of \sLORETA, going from a peeling of \SI{0}{\milli\meter} to
\SI{0.5}{\milli\meter} results in a disappearance of at least few outlier
markers. Further, moving from \SI{0.5}{\milli\meter} to \SI{1.0}{\milli\meter}
peeling further reduces the outliers, which is seen in the sparsity of the
outlier point clouds. With \DipoleScan, the peeling results are similar, as
the extent of outliers is reduced with increased peeling.

To further observe how peeling affects the outliers,
Tables~\ref{tab:sloreta-above-mu-plus-2sigma} and
\ref{tab:dipolescan-above-mu-plus-2sigma} show the per-noise-level numbers of
sources, whose localisation errors satisfy \(\Delta > \mu + 2\sigma\), for a
sample of 20 reconstructions obtained with a given noise level. Here \(\mu\)
and \(\sigma\) are the expected sample mean and standard deviation values of
\(\Delta\), provided by Cuffin et al~\cite[Table 1]{cuffin-2001}. These are
contrasted against the numbers of columns in each respective \num{10
000}-source lead field \(L\) by color-mapping each data point against the
largest number of columns. Each Cartesian source has coordinates pointing in
the \(x\)-, \(y\)- and \(z\)-directions, and hence the number of columns in
\(L\) is threefold, when compared to the mentioned number of sources.

\begin{table}[!h]
	\def\maxval{30 459}%
	\caption{The number of outliers with \(\Delta > \mu + 2\sigma\), for a sample of \num{20} \sLORETA\ reconstructions. The expected means \(\mu\) and standard deviations \(\sigma\) for each noise level have been gathered from~\cite[Table 1]{cuffin-2001}. Again, \(\peeld\) refers to the peeling depth.}
	\label{tab:sloreta-above-mu-plus-2sigma}
	\resizebox{\linewidth}{!}{
		\begin{tabular}{|c|c|c|ccc|ccc|}
			\toprule
			\multicolumn{3}{|c|}{\(\peeld\) (\unit{mm})} & \num{0.0} & \num{0.5} & \num{1.0} \\ % & \num{0.0} & \num{0.5} & \num{1.0} \\
			\midrule
			\multicolumn{3}{|c|}{Number of columns in \(L\)} & \hmcell{30 288}{\maxval} & \hmcell{30 459}{\maxval} & \hmcell{29 367}{\maxval} \\ %& \hmcell{29 835}{\maxval} & \hmcell{30 003}{\maxval} & \hmcell{28 926}{\maxval}  \\
			\midrule
			SNR (\unit{dB}) & \(\mu\) (\unit{mm}) & \(\sigma\) (\unit{mm}) & \multicolumn{3}{c|}{Number of sources} \\ % & \multicolumn{3}{c|}{Filtered} \\
			\midrule
			5  & \num{10.3} & \num{5.3} & \hmcell{29 535}{\maxval} & \hmcell{29 710}{\maxval} & \hmcell{28 646}{\maxval} \\ % & \hmcell{29 190}{\maxval} & \hmcell{29 339}{\maxval} & \hmcell{28 280}{\maxval}  \\
			10 & \num{10.4} & \num{5.4} & \hmcell{21 029}{\maxval} & \hmcell{20 965}{\maxval} & \hmcell{20 154}{\maxval} \\ % & \hmcell{20 883}{\maxval} & \hmcell{20 760}{\maxval} & \hmcell{19 970}{\maxval}  \\
			15 & \num{10.3} & \num{4.6} & \hmcell{12 079}{\maxval} & \hmcell{12 104}{\maxval} & \hmcell{12 050}{\maxval} \\ % & \hmcell{12 046}{\maxval} & \hmcell{12 174}{\maxval} & \hmcell{11 988}{\maxval}  \\
			20 & \num{10.6} & \num{4.1} & \hmcell{05 507}{\maxval} & \hmcell{05 552}{\maxval} & \hmcell{05 936}{\maxval} \\ % & \hmcell{05 472}{\maxval} & \hmcell{05 522}{\maxval} & \hmcell{05 937}{\maxval}  \\
			25 & \num{10.2} & \num{3.7} & \hmcell{03 039}{\maxval} & \hmcell{03 120}{\maxval} & \hmcell{03 355}{\maxval} \\ % & \hmcell{03 024}{\maxval} & \hmcell{03 146}{\maxval} & \hmcell{03 332}{\maxval}  \\
			30 & \num{09.8} & \num{3.6} & \hmcell{02 494}{\maxval} & \hmcell{02 582}{\maxval} & \hmcell{02 652}{\maxval} \\ % & \hmcell{02 478}{\maxval} & \hmcell{02 594}{\maxval} & \hmcell{02 582}{\maxval}  \\
			\bottomrule
		\end{tabular}
	}
\end{table}

\begin{table}[!h]
	\def\maxval{30 459}%
	\caption{The number of outliers with \(\Delta > \mu + 2\sigma\), for a sample of \num{20} \DipoleScan\ reconstructions. The definitions of \(\mu\), \(\sigma\) and \(\peeld\) are the same as in \tabref{tab:sloreta-above-mu-plus-2sigma}.}
	\label{tab:dipolescan-above-mu-plus-2sigma}
	\resizebox{\linewidth}{!}{
		\begin{tabular}{|c|c|c|ccc|ccc|}
			\toprule
			\multicolumn{3}{|c|}{\(\peeld\) (\unit{mm})} & \num{0.0} & \num{0.5} & \num{1.0} \\ % & \num{0.0} & \num{0.5} & \num{1.0} \\
			\midrule
			\multicolumn{3}{|c|}{Number of columns in \(L\)} & \hmcell{30 288}{\maxval} & \hmcell{30 459}{\maxval} & \hmcell{29 367}{\maxval} \\ % & \hmcell{29 835}{\maxval} & \hmcell{30 003}{\maxval} & \hmcell{28 926}{\maxval} \\
			\midrule
			SNR (\unit{dB}) & \(\mu\) (\unit{mm}) & \(\sigma\) (\unit{mm}) & \multicolumn{3}{c|}{Number of sources} \\ % & \multicolumn{3}{c|}{Filtered} \\
			\midrule
			5  & \num{10.3} & \num{5.3} & \hmcell{19 546}{\maxval} & \hmcell{19 652}{\maxval} & \hmcell{18 923}{\maxval} \\ % & \hmcell{19 384}{\maxval} & \hmcell{19 424}{\maxval} & \hmcell{18 713}{\maxval}   \\
			10 & \num{10.4} & \num{5.4} & \hmcell{06 935}{\maxval} & \hmcell{06 975}{\maxval} & \hmcell{06 866}{\maxval} \\ % & \hmcell{06 897}{\maxval} & \hmcell{06 902}{\maxval} & \hmcell{06 840}{\maxval}   \\
			15 & \num{10.3} & \num{4.6} & \hmcell{01 893}{\maxval} & \hmcell{01 857}{\maxval} & \hmcell{01 836}{\maxval} \\ % & \hmcell{01 879}{\maxval} & \hmcell{01 859}{\maxval} & \hmcell{01 859}{\maxval}   \\
			20 & \num{10.6} & \num{4.1} & \hmcell{00 429}{\maxval} & \hmcell{00 411}{\maxval} & \hmcell{00 411}{\maxval} \\ % & \hmcell{00 430}{\maxval} & \hmcell{00 421}{\maxval} & \hmcell{00 405}{\maxval}   \\
			25 & \num{10.2} & \num{3.7} & \hmcell{00 036}{\maxval} & \hmcell{00 031}{\maxval} & \hmcell{00 024}{\maxval} \\ % & \hmcell{00 037}{\maxval} & \hmcell{00 030}{\maxval} & \hmcell{00 025}{\maxval}   \\
			30 & \num{09.8} & \num{3.6} & \hmcell{00 001}{\maxval} & \hmcell{00 000}{\maxval} & \hmcell{00 000}{\maxval} \\ % & \hmcell{00 000}{\maxval} & \hmcell{00 000}{\maxval} & \hmcell{00 000}{\maxval}   \\
			\bottomrule
		\end{tabular}
	}
\end{table}

For \sLORETA\ we observed the following: the more tetra are peeled, the more
outliers occur at noise levels below \SI{15}{\decibel}. On the other hand,
\DipoleScan\ displays a rather consistent outcome at and below
\SI{15}{\decibel} noise level, with increased peeling reducing the numbers of
outliers. In either case, it seems that noise might be the major contributing
factor below the \SI{15}{\decibel} mark, with \DipoleScan \ being slightly
more resistant to noise effects.

To get a sense of where the statistical outliers are located,
Figures~\ref{fig:sloreta-nofilter-localisation-error-distribution} and
\ref{fig:dipolescan-nofilter-localisation-error-distribution} display the
localisation error \(\Delta\) of \sLORETA\ and \DipoleScan\ as a function of
position, in the active gray matter compartment. Colors in the figures have
been adjusted, such that the obvious outlier positions with \(\Delta\geq\mu +
2\sigma\) are highlighted in red. The displayed noise levels are the ones
where peeling was discovered to have an observable reduction in \(\Delta\).

\def\subfigwidth{0.31\linewidth}

\def\cbarwidth{0.04\linewidth}

\def\cbvspace{-0.6cm}

\begin{figure}[h!]

	\begin{minipage}[t]{0.9\linewidth}
		\begin{subfigure}[\subfigwidth]{\(d\subt p = \SI{0}{\milli\meter}\) at \SI{25}{\decibel}.}{%
			fig:sloreta-localisation-error-distribution-d0mm-25db%
		}
			\includegraphics[width=\linewidth]{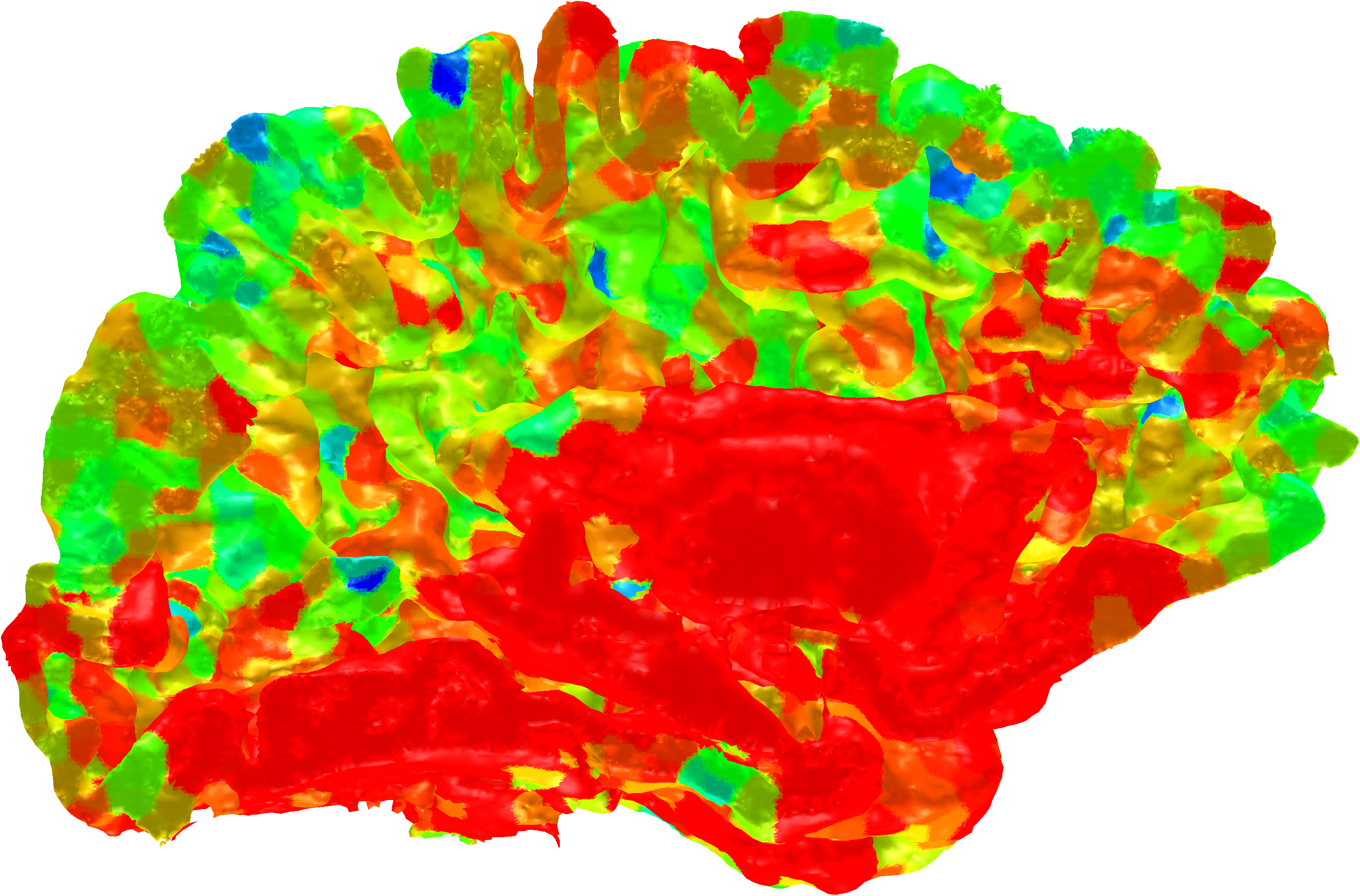}
		\end{subfigure}
		\hfill
		\begin{subfigure}[\subfigwidth]{\(d\subt p = \SI{0.5}{\milli\meter}\) at \SI{25}{\decibel}.}{%
			fig:sloreta-localisation-error-distribution-d05mm-25db%
		}
			\includegraphics[width=\linewidth]{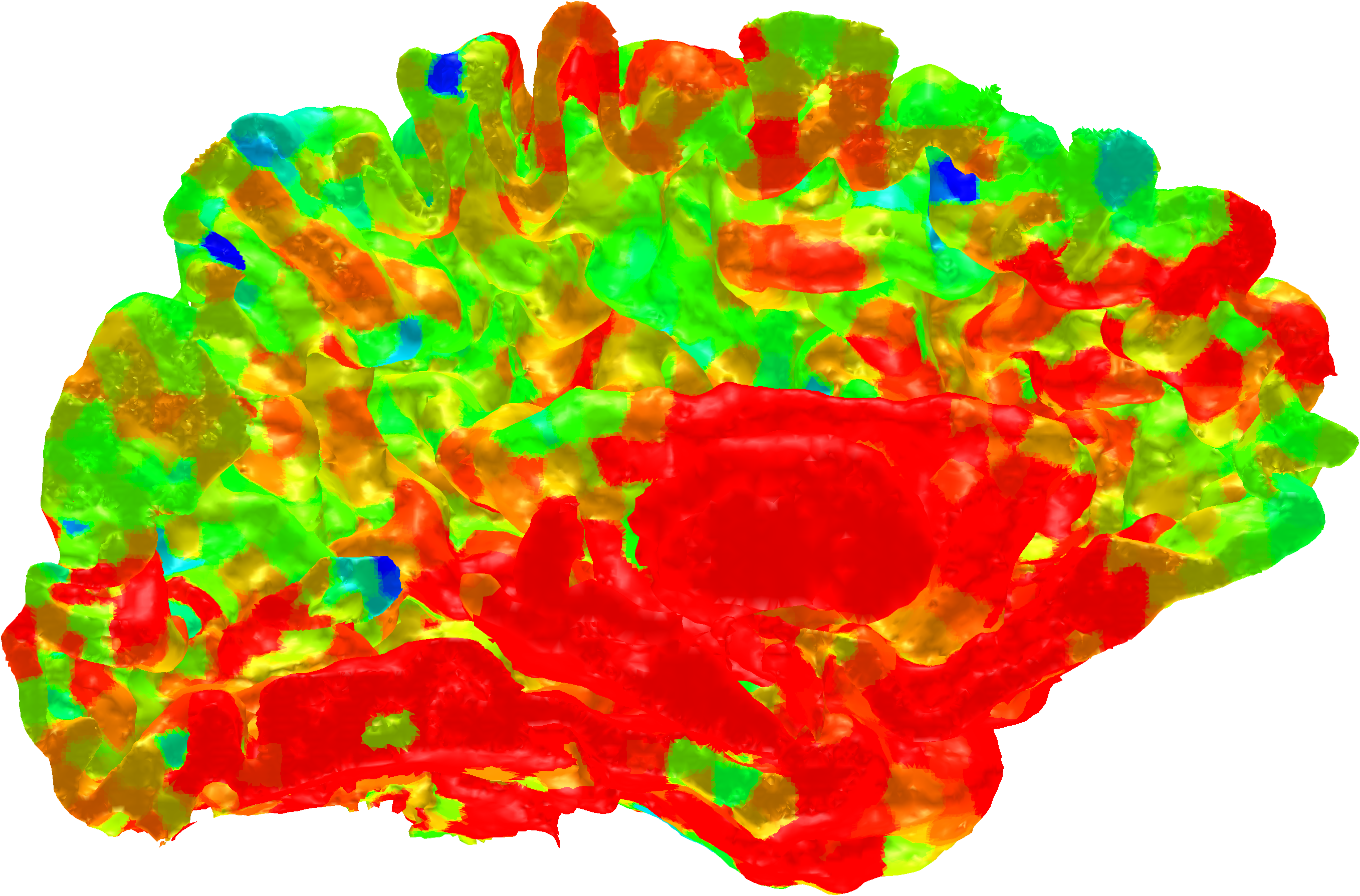}
		\end{subfigure}
		\hfill
		\begin{subfigure}[\subfigwidth]{\(d\subt p = \SI{1}{\milli\meter}\) at \SI{25}{\decibel}.}{%
			fig:sloreta-localisation-error-distribution-d1mm-25db%
		}
			\includegraphics[width=\linewidth]{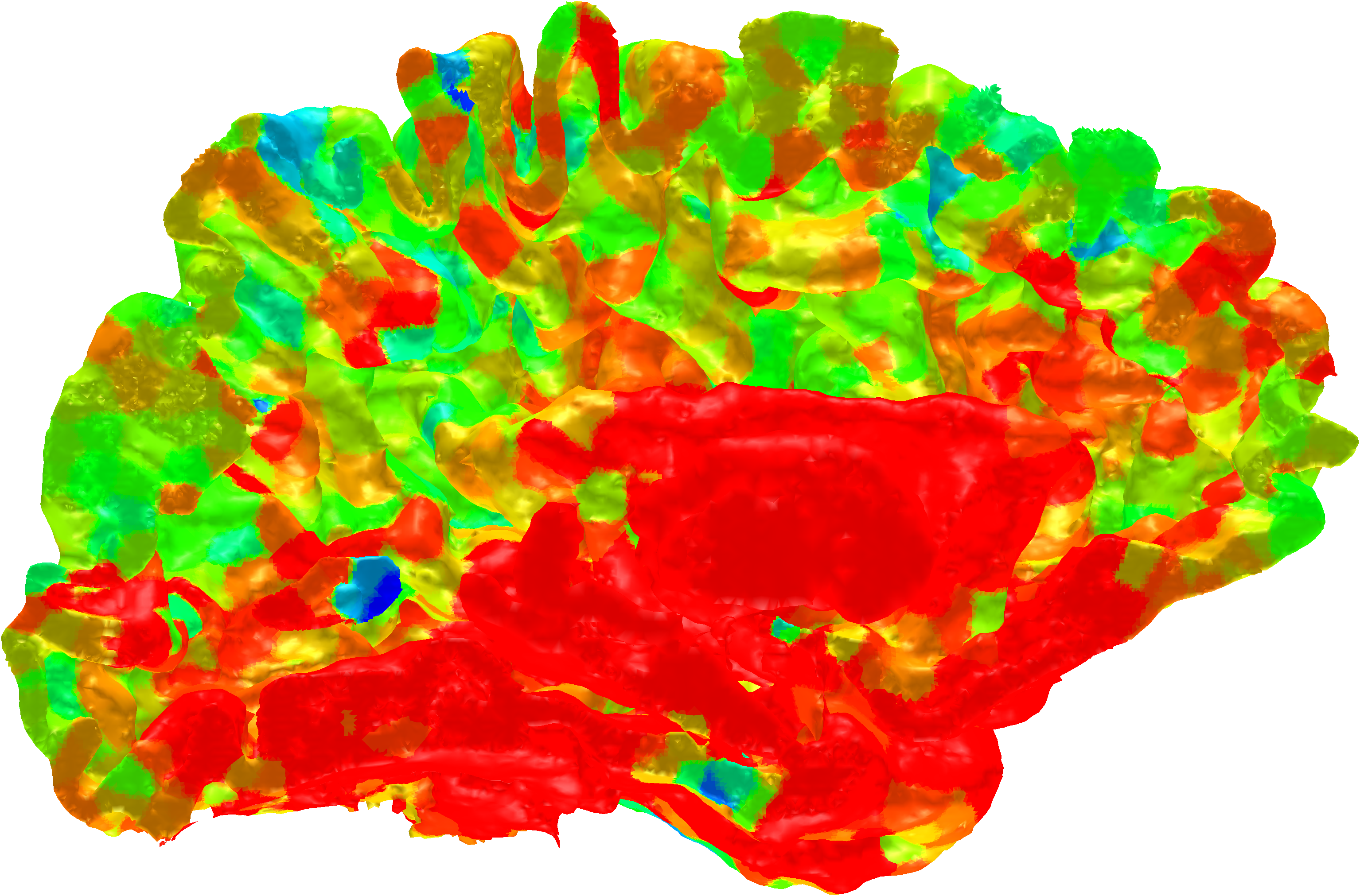}
		\end{subfigure}
	\end{minipage}
	\hfill
	\begin{minipage}[t]{\cbarwidth}
		\centering%
		\vspace*{\cbvspace}%
		\includegraphics[width=\linewidth]{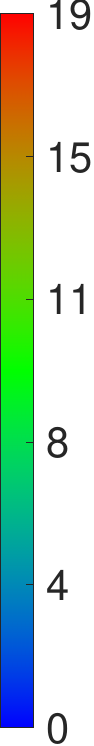}
	\end{minipage}

	\bigskip

	\begin{minipage}[t]{0.9\linewidth}
		\begin{subfigure}[\subfigwidth]{\(d\subt p = \SI{0}{\milli\meter}\) at \SI{30}{\decibel}.}{%
			fig:sloreta-localisation-error-distribution-d0mm-30db%
		}
			\includegraphics[width=\linewidth]{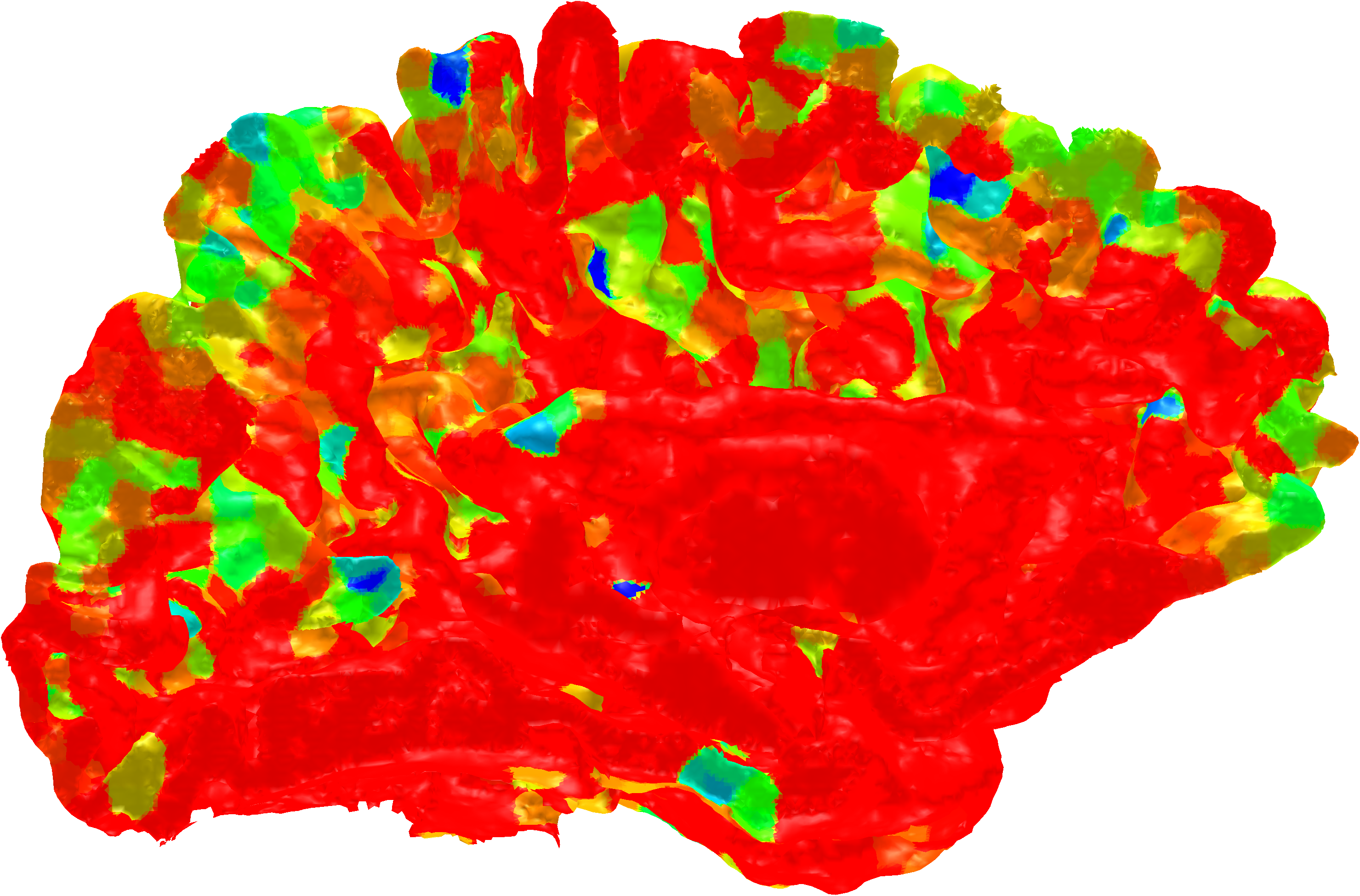}
		\end{subfigure}
		\hfill
		\begin{subfigure}[\subfigwidth]{\(d\subt p = \SI{0.5}{\milli\meter}\) at \SI{30}{\decibel}.}{%
			fig:sloreta-localisation-error-distribution-d05mm-30db%
		}
			\includegraphics[width=\linewidth]{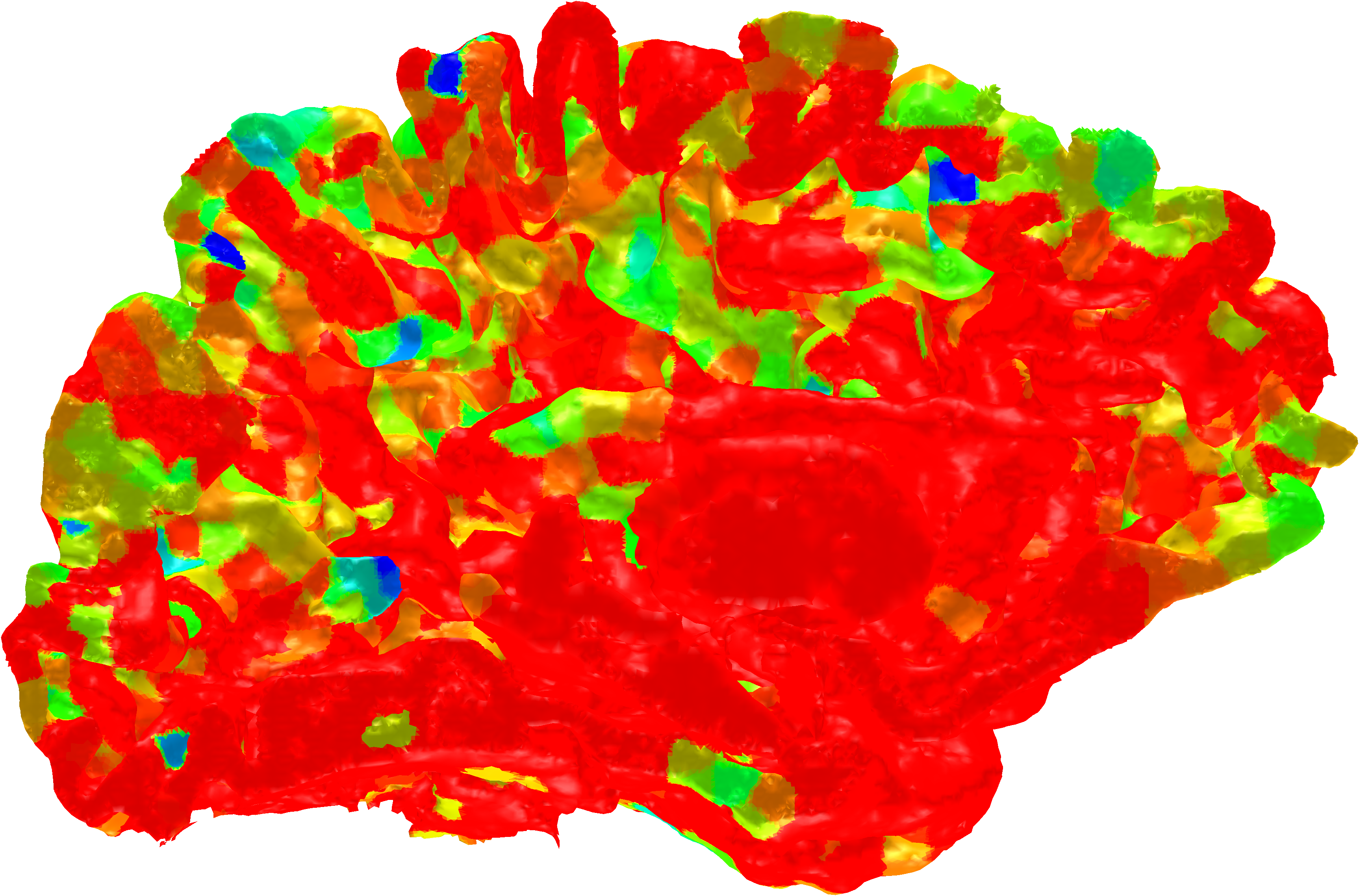}
		\end{subfigure}
		\hfill
		\begin{subfigure}[\subfigwidth]{\(d\subt p = \SI{1}{\milli\meter}\) at \SI{30}{\decibel}.}{%
			fig:sloreta-localisation-error-distribution-d1mm-30db%
		}
			\includegraphics[width=\linewidth]{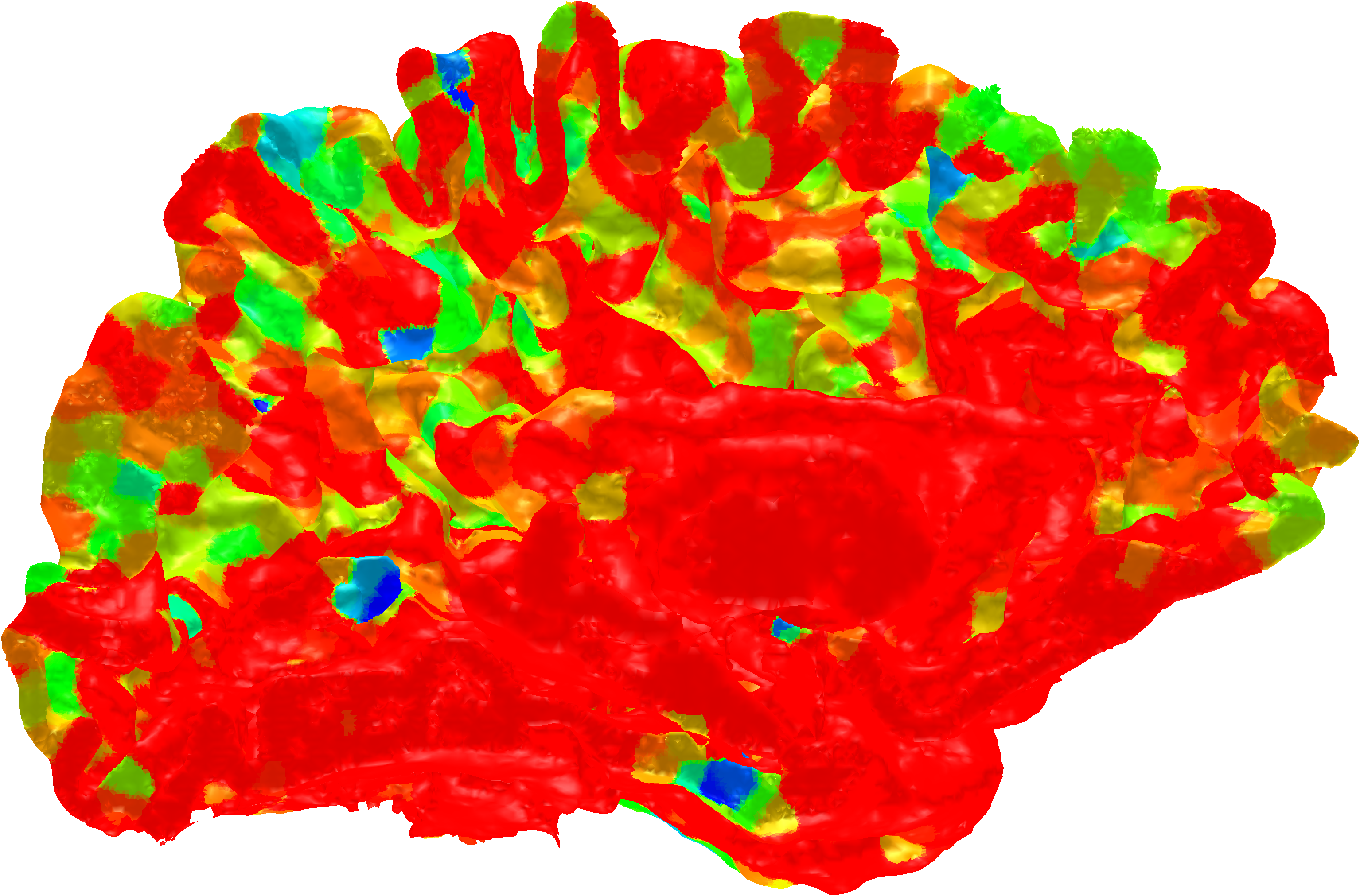}
		\end{subfigure}
	\end{minipage}
	\hfill
	\begin{minipage}[t]{\cbarwidth}
		\centering%
		\vspace*{\cbvspace}%
		\includegraphics[width=\linewidth]{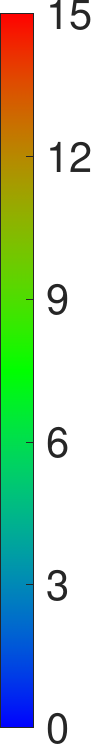}
	\end{minipage}

	\caption{%
		Sagittal views of localisation error \(\Delta = \Delta\arcs{\vec x}\)
		(\si{\milli\meter}) in the active gray matter layer in the case of
		\sLORETA, contrasted against the different peeling depths \(\peeld\)
		of \tabref{tab:sloreta-unfiltered-inverse-statistics}. The red color
		indicates where \(\Delta\geq\mu+2\sigma\), as in outliers in the
		entire distribution.%
	}
	\label{fig:sloreta-nofilter-localisation-error-distribution}
\end{figure}

\def\cbvspace{\vspace*{-0.2cm}}

\def\cbarwidth{0.03\linewidth}

\begin{figure}[h!]

	\begin{minipage}[t]{0.9\linewidth}
		\begin{subfigure}[\subfigwidth]{\(d\subt p = \SI{0}{\milli\meter}\) at \SI{25}{\decibel}.}{%
			fig:dipolescan-localisation-error-distribution-d0mm-25db%
		}
			\includegraphics[width=\linewidth]{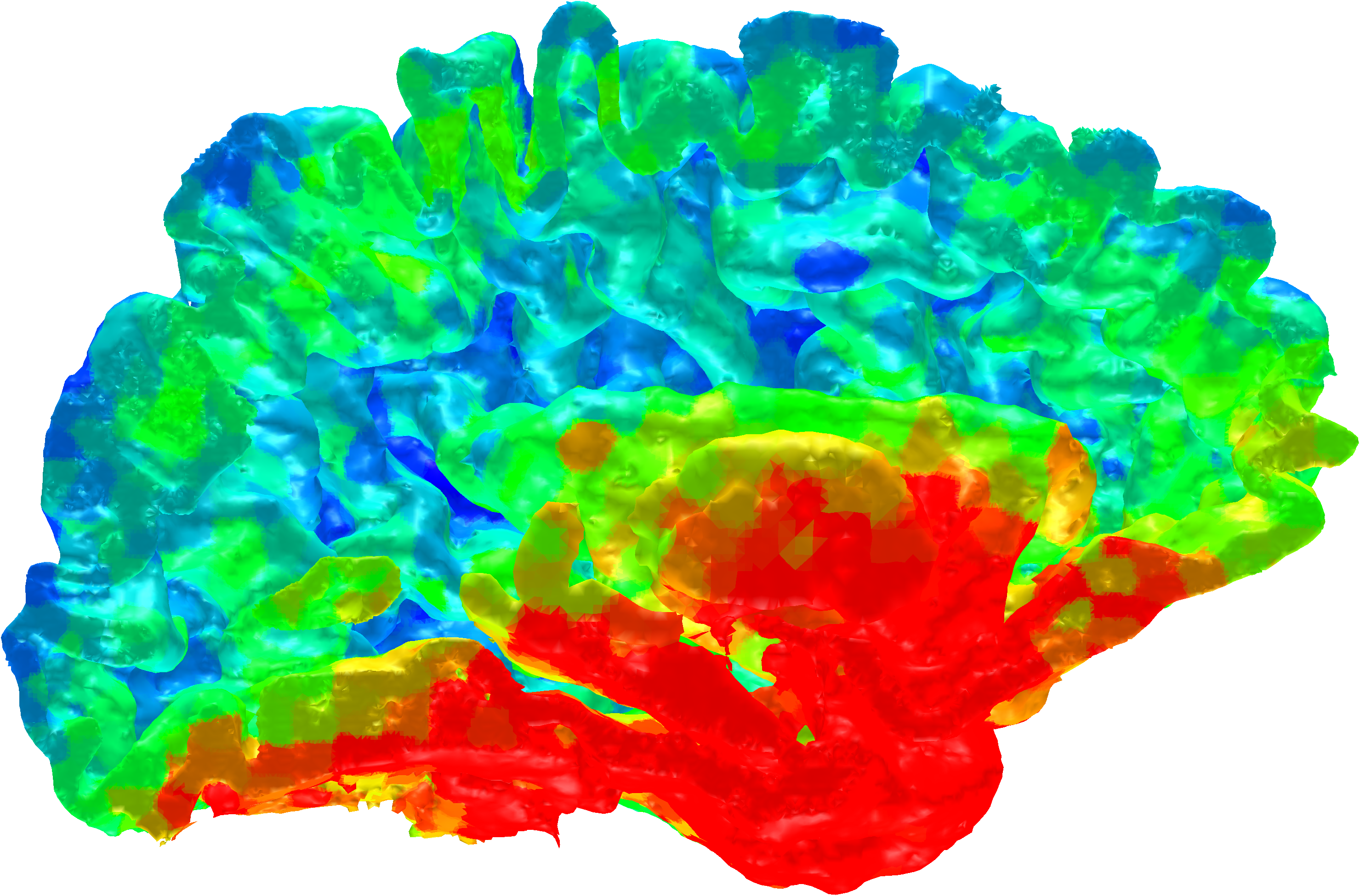}
		\end{subfigure}
		\hfill
		\begin{subfigure}[\subfigwidth]{\(d\subt p = \SI{0.5}{\milli\meter}\) at \SI{25}{\decibel}.}{%
			fig:dipolescan-localisation-error-distribution-d05mm-25db%
		}
			\includegraphics[width=\linewidth]{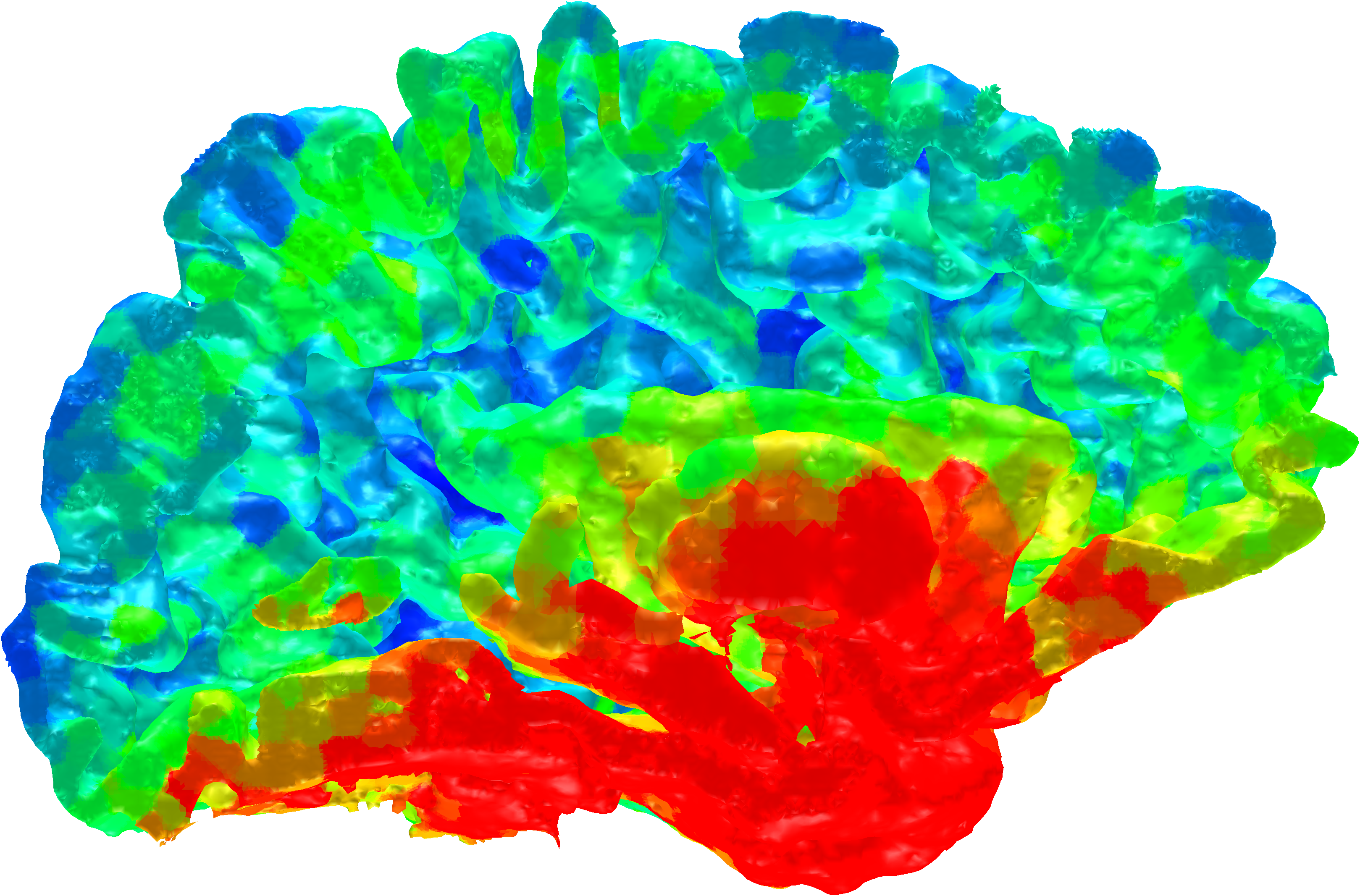}
		\end{subfigure}
		\hfill
		\begin{subfigure}[\subfigwidth]{\(d\subt p = \SI{1}{\milli\meter}\) at \SI{25}{\decibel}.}{%
			fig:dipolescan-localisation-error-distribution-d1mm-25db%
		}
			\includegraphics[width=\linewidth]{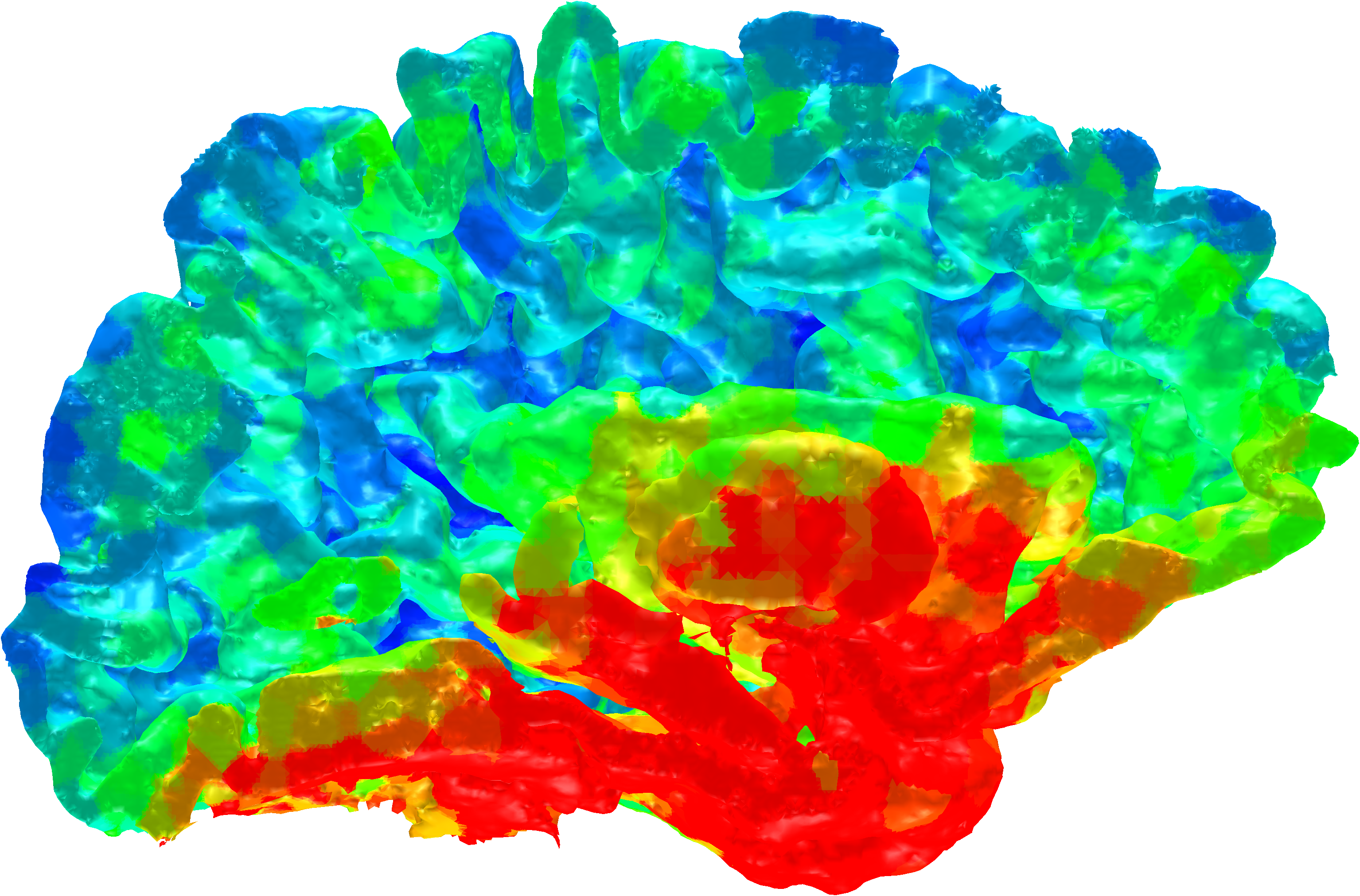}
		\end{subfigure}
	\end{minipage}
	\hfill
	\begin{minipage}[t]{\cbarwidth}
		\centering%
		\cbvspace%
		\includegraphics[width=\linewidth]{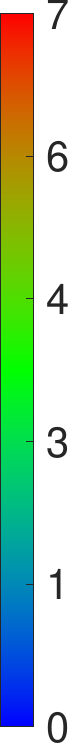}
	\end{minipage}

	\bigskip

	\begin{minipage}[t]{0.9\linewidth}
		\begin{subfigure}[\subfigwidth]{\(d\subt p = \SI{0}{\milli\meter}\) at \SI{30}{\decibel}.}{%
			fig:dipolescan-localisation-error-distribution-d0mm-30db%
		}
			\includegraphics[width=\linewidth]{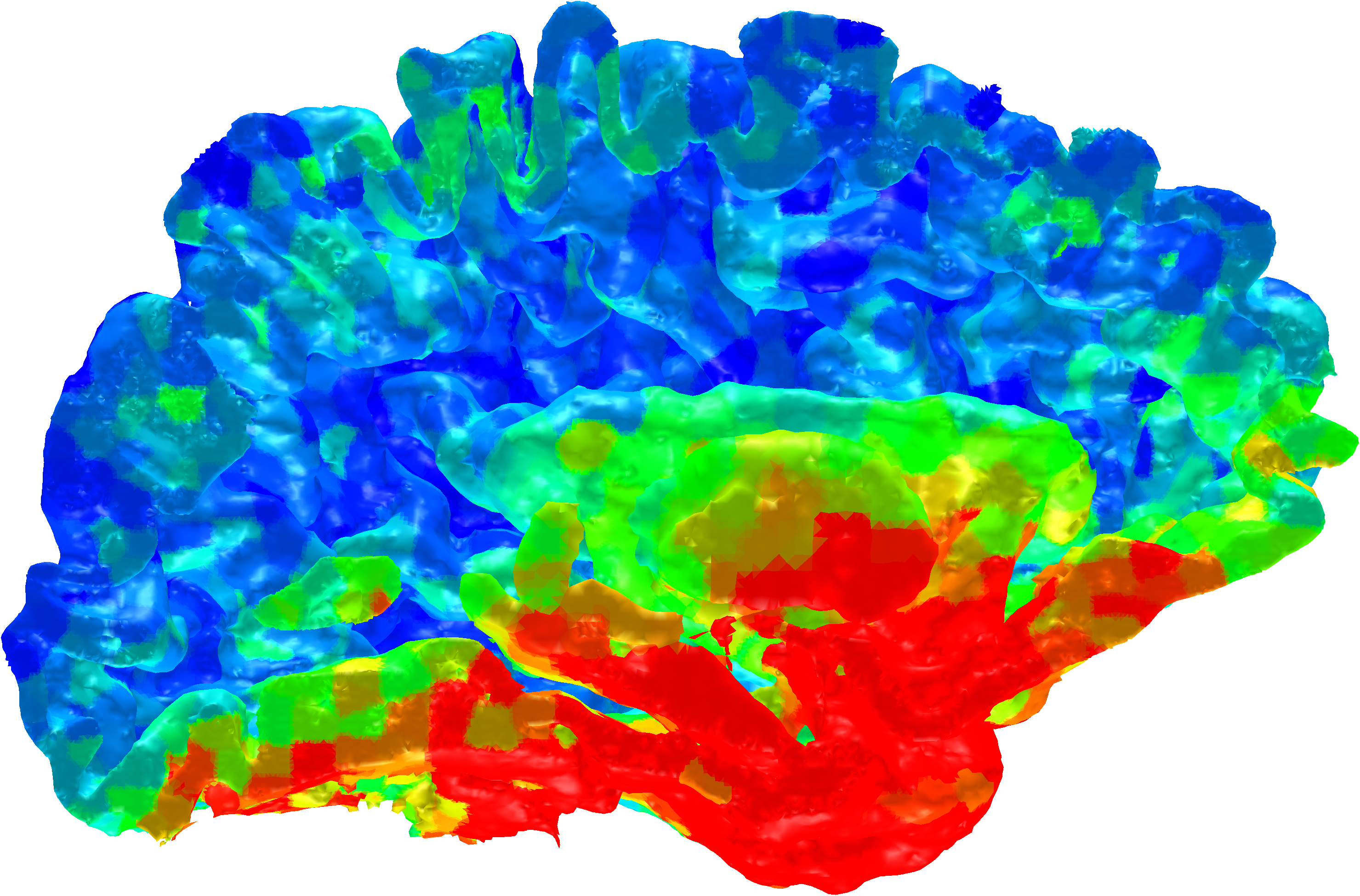}
		\end{subfigure}
		\hfill
		\begin{subfigure}[\subfigwidth]{\(d\subt p = \SI{0.5}{\milli\meter}\) at \SI{30}{\decibel}.}{%
			fig:dipolescan-localisation-error-distribution-d05mm-30db%
		}
			\includegraphics[width=\linewidth]{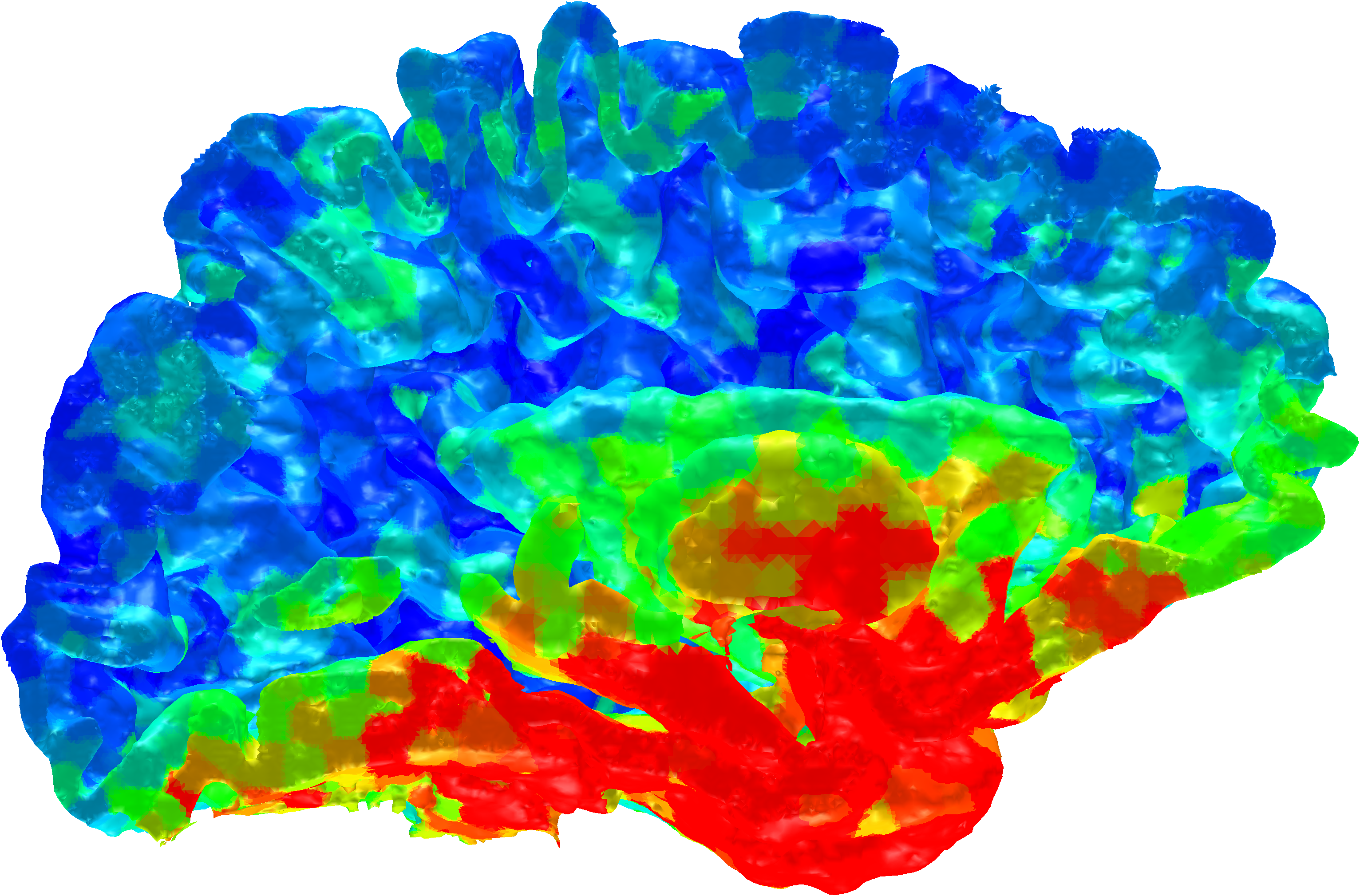}
		\end{subfigure}
		\hfill
		\begin{subfigure}[\subfigwidth]{\(d\subt p = \SI{1}{\milli\meter}\) at \SI{30}{\decibel}.}{%
			fig:dipolescan-localisation-error-distribution-d1mm-30db%
		}
			\includegraphics[width=\linewidth]{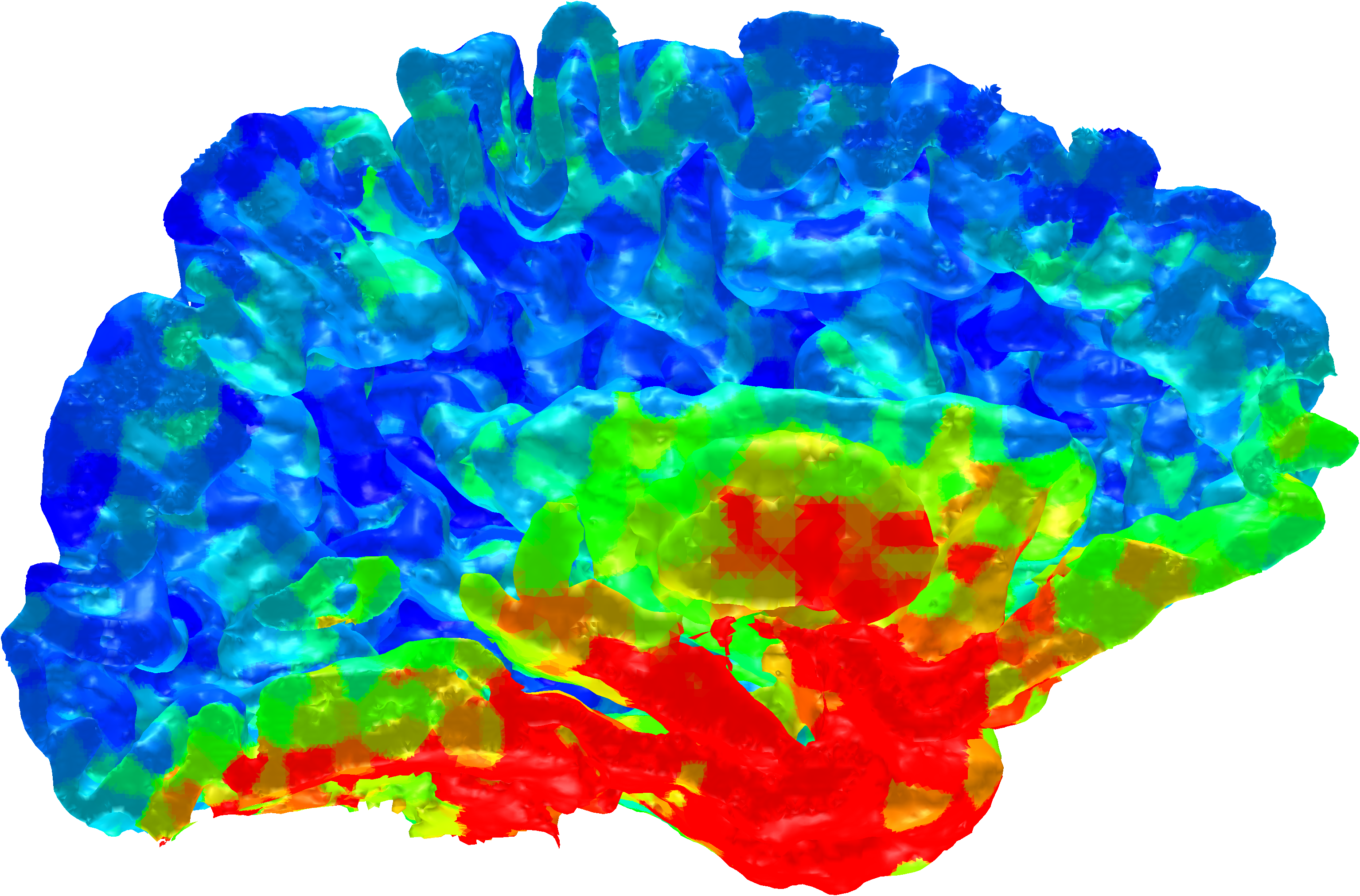}
		\end{subfigure}
	\end{minipage}
	\hfill
	\begin{minipage}[t]{\cbarwidth}
		\centering%
		\cbvspace%
		\includegraphics[width=\linewidth]{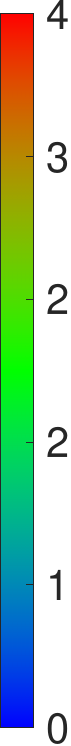}
	\end{minipage}

	\caption{%
		Sagittal views of localisation error \(\Delta = \Delta\arcs{\vec x}\)
		(\si{\milli\meter}) in the active gray matter layer in the case of
		\DipoleScan, contrasted against the different peeling depths
		\(\peeld\) of \tabref{tab:dipolescan-unfiltered-inverse-statistics}.
		The meaning of the colors is the same as in
		\figref{fig:sloreta-nofilter-localisation-error-distribution}.%
	}
	\label{fig:dipolescan-nofilter-localisation-error-distribution}
\end{figure}

To illustrate how the localisation error fluctuates locally at lower noise
levels, both deep in the brain and more superficially,
Figures~\ref{fig:focused-thalamus-sloreta-nofilter-localisation-error-distribution}
and \ref{fig:focused-thalamus-dipolescan-nofilter-localisation-error-distribution}
display thalamically focused views of the above images, whereas
Figures~\ref{fig:focused-cerebral-sloreta-nofilter-localisation-error-distribution}
and \ref{fig:focused-cerebral-dipolescan-nofilter-localisation-error-distribution} do
the same for the cortex. The displayed noise levels are chosen to be the ones, where
visible improvements can be seen.

\edef\thalamictrimsettings{trim={150px 0 80px 150px}, clip}

\edef\parietaltrimsettings{trim={50px 200px 200px 0px}, clip}

\def\cbarwidth{0.03\linewidth}

\def\cbvspace{-0.2cm}

\begin{figure}[h!]

	\begin{minipage}[t]{0.9\linewidth}
		\begin{subfigure}[\subfigwidth]{\(d\subt p = \SI{0}{\milli\meter}\) at \SI{30}{\decibel}.}{%
			fig:parietal-focused-sloreta-localisation-error-distribution-d0mm-30db%
		}
			\expandafter\includegraphics\expandafter[\parietaltrimsettings, width=\linewidth]{images/sloreta_d0mm_nofilter_30db_dist_err.png}
		\end{subfigure}
		\hfill
		\begin{subfigure}[\subfigwidth]{\(d\subt p = \SI{0.5}{\milli\meter}\) at \SI{30}{\decibel}.}{%
			fig:parietal-focused-sloreta-localisation-error-distribution-d05mm-30db%
		}
			\expandafter\includegraphics\expandafter[\parietaltrimsettings, width=\linewidth]{images/sloreta_d05mm_nofilter_30db_dist_err.png}
		\end{subfigure}
		\hfill
		\begin{subfigure}[\subfigwidth]{\(d\subt p = \SI{1}{\milli\meter}\) at \SI{30}{\decibel}.}{%
			fig:parietal-focused-sloreta-localisation-error-distribution-d1mm-30db%
		}
			\expandafter\includegraphics\expandafter[\parietaltrimsettings, width=\linewidth]{images/sloreta_d1mm_nofilter_30db_dist_err.png}
		\end{subfigure}
	\end{minipage}
	\hfill
	\begin{minipage}[t]{\cbarwidth}%
		\centering%
		\vspace*{\cbvspace}%
		\includegraphics[width=\linewidth]{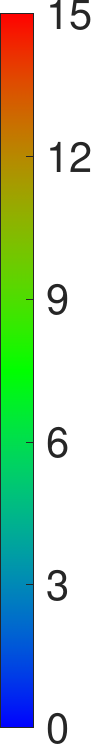}%
	\end{minipage}

	\caption{Focused views of \(\Delta = \Delta\arcs{\vec x}\) in the parietal region, in the case of \sLORETA.}

	\label{fig:focused-cerebral-sloreta-nofilter-localisation-error-distribution}

\end{figure}

\begin{figure}[h!]

	\begin{minipage}[t]{0.9\linewidth}
		\begin{subfigure}[\subfigwidth]{\(d\subt p = \SI{0}{\milli\meter}\) at \SI{30}{\decibel}.}{%
			fig:thalamic-focused-sloreta-localisation-error-distribution-d0mm-30db%
		}
			\expandafter\includegraphics\expandafter[\thalamictrimsettings, width=\linewidth]{images/sloreta_d0mm_nofilter_30db_dist_err.png}
		\end{subfigure}
		\hfill
		\begin{subfigure}[\subfigwidth]{\(d\subt p = \SI{0.5}{\milli\meter}\) at \SI{30}{\decibel}.}{%
			fig:thalamic-focused-sloreta-localisation-error-distribution-d05mm-30db%
		}
			\expandafter\includegraphics\expandafter[\thalamictrimsettings, width=\linewidth]{images/sloreta_d05mm_nofilter_30db_dist_err.png}
		\end{subfigure}
		\hfill
		\begin{subfigure}[\subfigwidth]{\(d\subt p = \SI{1}{\milli\meter}\) at \SI{30}{\decibel}.}{%
			fig:thalamic-focused-sloreta-localisation-error-distribution-d1mm-30db%
		}
			\expandafter\includegraphics\expandafter[\thalamictrimsettings, width=\linewidth]{images/sloreta_d1mm_nofilter_30db_dist_err.png}
		\end{subfigure}
	\end{minipage}
	\hfill
	\begin{minipage}[t]{\cbarwidth}%
		\centering%
		\vspace*{\cbvspace}%
		\includegraphics[width=\linewidth]{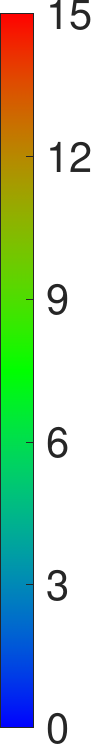}%
	\end{minipage}

	\caption{Focused views of \(\Delta = \Delta\arcs{\vec x}\) in the thalamic region, in the case of \sLORETA.}

	\label{fig:focused-thalamus-sloreta-nofilter-localisation-error-distribution}

\end{figure}

\def\cbarwidth{0.025\linewidth}

\begin{figure}[h!]

	\begin{minipage}[t]{0.9\linewidth}
		\begin{subfigure}[\subfigwidth]{\(d\subt p = \SI{0}{\milli\meter}\) at \SI{30}{\decibel}.}{%
			fig:parietal-focused-dipolescan-localisation-error-distribution-d0mm-30db%
		}
			\expandafter\includegraphics\expandafter[\parietaltrimsettings, width=\linewidth]{images/dipolescan_d0mm_nofilter_30db_dist_err.png}
		\end{subfigure}
		\hfill
		\begin{subfigure}[\subfigwidth]{\(d\subt p = \SI{0.5}{\milli\meter}\) at \SI{30}{\decibel}.}{%
			fig:parietal-focused-dipolescan-localisation-error-distribution-d05mm-30db%
		}
			\expandafter\includegraphics\expandafter[\parietaltrimsettings, width=\linewidth]{images/dipolescan_d05mm_nofilter_30db_dist_err.png}
		\end{subfigure}
		\hfill
		\begin{subfigure}[\subfigwidth]{\(d\subt p = \SI{1}{\milli\meter}\) at \SI{30}{\decibel}.}{%
			fig:parietal-focused-dipolescan-localisation-error-distribution-d1mm-30db%
		}
			\expandafter\includegraphics\expandafter[\parietaltrimsettings, width=\linewidth]{images/dipolescan_d1mm_nofilter_30db_dist_err.png}
		\end{subfigure}
	\end{minipage}
	\hfill
	\begin{minipage}[t]{\cbarwidth}%
		\centering%
		\vspace*{\cbvspace}%
		\includegraphics[width=\linewidth]{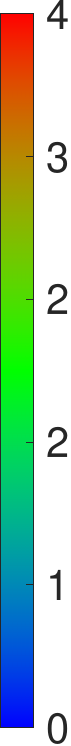}%
	\end{minipage}

	\caption{Focused views of \(\Delta = \Delta\arcs{\vec x}\) in the parietal region, in the case of \DipoleScan.}

	\label{fig:focused-cerebral-dipolescan-nofilter-localisation-error-distribution}

\end{figure}

\begin{figure}[h!]

	\begin{minipage}[t]{0.9\linewidth}
		\begin{subfigure}[\subfigwidth]{\(d\subt p = \SI{0}{\milli\meter}\) at \SI{30}{\decibel}.}{%
			fig:thalamic-focused-dipolescan-localisation-error-distribution-d0mm-30db%
		}
			\expandafter\includegraphics\expandafter[\thalamictrimsettings, width=\linewidth]{images/dipolescan_d0mm_nofilter_30db_dist_err.png}
		\end{subfigure}
		\hfill
		\begin{subfigure}[\subfigwidth]{\(d\subt p = \SI{0.5}{\milli\meter}\) at \SI{30}{\decibel}.}{%
			fig:thalamic-focused-dipolescan-localisation-error-distribution-d05mm-30db%
		}
			\expandafter\includegraphics\expandafter[\thalamictrimsettings, width=\linewidth]{images/dipolescan_d05mm_nofilter_30db_dist_err.png}
		\end{subfigure}
		\hfill
		\begin{subfigure}[\subfigwidth]{\(d\subt p = \SI{1}{\milli\meter}\) at \SI{30}{\decibel}.}{%
			fig:thalamic-focused-dipolescan-localisation-error-distribution-d1mm-30db%
		}
			\expandafter\includegraphics\expandafter[\thalamictrimsettings, width=\linewidth]{images/dipolescan_d1mm_nofilter_30db_dist_err.png}
		\end{subfigure}
	\end{minipage}
	\hfill
	\begin{minipage}[t]{\cbarwidth}%
		\centering%
		\vspace*{\cbvspace}%
		\includegraphics[width=\linewidth]{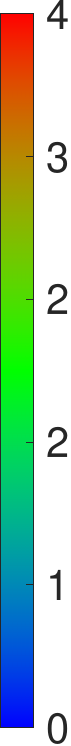}%
	\end{minipage}

	\caption{Focused views of \(\Delta = \Delta\arcs{\vec x}\) in the thalamic region, in the case of \DipoleScan.}

	\label{fig:focused-thalamus-dipolescan-nofilter-localisation-error-distribution}

\end{figure}

These mappings demonstrate that peeling can help in reducing localisation error
universally in the active domain, when the noise level is low enough: for example at
\SI{30}{\decibel} \SNR\ in
\figref{fig:focused-thalamus-dipolescan-nofilter-localisation-error-distribution},
the localisation error in the part of the frontal lobe in front of the thalamus
decreases, as \(d\subt p\) is increased. A similar effect can be seen in
\figref{fig:focused-cerebral-sloreta-nofilter-localisation-error-distribution}, where
localisation error in the cerebrum, in front of the central sulcus decreases.

To observe how the peeling of the active layers affects the accuracy of localisation
with \sLORETA\ and \DipoleScan\ locally,
Figures~\ref{fig:tan-normal-amplitude-dipolescan-single-source}--\ref{fig:tan-normal-sloreta-single-source}
display relative strenghts of dipole moments of the reconstructed dipolar
distribution \(\rec\Jp = \rec\Jp\arcs{\vec x}\) for a single dipole near Brodmann area
3b, around the central sulcus. Both tangential and radial sources are observed, as
these correspond to the ideal cases of \MEG\ and \EEG, respectively.

\def\subfigwidth{0.30\linewidth}

\begin{figure}[h!]
	\begin{minipage}[t]{0.8\columnwidth}
		\begin{subfigure}[\subfigwidth]{}{%
			fig:tan-amplitude-single-source-dipolescan-d0mm%
		}
			\includegraphics[width=\linewidth]{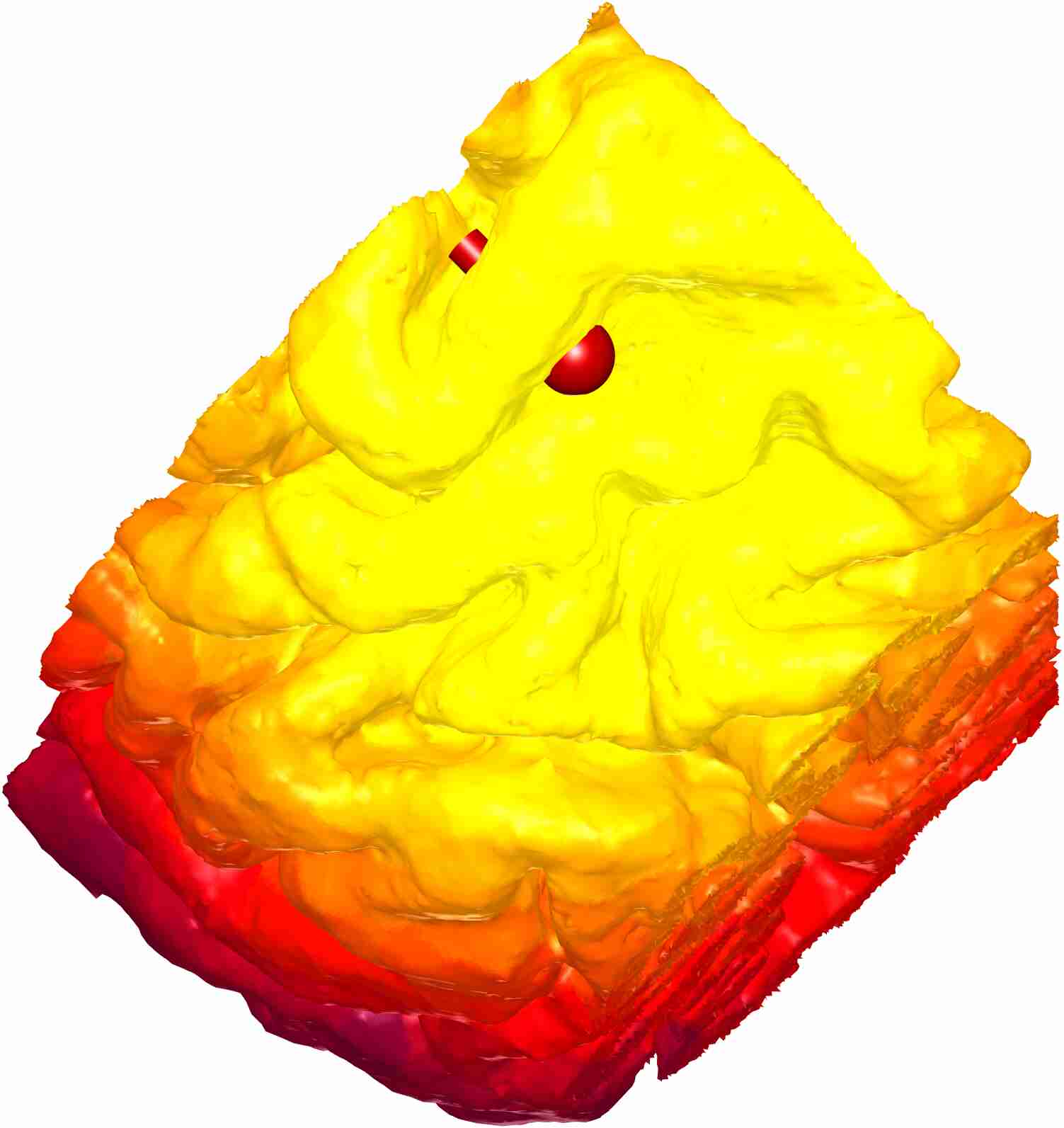}
		\end{subfigure}
		\hfill
		\begin{subfigure}[\subfigwidth]{}{%
			fig:tan-amplitude-single-source-dipolescan-d05mm%
		}
			\includegraphics[width=\linewidth]{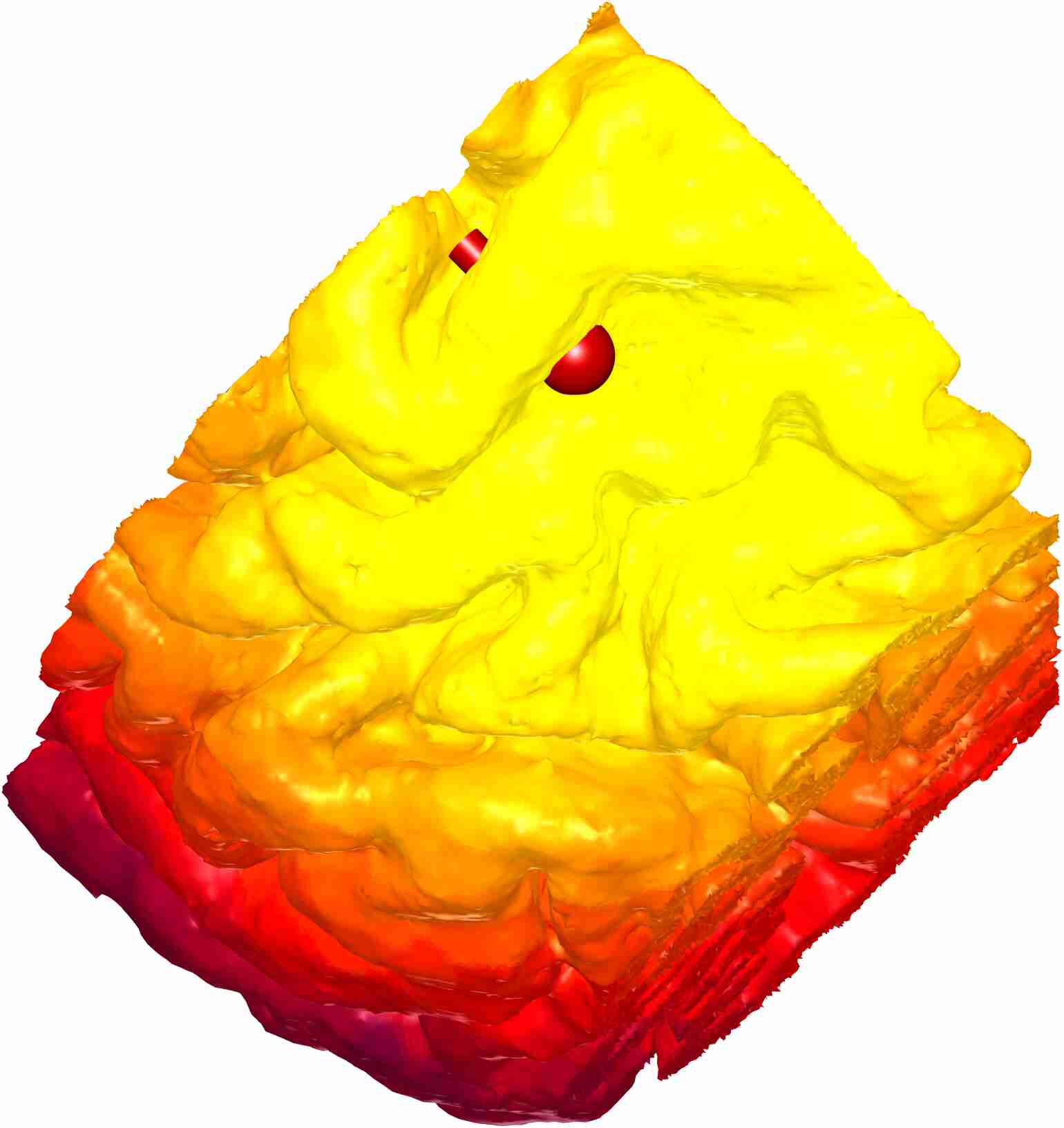}
		\end{subfigure}
		\hfill
		\begin{subfigure}[\subfigwidth]{}{%
			fig:tan-amplitude-single-source-dipolescan-d1mm%
		}
			\includegraphics[width=\linewidth]{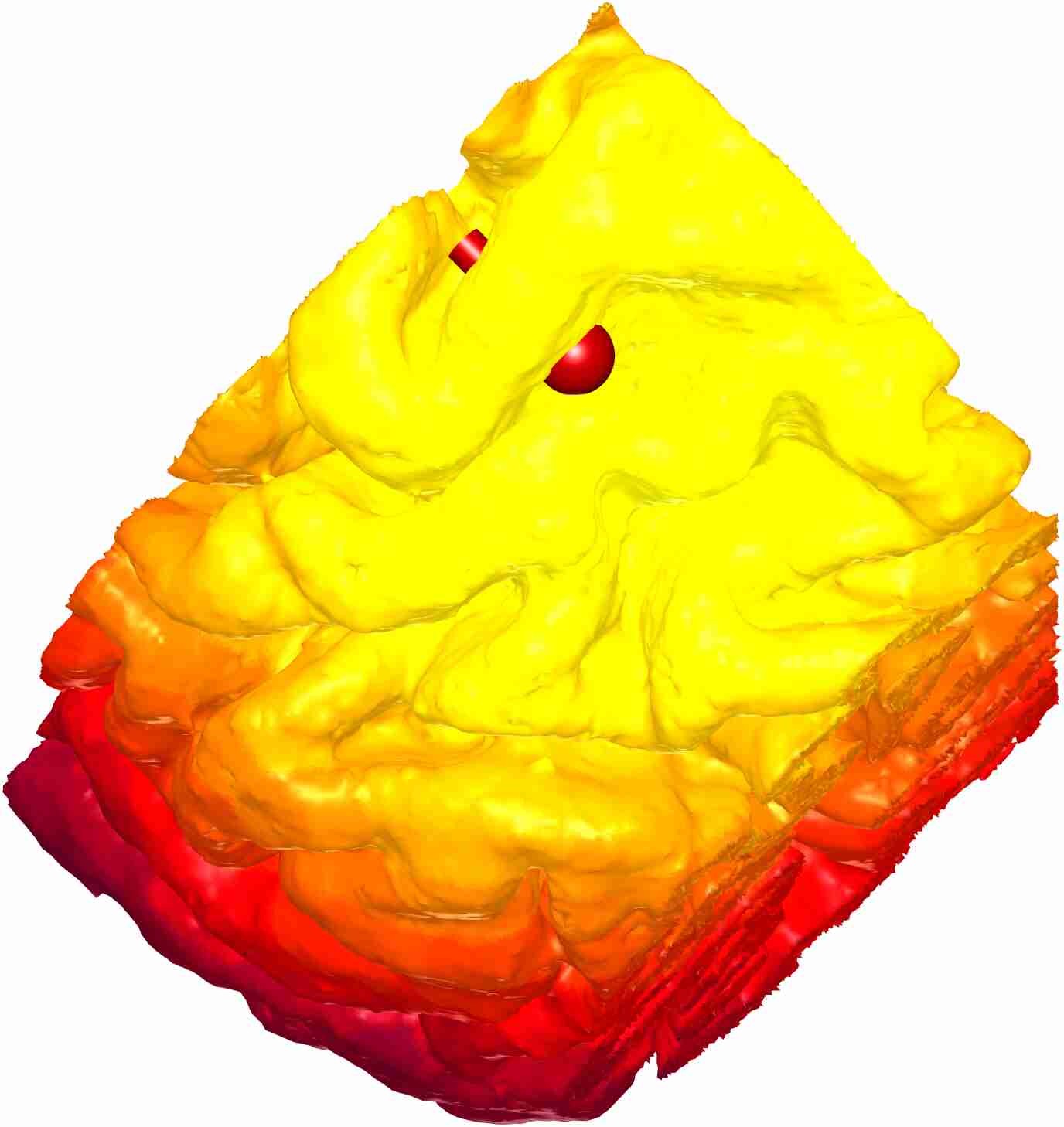}
		\end{subfigure}

		\bigskip

		\begin{subfigure}[\subfigwidth]{}{%
			fig:normal-amplitude-single-source-dipolescan-d0mm%
		}
			\includegraphics[width=\linewidth]{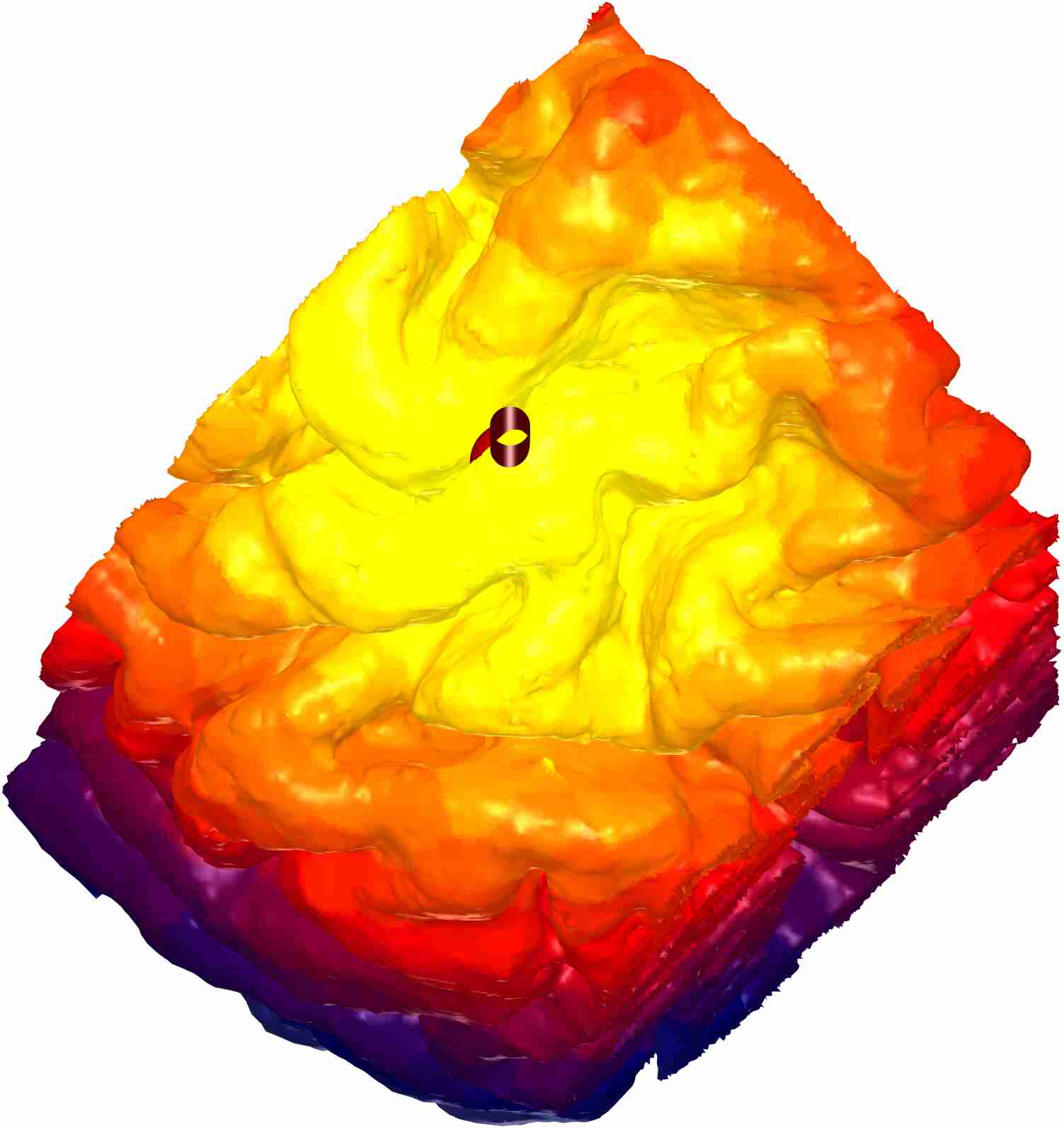}
		\end{subfigure}
		\hfill
		\begin{subfigure}[\subfigwidth]{}{%
			fig:normal-amplitude-single-source-dipolescan-d05mm%
		}
			\includegraphics[width=\linewidth]{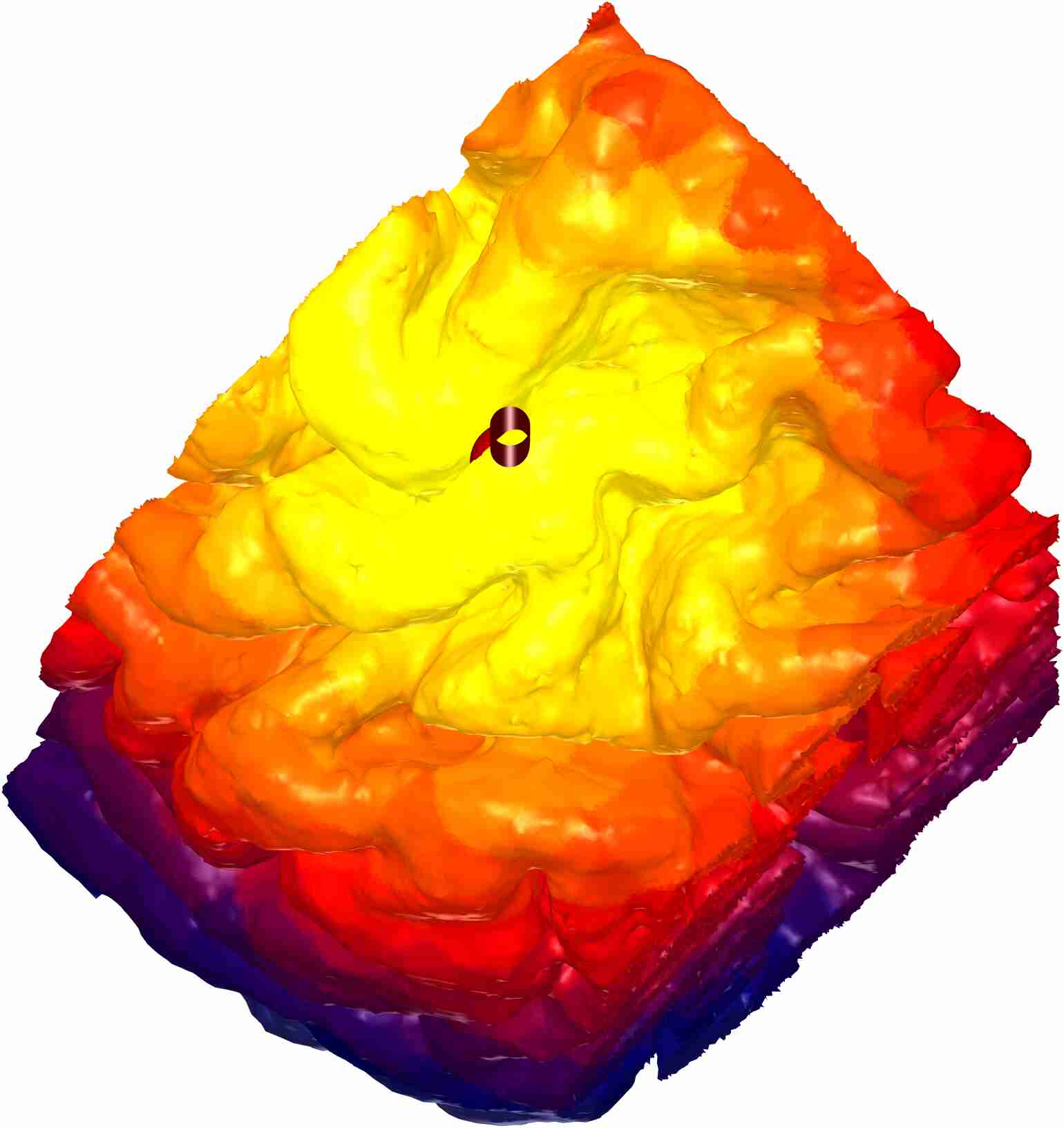}
		\end{subfigure}
		\hfill
		\begin{subfigure}[\subfigwidth]{}{%
			fig:normal-amplitude-single-source-dipolescan-d1mm%
		}
			\includegraphics[width=\linewidth]{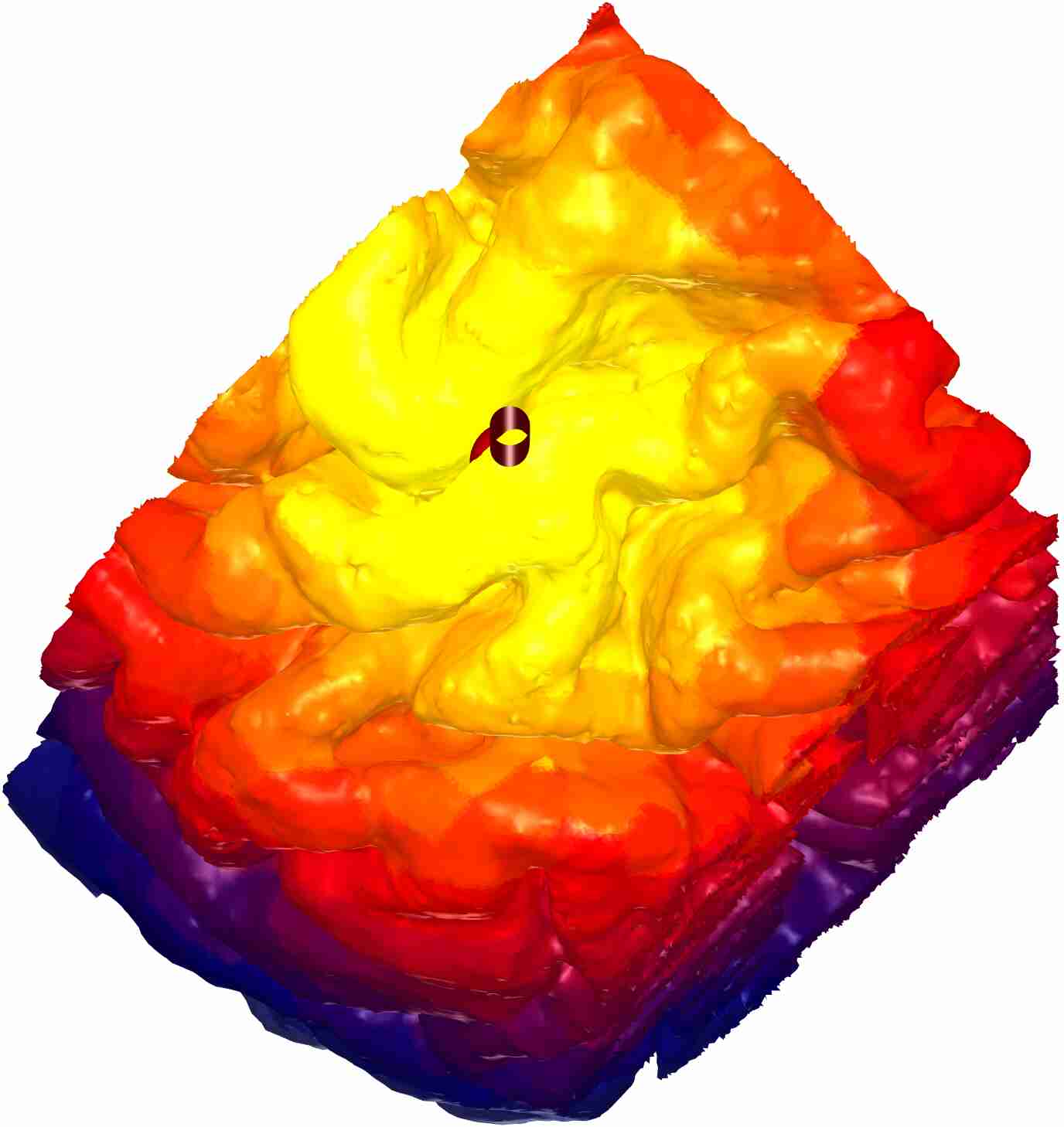}
		\end{subfigure}
	\end{minipage}
	\hfill
	\begin{minipage}[t]{0.09\columnwidth}%
		\centering%
		\vspace*{1cm}%
		\includegraphics[width=\linewidth]{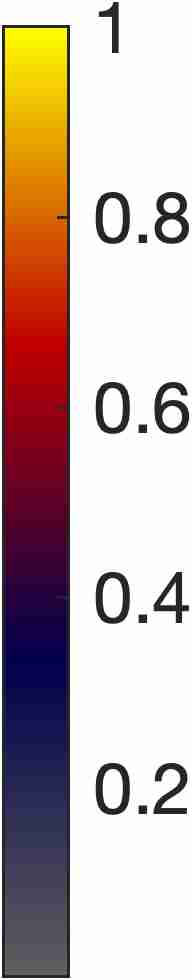}%
	\end{minipage}

	\caption{
		Relative strengths of \(\rec\Jp\) in the case of \DipoleScan, near
		Brodmann area 3b with \SI{30}{\decibel} noise.
		Subfigures~\subfigref{fig:tan-single-source-sloreta-d0mm}--\subfigref{fig:tan-single-source-sloreta-d1mm}
		display the case of a tangential source and
		Subfigures~\subfigref{fig:normal-single-source-sloreta-d0mm}--\subfigref{fig:normal-single-source-sloreta-d1mm}
		show normal sources, with \(d\subt p = \num{0}, \num{0.5}\) and
		\(\SI{1.0}{\milli\meter}\) respectively. The stronger the
		reconstruction located in a particular tetrahedron, the more yellow
		the tetrahedron is.%
	}

	\label{fig:tan-normal-amplitude-dipolescan-single-source}

\end{figure}

\begin{figure}[h!]

	\begin{minipage}[t]{0.8\columnwidth}
		\begin{subfigure}[\subfigwidth]{}{%
			fig:tan-single-source-sloreta-d0mm%
		}
			\includegraphics[width=\linewidth]{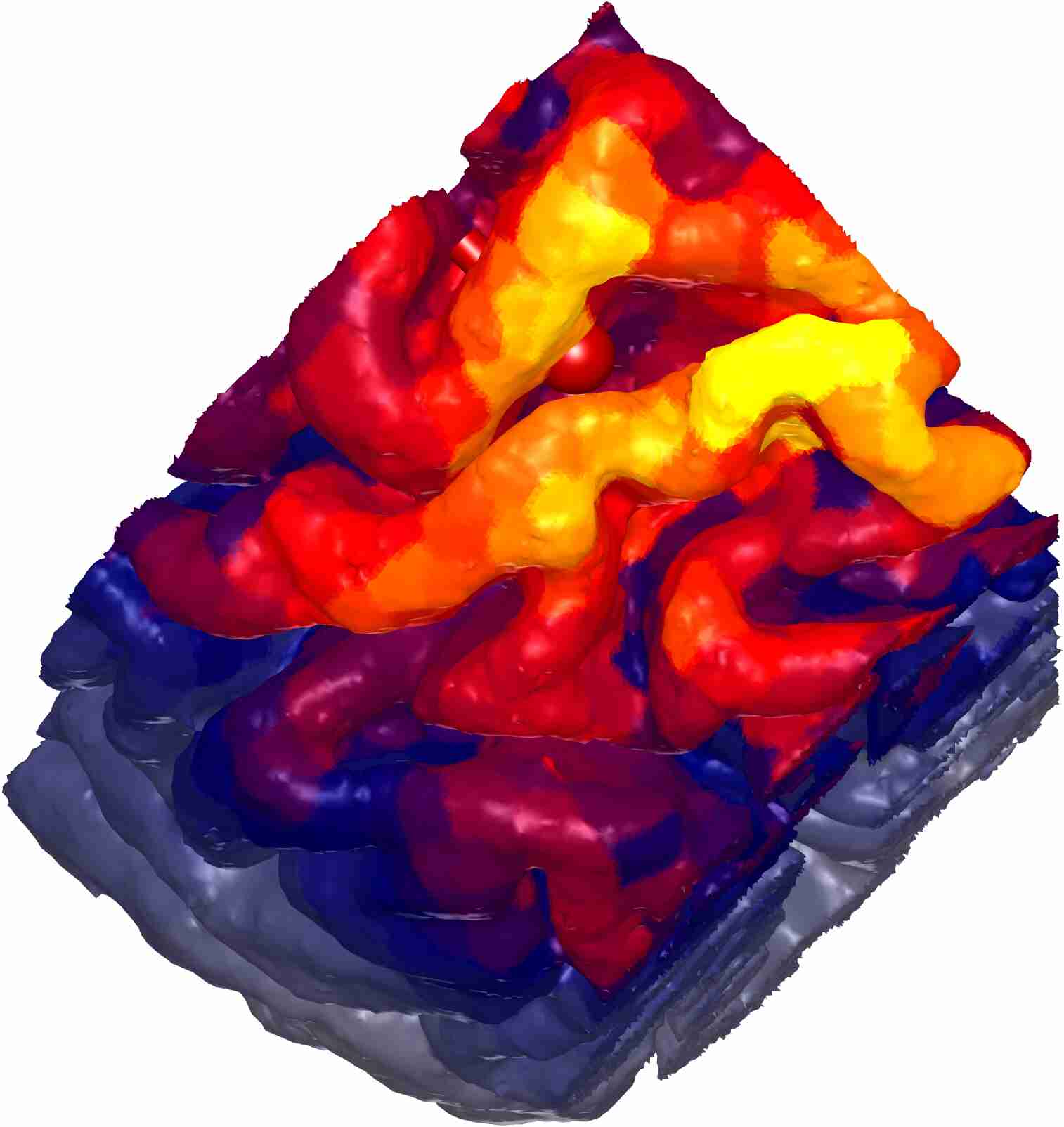}
		\end{subfigure}
		\hfill
		\begin{subfigure}[\subfigwidth]{}{%
			fig:tan-single-source-sloreta-d05mm%
		}
			\includegraphics[width=\linewidth]{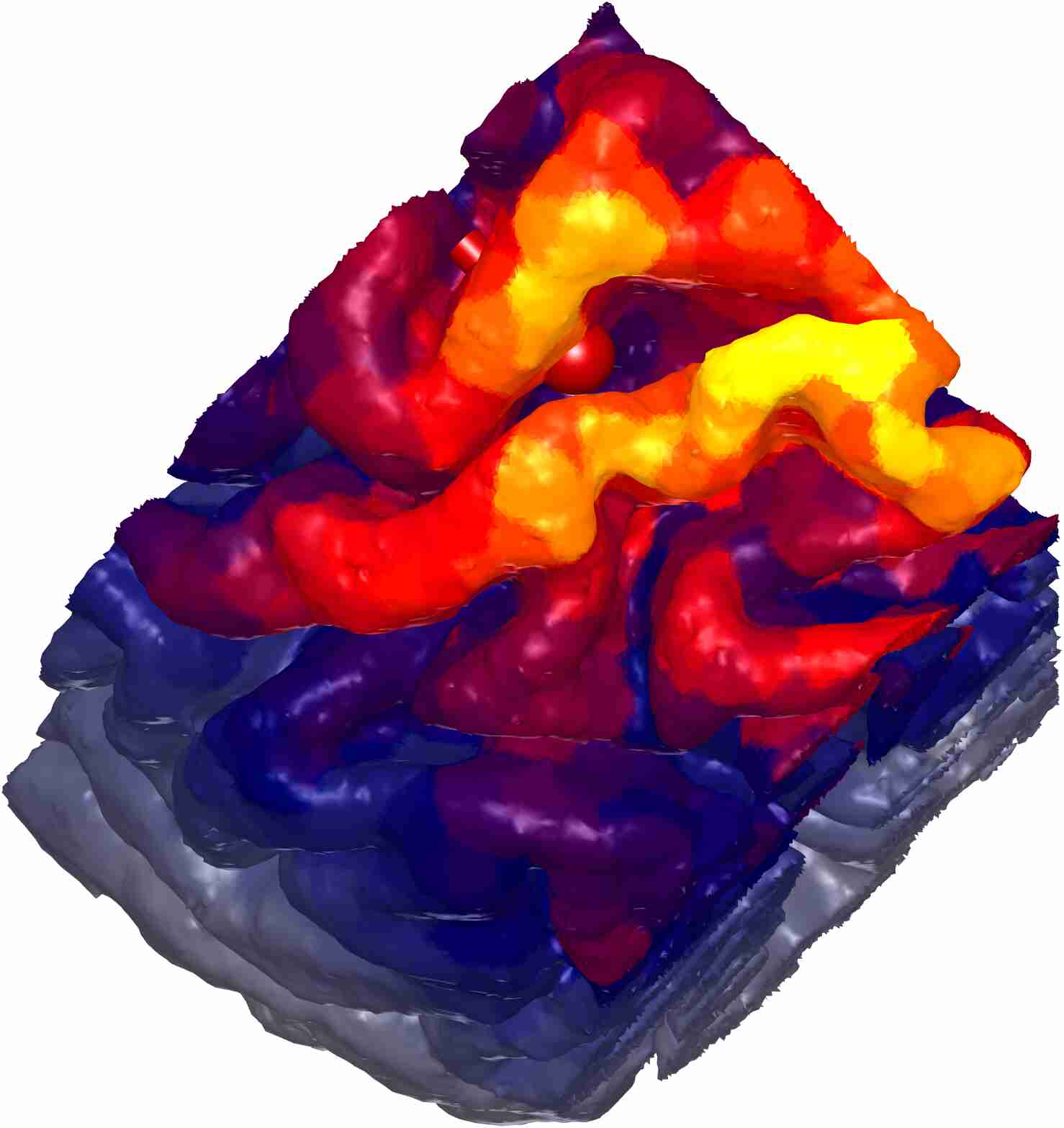}
		\end{subfigure}
		\hfill
		\begin{subfigure}[\subfigwidth]{}{%
			fig:tan-single-source-sloreta-d1mm%
		}
			\includegraphics[width=\linewidth]{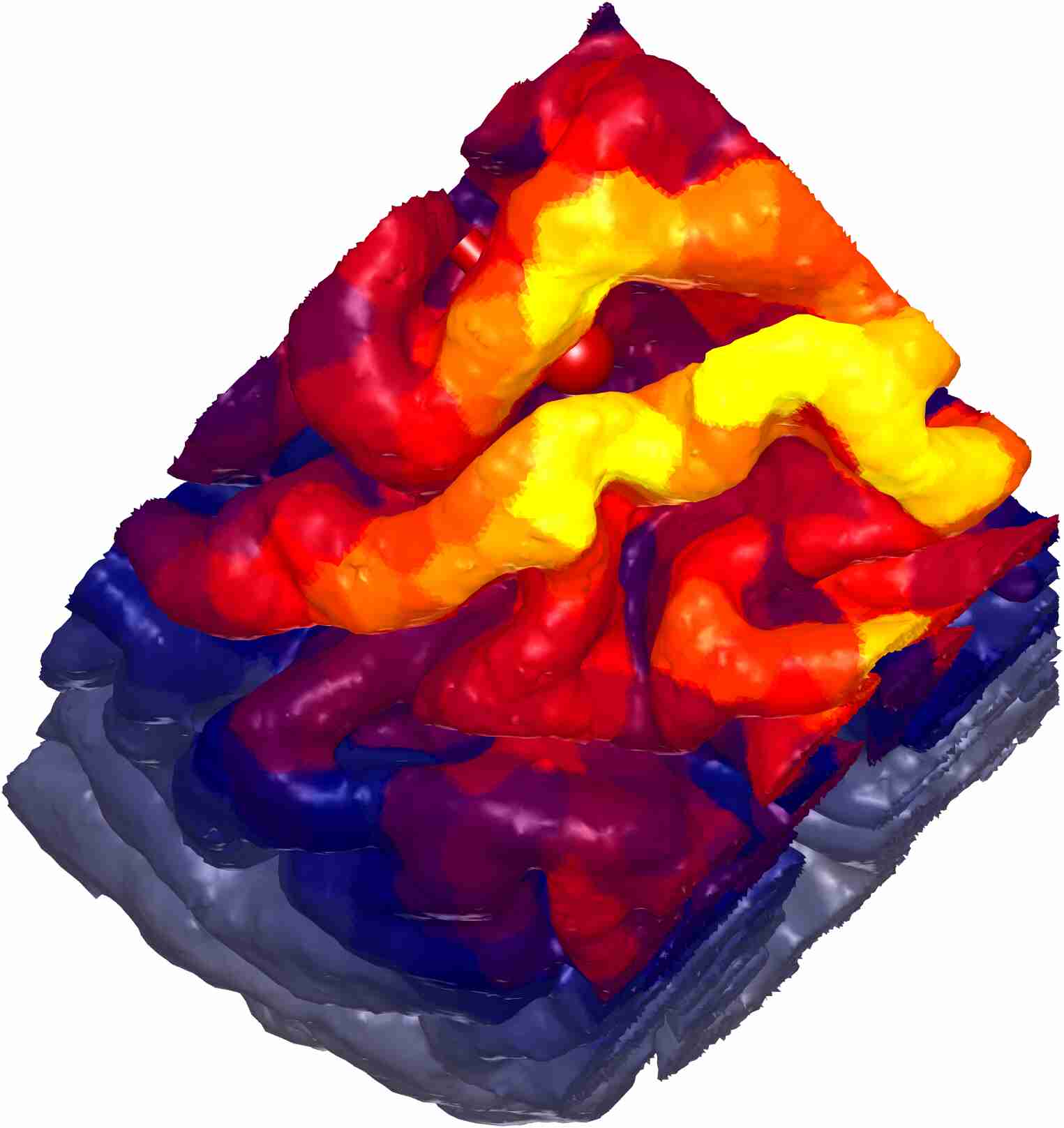}
		\end{subfigure}

		\bigskip

		\begin{subfigure}[\subfigwidth]{}{%
			fig:normal-single-source-sloreta-d0mm%
		}
			\includegraphics[width=\linewidth]{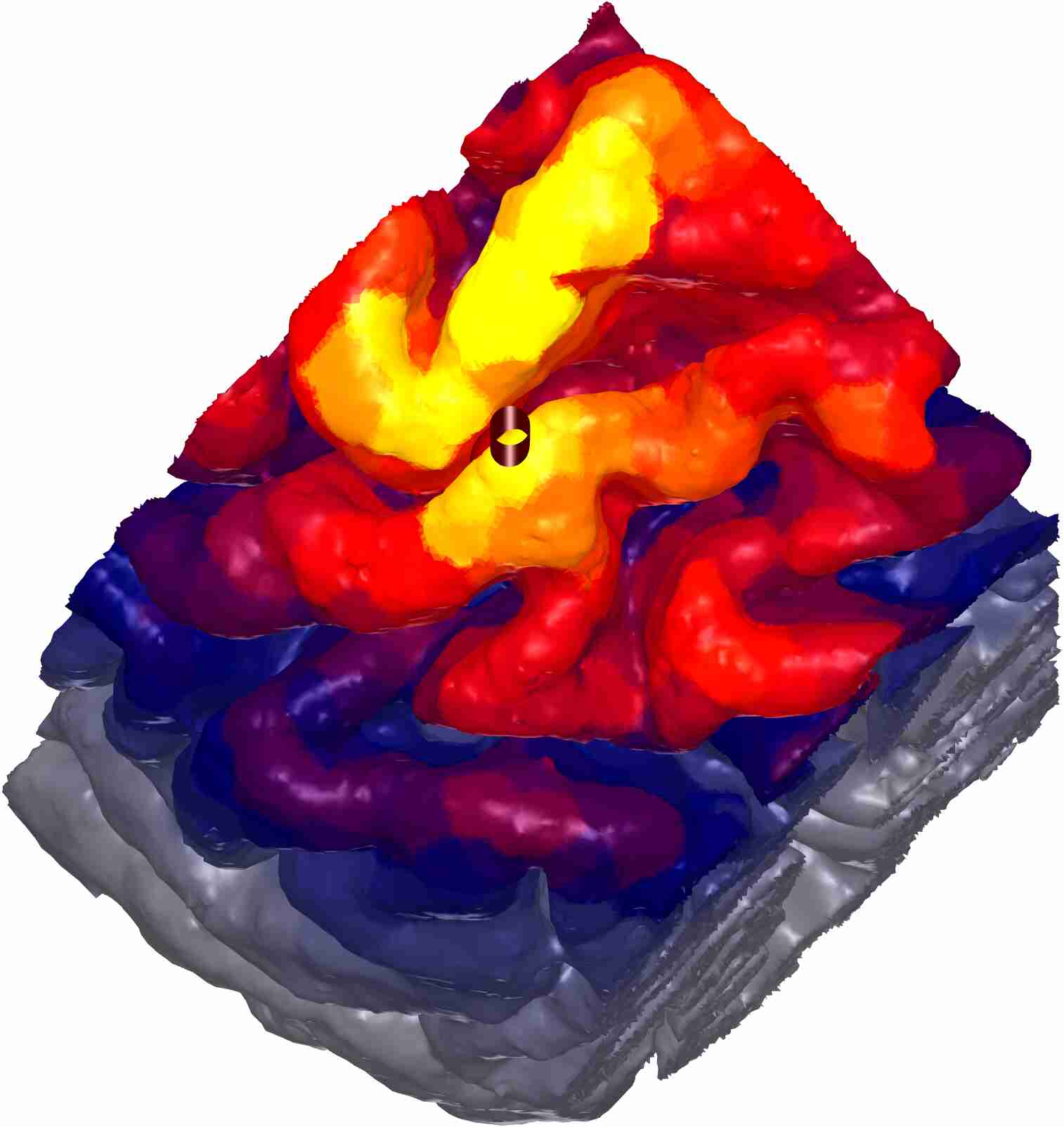}
		\end{subfigure}
		\hfill
		\begin{subfigure}[\subfigwidth]{}{%
			fig:normal-single-source-sloreta-d05mm%
		}
			\includegraphics[width=\linewidth]{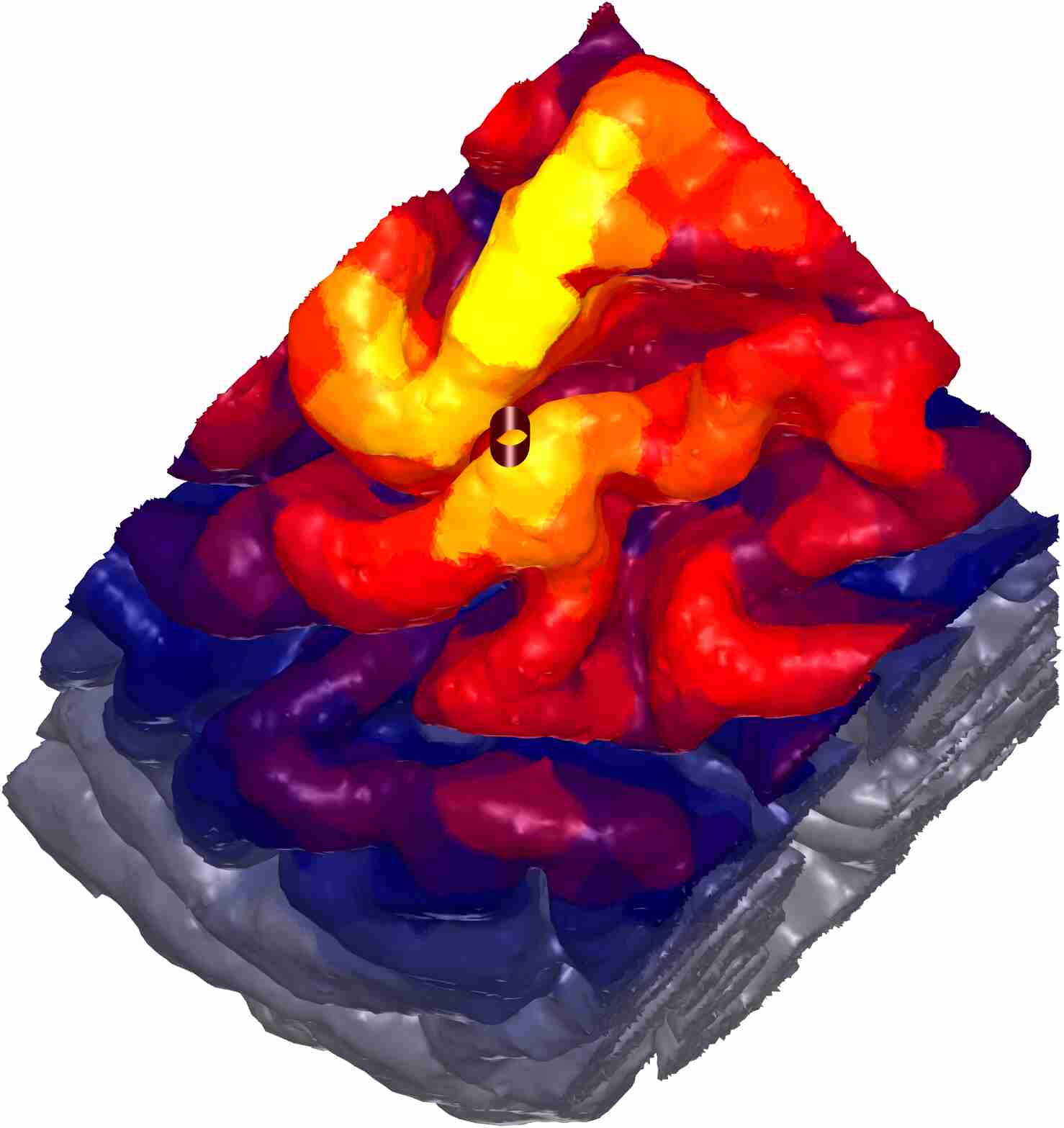}
		\end{subfigure}
		\hfill
		\begin{subfigure}[\subfigwidth]{}{%
			fig:normal-single-source-sloreta-d1mm%
		}
			\includegraphics[width=\linewidth]{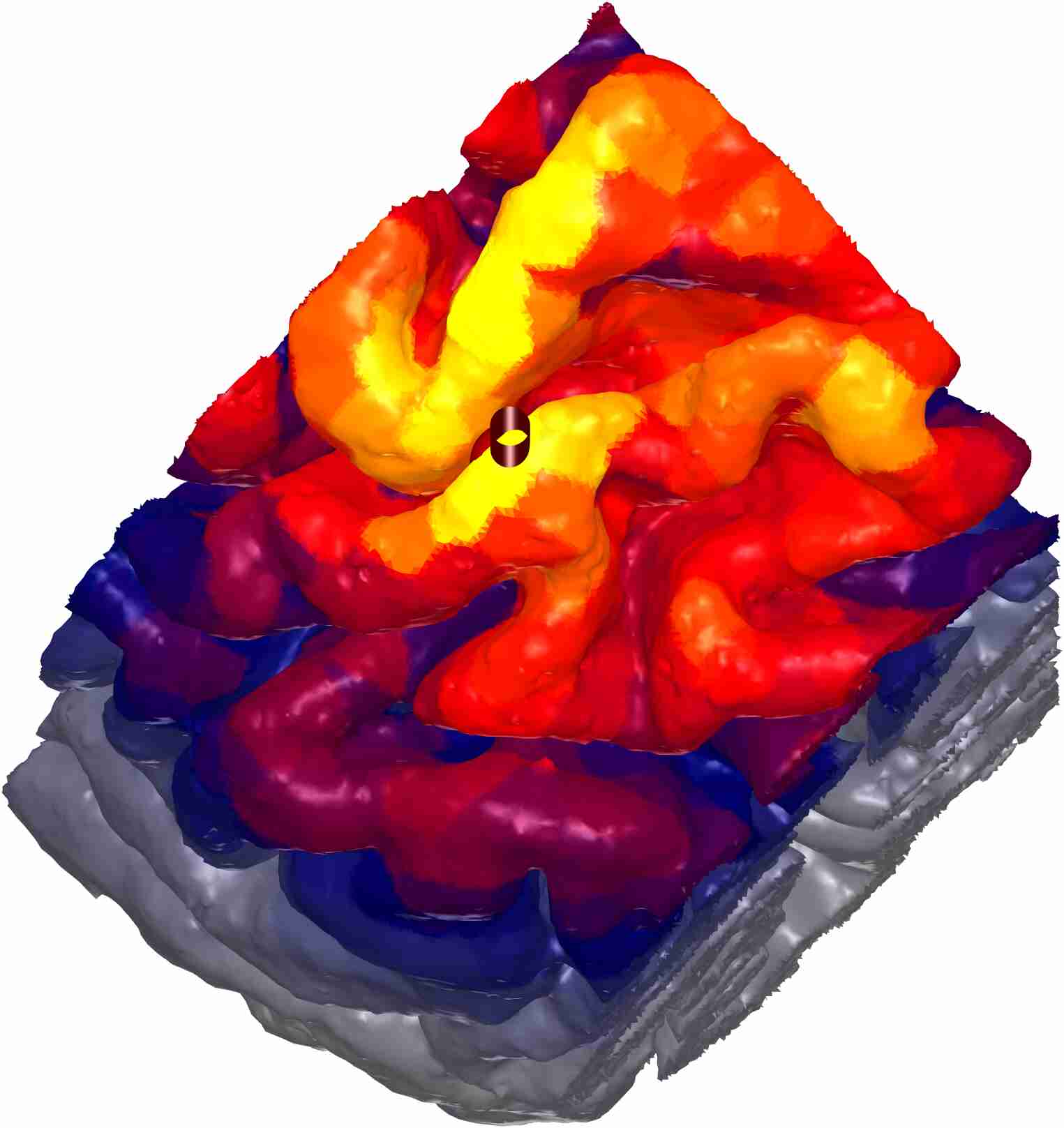}
		\end{subfigure}
	\end{minipage}
	\hfill
	\begin{minipage}[t]{0.09\columnwidth}%
		\centering%
		\vspace*{1cm}%
		\includegraphics[width=\linewidth]{images/single_source_colorbar.png.jpg}%
	\end{minipage}

	\caption{
		Relative strengths of \(\rec\Jp\) in the case of \sLORETA, near
		Brodmann area 3b with \SI{30}{\decibel} noise. The subfigures and
		colors are as in \figref{fig:tan-normal-amplitude-dipolescan-single-source}.
	}

	\label{fig:tan-normal-sloreta-single-source}

\end{figure}

In the case of \sLORETA, the region around the single reconstructed tangential
source increases in amplitude, as peeling depth \(\peeld\) is increased, which
can be seen as the moving of the bright yellow area from right to left,
towards the dipole position itself. The difference in the normal case is much
smaller. With \DipoleScan, the effects of \(\peeld\) on the inversion of the
tangential measurements are not obvious. The case of a radial or normal source
allows the differences to become more apparent: the area that contains the
strongest intesity shrinks, which makes the reconstruction more localised.

Finally, we observe the dispersion measures of equation~\eqref{eq:dispersion}
for \sLORETA\ and \DipoleScan\ at \SI{30}{\decibel} noise level. These are
presented in Figures~\ref{fig:dispersion-sloreta} and
\ref{fig:dispersion-dipolescan}, respectively. The ROI of dispersion around
each source position was set at \SI{30}{\milli\meter}.

\def\subfigwidth{0.32\linewidth}

\def\cbarwidth{0.05\linewidth}

\def\cbvspace{-0.5cm}

\begin{figure}[h!]
	\begin{minipage}[t]{0.9\linewidth}
		\begin{subfigure}[\subfigwidth]{\(\peeld = \SI{0}{\milli\meter}\).}{sfig:dispersion-sloreta-d0mm}
			\includegraphics[width=\linewidth]{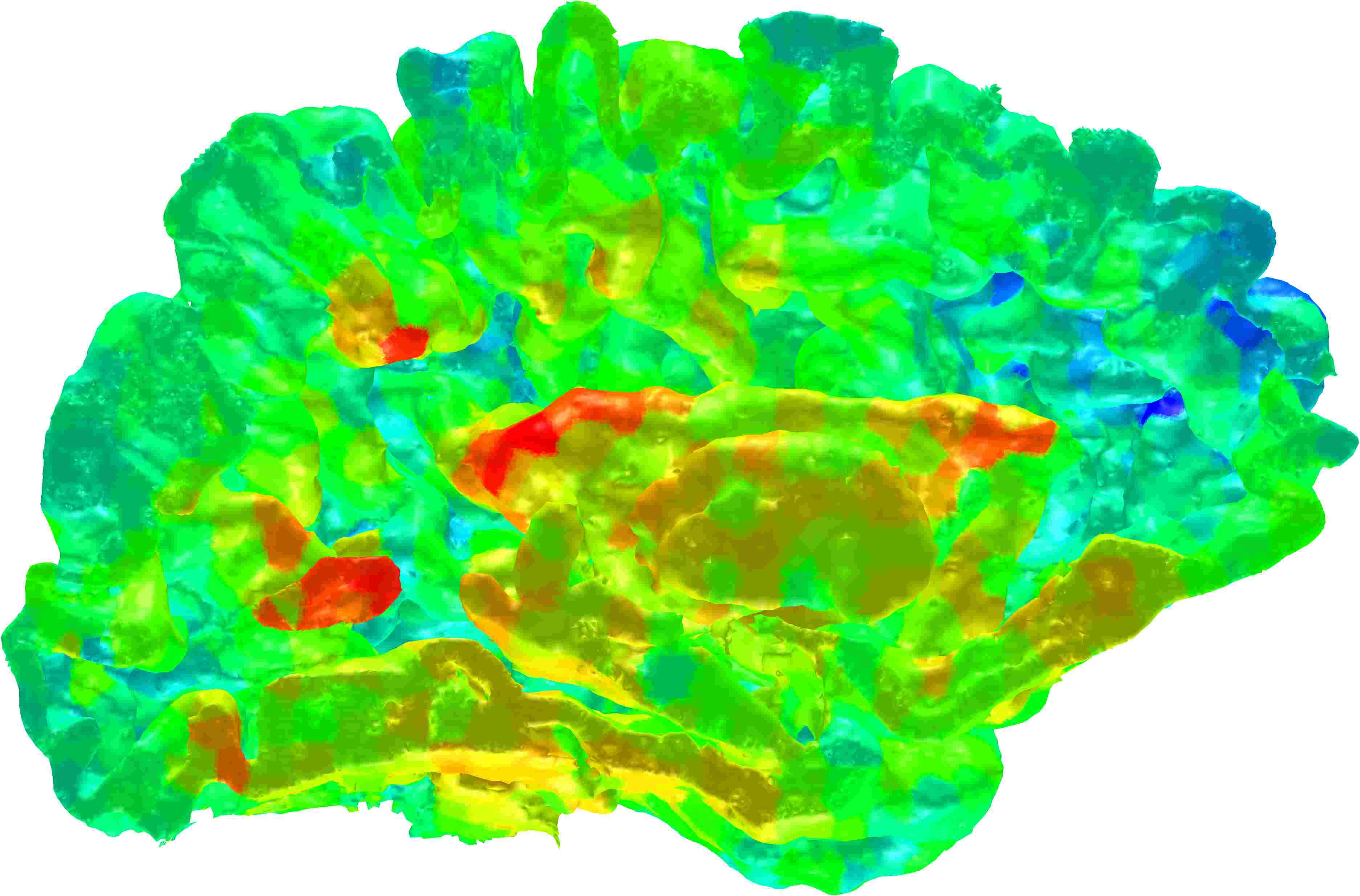}
		\end{subfigure}
		\hfill
		\begin{subfigure}[\subfigwidth]{\(\peeld = \SI{0.5}{\milli\meter}\).}{sfig:dispersion-sloreta-d0.5mm}
			\includegraphics[width=\linewidth]{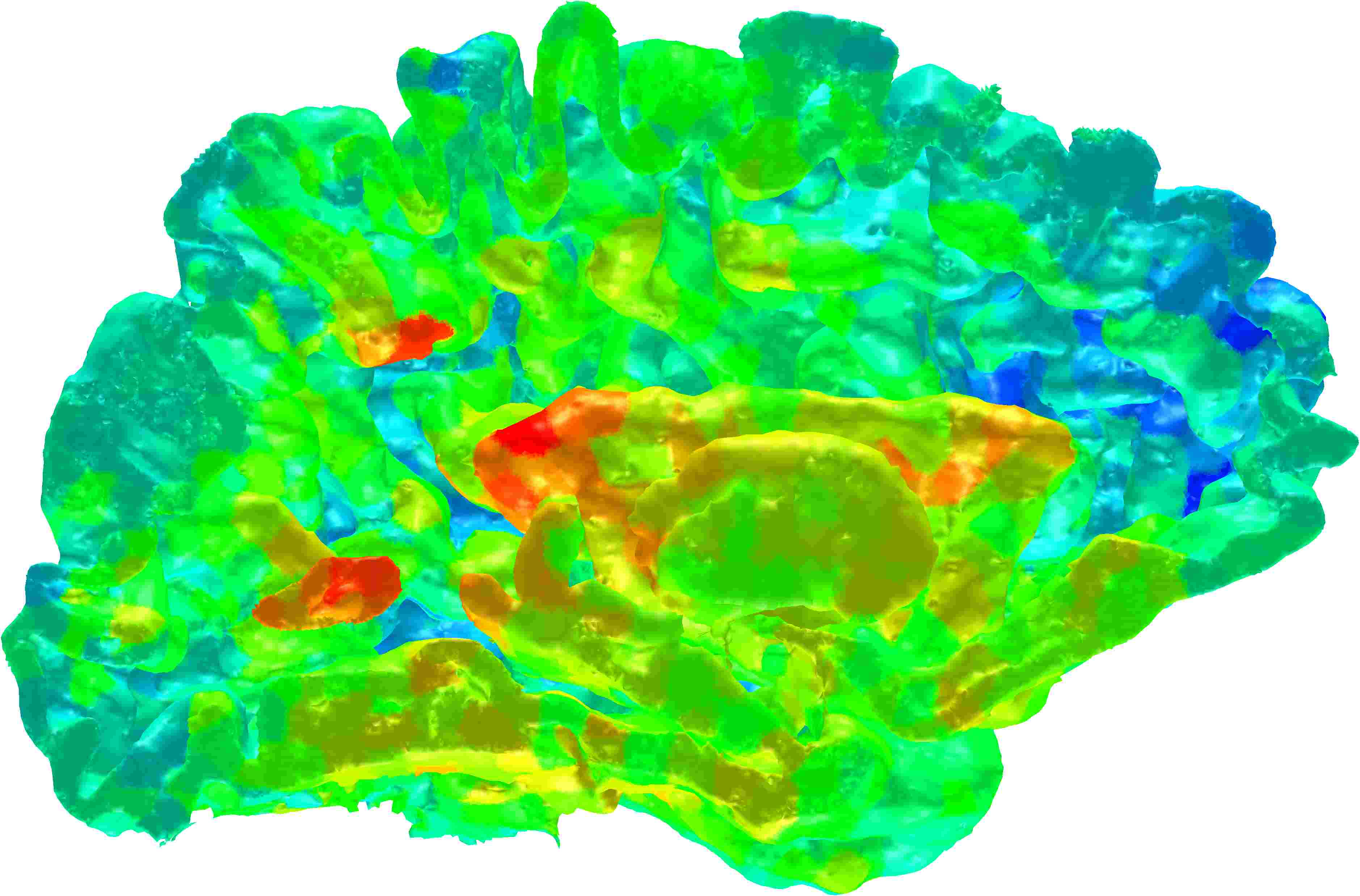}
		\end{subfigure}
		\hfill
		\begin{subfigure}[\subfigwidth]{\(\peeld = \SI{1}{\milli\meter}\).}{sfig:dispersion-sloreta-d1mm}
			\includegraphics[width=\linewidth]{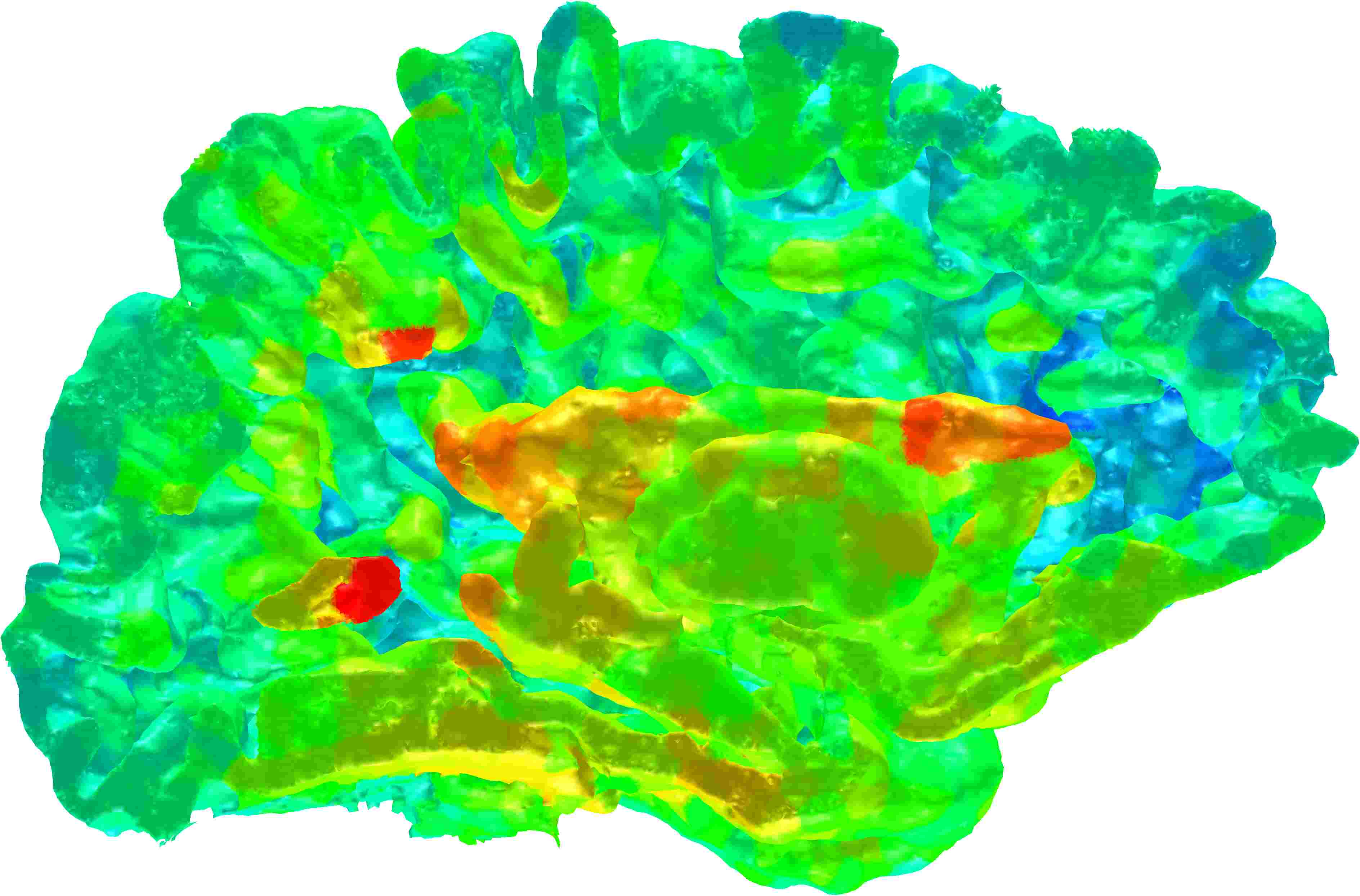}
		\end{subfigure}
	\end{minipage}
	\hfill
	\begin{minipage}[t]{\cbarwidth}%
		\centering%
		\vspace*{\cbvspace}%
		\includegraphics[width=\linewidth]{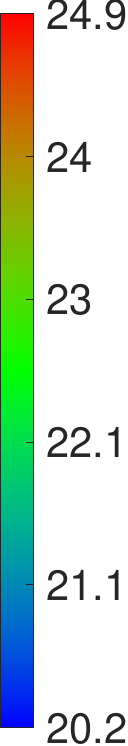}%
	\end{minipage}

	\caption{%
		Dispersion of equation \eqref{eq:dispersion} (in \unit{\milli\meter})
		of \sLORETA\ at \SI{30}{\decibel} noise level in each position in the
		active gray matter layer.
	}

	\label{fig:dispersion-sloreta}

\end{figure}

\begin{figure}[h!]
	\begin{minipage}[t]{0.9\linewidth}
		\begin{subfigure}[\subfigwidth]{\(\peeld = \SI{0}{\milli\meter}\).}{sfig:dispersion-dipolescan-d0mm}
			\includegraphics[width=\linewidth]{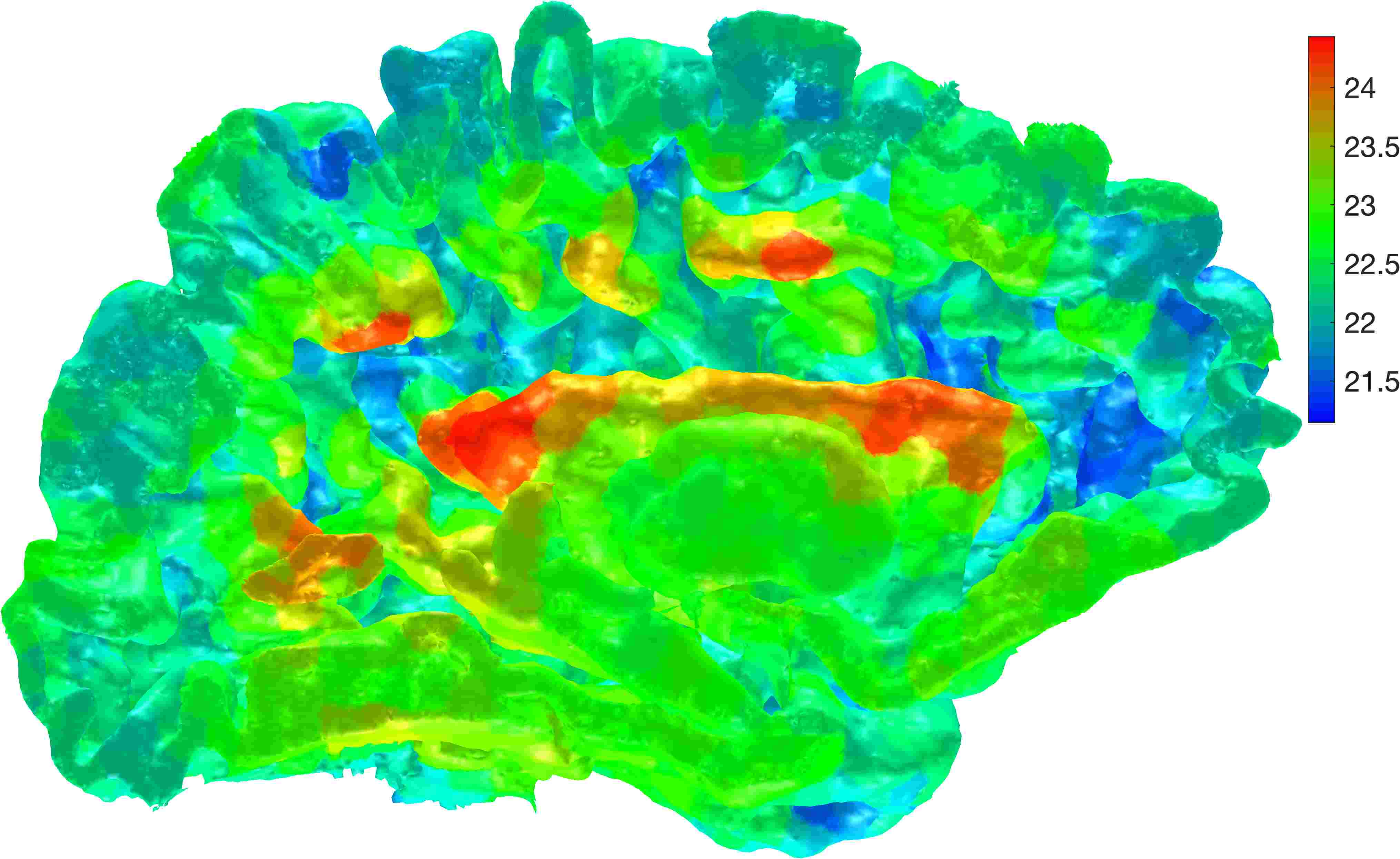}
		\end{subfigure}
		\hfill
		\begin{subfigure}[\subfigwidth]{\(\peeld = \SI{0.5}{\milli\meter}\).}{sfig:dispersion-dipolescan-d0.5mm}
			\includegraphics[width=\linewidth]{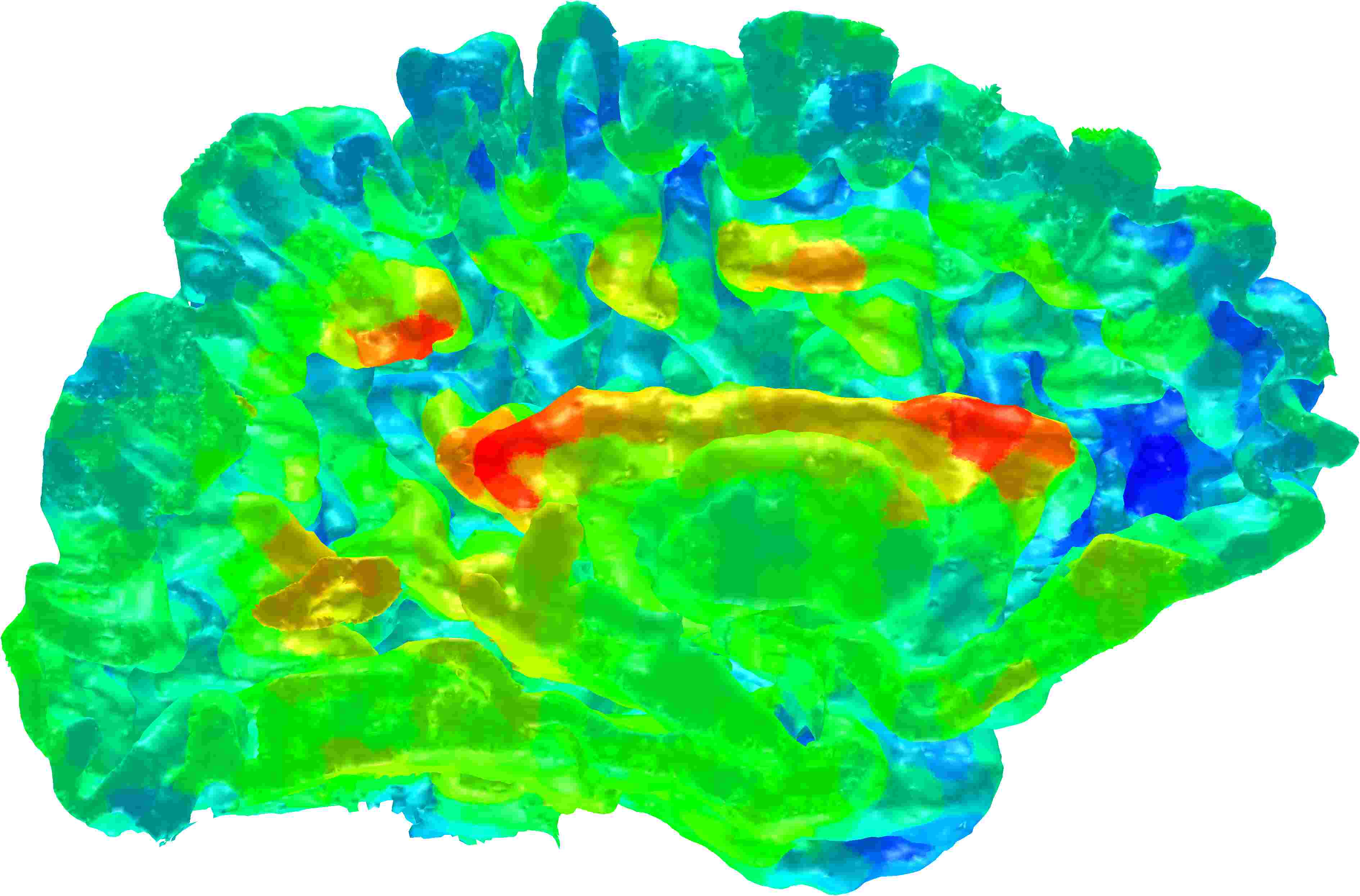}
		\end{subfigure}
		\hfill
		\begin{subfigure}[\subfigwidth]{\(\peeld = \SI{1}{\milli\meter}\).}{sfig:dispersion-dipolescan-d1mm}
			\includegraphics[width=\linewidth]{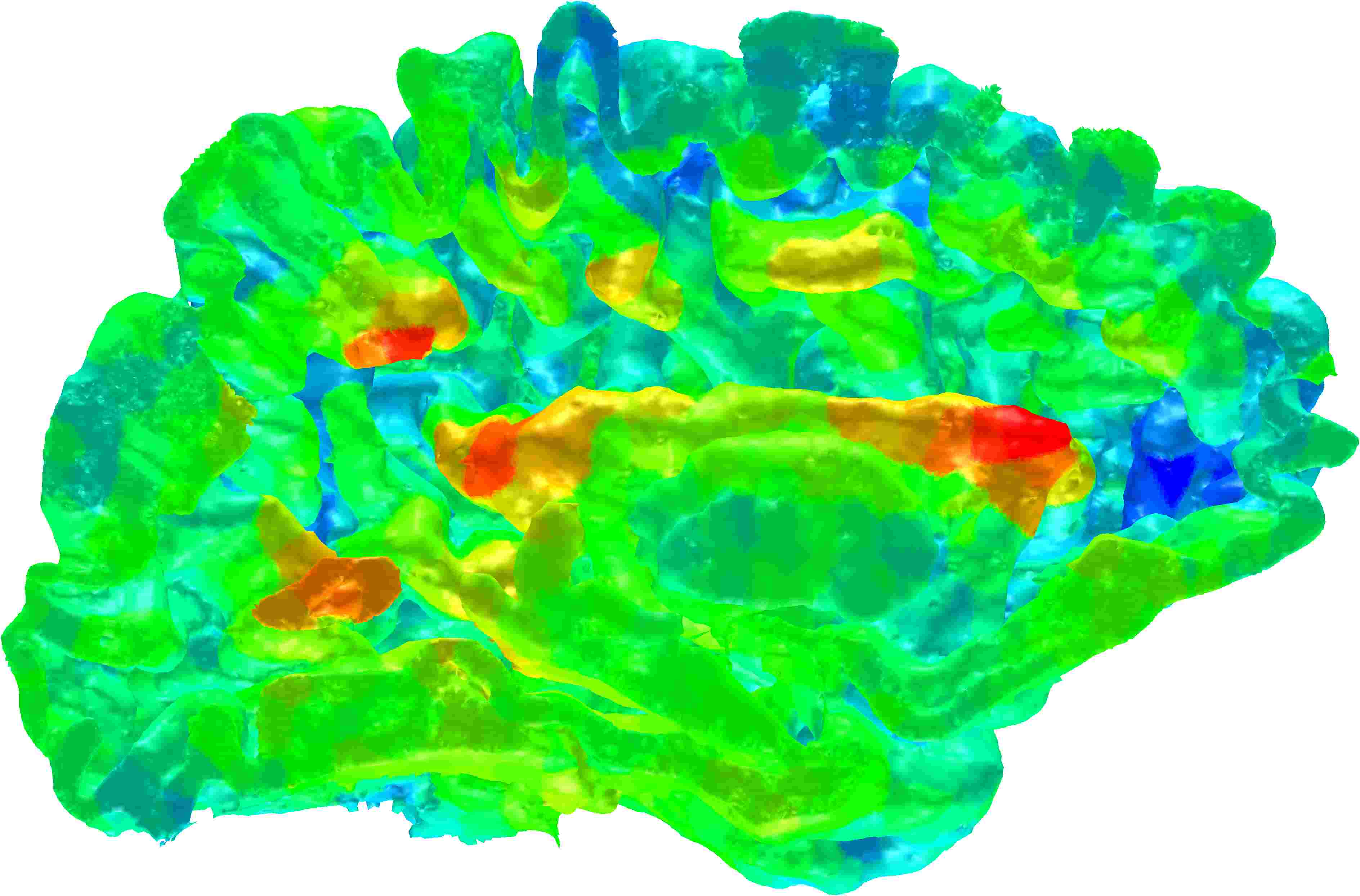}
		\end{subfigure}
	\end{minipage}
	\hfill
	\begin{minipage}[t]{\cbarwidth}%
		\centering%
		\vspace*{\cbvspace}%
		\includegraphics[width=\linewidth]{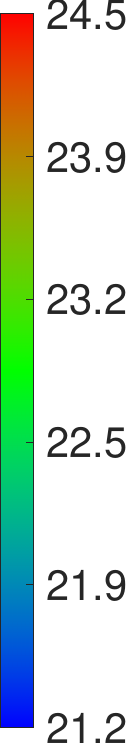}%
	\end{minipage}
	\caption{%
		Dispersion of equation \eqref{eq:dispersion} (in \unit{\milli\meter})
		of \DipoleScan\ at \SI{30}{\decibel} noise level in each position in
		the active gray matter layer.
	}
	\label{fig:dispersion-dipolescan}
\end{figure}

We see a clear reduction in the width of the peak around the reconstructed
dipoles, both around the thalamic region and in the cortex. For example the
dipersion at the back of the thalamic region decreases from roughly
\SI{24.5}{\milli\meter} to \SI{23.5}{\milli\meter} between
Subfigures~\subfigref{sfig:dispersion-sloreta-d0mm}--\subfigref{sfig:dispersion-sloreta-d1mm}
of \figref{fig:dispersion-sloreta}. A similar effect can be observed with
\DipoleScan\ in the subfigures of \figref{fig:dispersion-dipolescan}.

\section{Discussion}\label{sec:discussion}

\EEG\ source modelling is often performed on spherical head
models~\cite{ary-etal-1981}\cite{roth-1993}\cite{meneghini--etal-2010}, or
head models that are realistically shaped, but which only contain \num 3 or
\num 4 brain compartments, such as gray matter, skull, skin and maybe the
cerebrospinal fluid layers~\cite{vorwerk-2014}. In addition, instead of these
models being volumetric, only boundaries of the modelled brain compartments
might be included, if forward modelling is based on boundary element
methods~\cite{hamalainen-etal-1993}\cite{kybic--etal-2005}. Those are less
resource-intensive than finite element methods, but also do not model the
volumetric aspects of the problem as accurately.

In a realistic situation, there are at least 19 different brain tissue types,
and this is only if one does not take into account the inhomogeneities within
these major compartments, such as different types of bone in the
skull~\cite{dannhauer-2011}. Each compartment can have its own isotropic or
anisotropic conductivity structure, which affects the forward solution \(L\).
Simplified approaches thus have little hope of modelling the subtleties
introduced by a realistic head model~\cite{meneghini--etal-2010}. Therefore in
this paper, we performed our investigations in a realistically shaped
high-resolution finite element head model~\cite{duneuro-head-model}, with
\num{19} major brain compartments present. Most importantly, the modelled
brain activity was restricted to realistically shaped and volumetrically
refined compartments of gray matter, where especially the cortical portion is
only a few \si{\milli\meter} thick. This again required (as one option) the
use of the tightly supported \(\Hdiv\) source model, which allowed the
placement of dipoles with nonzero lengths into these thin gray matter layers.

In this setting, we investigated how restricting tetrahedra, into which a set
of \(\Hdiv\) dipoles \(\Jp = \Jp\arcs{\vec x}\) might be placed, affects the
source localisation error \(\Delta\) of \sLORETA\ and \DipoleScan. The
restriction was based on disallowing the placement of dipoles near surfaces of
active layers. The relevance of this study comes from the well-known fact,
that the convergence of the forward \EEG\ solution is negatively affected,
when sources are placed near conductivity discontinuities at said tissue
boundaries. This applies both to boundary element
methods~\cite{he-1999}\cite{fuchs-2001}, and different finite element
approaches, such as the Subtraction Method~\cite{wolters--etal-2007} and the
\(\Hdiv\) approach~\cite{pursiainen--vorwerk--wolters-2016}, which was used in
this paper. This again negatively affects the localisation of synthetic
sources via different inverse methods.

Before presenting the above results, we observed how the peeling algorithm
works. In addition to this, our numerical forward solver was compared to a
semi-analytical one, to ensure that it works appropriately in a simplified
spherical domain and at higher eccentricities, where cortical sources are
known to reside, but where the errors are also known to be more extensive.

The peeling algorithm, which restricts the positions \(\vec x\) where sources
can be placed, was found to function as intended: it recognized the intended
tetrahedra and always removed at least one layer from the surfaces of the
specified brain layers, to prevent singularities from appearing in the forward
solution because of discontinuities or jumps in conductivities between
neighbouring source tetra. The only anomaly related to the source positioning
after peeling was that the \(d\subt p = \SI{0.5}{\milli\meter}\) case seemed
to contain \emph{more} sources than the \(d\subt p = \SI{0}{\milli\meter}\)
case. This can be explained with \zi's iterative way of evenly distributing
sources into the active volume, if the initial guess is not near the
user-given amount.

Comparing the numerical lead field or forward solution \(L\) to that of the
semi-analytical \(3\)-layer Ary model~\cite{ary-etal-1981} produced
appropriately good results, even at above \SI{98}{\percent} eccentricities,
with upper \SI{75}{\percent} quantiles of magnitude measure~(\MAG) and relative
difference measure~(\RDM) being at most \num{0.06} and \num{0.04},
respectively. It was not obvious, whether a more refined mesh would diminish
differences between the semi-analytical and numerical solutions \(\La\) and
\(\Ln\). In fact, having a more refined mesh seems to produce deteriorated
median results in the case of \MAG, whereas with \RDM\ the outcome is slightly
improved. This was especially true when \PBO~\cite{bauer--etal-2015} was used
for synthetic source interpolation, and hence it was chosen as the optimization
method in the case of the realistic head model. \MPO~\cite{atena-2021} was also
considered, but was deemed less suitable, as according
to~\cite{pursiainen--vorwerk--wolters-2016} \MPO\ is less stable than \PBO.

When looking at the average localisation errors \(\Delta\), and comparing the
results to those of Cuffin et. al.~\cite{cuffin-2001}, one can observe that
the performance of \sLORETA\ is inferior to that of \DipoleScan. The reported
\(\Delta < \SI{10.5}{\milli\meter}\) is only reached at lower noise levels,
when \(\SNR\geq\SI{30}{\decibel}\), whereas with \DipoleScan\ an average
within acceptable boundaries can be obtained all the way down to \(\SNR =
\SI{15}{\decibel}\). In the case of \sLORETA, the peeling cannot be said to
improve the average localisation accuracy, as utilizing it seems to first
improve the result when going from a peeling depth \(\peeld =
\SI{0}{\milli\meter}\) to \(\peeld = \SI{0.5}{\milli\meter}\), but then
localisation accuracy is decreased again, when \(\peeld\) is increased to
\SI{1}{\milli\meter}. Only at the very highest \SNR\ can a slight but
systematic decrease in the average \(\Delta\) be seen. \DipoleScan\ is
slightly more consistent in its performance with respect to peeling, as with
\(\SNR\geq\SI{20}{\decibel}\) the added peeling no longer increases the
average \(\Delta\). It should be noted, that Cuffin used a 3-compartment
boundary element model with a simplex search method to locate
sources~\cite{press-1989}, and hence our results are not directly comparable.

The behaviour of localisation error outliers also reflect the above state of
matters. In the case of \sLORETA, the initial impression is that peeling seems
to reduce the general denseness of the outliers clouds further above the box
plots, but at the same time some maximum outliers are further away from the
tops of the whiskers, especially in the case of \(d\subt p =
\SI{0.5}{\milli\meter}\). Peeling somewhat increases the number of outliers
for the inverse method at higher \SNR, although with \(\SNR <
\SI{20}{\decibel}\) there seem to be cases where the outliers are reduced with
added peeling. A possible explanation for this is that peeling reduces
statistical variability in the low-\SNR\ cases of \sLORETA. \DipoleScan\ again
performs more consistently, with the numbers of outliers being reduced at
every \(\SNR\geq\SI{15}{\decibel}\), when the peeling depth \(d\subt p\) is
increased. It might then be said that for \DipoleScan, the effects of random
noise seem to dominate when \SNR\ is reduced beyond this point.

When observing the localisation error as a function of primary dipole
position, \(\Delta = \Delta\arcs{\vec x}\) of \sLORETA\ and \DipoleScan, the
superior performance of the latter becomes obvious. With \sLORETA, while the
actual scale of the localisation error is reduced as \SNR\ increases, we end
up having more outliers with a relatively large \(\Delta\) in the entire
volume. In contrast, the behaviour of \DipoleScan\ is again more consistent:
not only does the localisation error go down with increasing \SNR, as one
would expect, the localisation error outliers become more focused into the
deep structures of the brain, where the sources are further away from the
sensors. The benefits of increasing the peeling depth \(\peeld\) are also not
unambiguous in either case. It seems that especially in the case of \sLORETA,
increasing \(\peeld\) simply moves the positions where outliers occur around
the volume, instead of eliminating them. There are still places where the
areas containing larger errors become less focal, such as in the parietal
region at \(\SNR = \SI{25}{\decibel}\). The fluctuation of the higher
localisation error patches around the volume is also visible in the case of
\DipoleScan, although maybe to a lesser extent. The consistent reduction in
\(\Delta\) with increased peeling seems to be slightly more prominent, as can
be seen in the environment of \emph{corpus callosum} and around the lower part
of the frontal lobe at \SNR\ between \num{20}--\SI{30}{\decibel}. The takeaway
from this is, that peeling the active brain layers does not guarantee
unambiguously better localisation results, but that it might do so on a
regional and \SNR\ basis.

A look at the relative strengths of reconstructions at
\(\SNR=\SI{30}{\decibel}\) also mostly supports this interpretation from the
point of view of varying the peeling depth \(\peeld\), although there seems to
be a hope of seeing an actual improvement in the case of \DipoleScan\ and a
normal source, which is the ideal case for \EEG: we can see the reconstruction
becoming slightly more focal around the plotted source. With \sLORETA, the
peeling seems to increase the general strength of the reconstruction around
the tangential source, whereas in the case of a normal source we see a
\emph{very} slight shift of the local peak of the distribution \(\rec\Jp\)
towards the source itself.

The dispersion measures of the results of different inverse methods were seen
to be improved, with parts of regions of largest dispersions seeing a
reduction of \SI{1.0}{\milli\meter} in the width of the peak around a
reconstructed dipole position. This is yet another indication of the
reconstruction becoming more focal in those regions, as a reduction in
dispersion indicates, that the amplitude of \(\rec\Jp\) is reduced in the
\SI{30}{\milli\meter} ROI for the dispersion around each respective source
position, with the most focal point being at the center of the ROI.

\section{Conclusions}\label{sec:summary}

\emph{Summa summarum}, while restricting the synthetic dipoles \(\Jp\), that
model brain activity to be further away from active brain compartment
boundaries, and therefore from conductivity jumps, the reconstructed
distribution of dipoles \(\rec\Jp\) is altered, but not definitely improved.
The peeling algorithm responsible for the application of the above restriction
was enabled by a local refinement or resolution increase of the \fe\ mesh, and
was shown to uniformly reduce statistical outliers in the forward solution
\(L\), which maps \(\Jp\) to the observed potentials at the \EEG\ electrodes,
thereby increasing its robustness.

This reduction in outliers is best observed with noise levels below
\SI{30}{\decibel}, which occur in medical studies, where possibly thousands of
stimuli might be applied over the duration of the experiment, such as when
measuring somatosensory-evoked potential responses~\cite{cruccu-2008}. In
addition to improving the robustness of the forward solution, peeling had the
same effect on the reconstruction, in addition to enhancing the regularity of
it. The importance of the regularity of a reconstruction comes from the need
to locate multiple or distributed sources, such as in the case of localising
epileptogenic zones or their irritative regions for the purposes of
treatment~\cite{van-mierlo-2020}\cite{fernandez-corazza-2021}.

\section{Future prospects}\label{sec:future-prospects}

Follow-up studies might include further increasing the base resolution and/or
the refinement of the active layers in the mesh, so that the possible
locations of dipoles would become more varied. Peeling depths might also be
observed in a wider or more refined range than was presented in this study, to
find out whether there is some optimal value for it. This optimal value might
depend on the mesh geometry, and therefore the tests might be performed on
multiple different head models with the same mesh resolution.

For case studies, the use of the peeling technique might be applied to head
models of patients with conditions such as
\emph{microlissencephaly}~\cite{abdel-razek-2009}, where the cortex is
flattened and possibly thickened. A
\emph{schizencephalic}~\cite{abdel-razek-2009}\cite{loturco-2013} brain model,
where a part of the brain has been displaced, and the displaced sections are
connected only by a thin ribbon of gray matter is another possibly interesting
target of peeling.

\section*{Acknowledgements}\label{sec:acknowledgements}

\def\AoF{Academy of Finland}
\def\CoE{Centre of Excellence}
\def\IMaI{Inverse Modelling and Imaging}
\def\PerEpi{Personalised diagnosis and treatment for refractory focal paediatric and adult epilepsy}
\def\DAAD{German Academic Exchange Service}
\def\Vaisala{Vilho, Yrjö and Kalle Väisälä Fund}

\def\aAoF{AoF}
\def\aCoE{CoE}
\def\aIMaI{IMaI}
\def\aPerEpi{PerEpi}
\def\aDAAD{DAAD}
\def\aVaisala{VYKVF}

\NewDocumentCommand\bfn{m}{\textbf{(#1)}}

This work has been supported by \bfn 1 the \AoF~(\aAoF) \CoE~(\aCoE) in \IMaI\ 2018-2025, \bfn 2 by the \aAoF~project \#344712 and the Bundesministerium für Gesundheit (BMG) project ZMI1-2521FSB006, under the frame of ERA PerMed as project ERAPERMED2020-227 (\aPerEpi), \bfn 3 by the bilateral \aAoF\ and \aDAAD\ projects  \#354976 and \#57663920, and \bfn 4 by \Vaisala~(\aVaisala). Per author, the funding sources were as follows: Santtu Söderholm~(\aAoF/PerEpi), Joonas Lahtinen~(\aVaisala), Carsten H. Wolters~(\aDAAD/PerEpi)  and Sampsa Pursiainen~(\aCoE, \aAoF/PerEpi).

\section*{Declaration of competing interests}

The authors declare that they have no known competing financial interests or
personal relationships, that could have appeared to influence the work
reported in this paper.

\appendix

\bibliographystyle{elsarticle-num}

{\interlinepenalty=10000\raggedright\bibliography{references}}

\end{document}